\documentclass[11pt]{article}
\usepackage{amsmath,amsfonts,amsthm,epsfig,latexsym}
\usepackage[english,francais]{babel}
 
 \def\NN{{\mathbb N}}  
 \def\RR{{\mathbb R}} \def\SS{{\mathbb S}} 
   
 \def\ZZ{{\mathbb Z}}

\def\Si{\Sigma}
\def\La{\Lambda}
\def\De{\Delta}
\def\Om{\Omega}
\def\Ga{\Gamma}

\def\cA{{\cal A}}  \def\cG{{\cal G}}  \def\cS{{\cal S}}
\def\cB{{\cal B}}    
\def\cC{{\cal C}}   \def\cO{{\cal O}} \def\cU{{\cal U}}
\def\cD{{\cal D}}   \def\cP{{\cal P}} \def\cV{{\cal V}}
\def\cE{{\cal E}}    
   \def\cR{{\cal R}}

\def\vdashv{\vdash\!\dashv}
\def\diff{\operatorname{Diff}}

\newtheorem{theo}{Th\'eor\`eme}[section]
\newtheorem{theor}{Th\'eor\`eme}
\newtheorem{affi}{Affirmation}
\newtheorem{lemme}{Lemme}
\newtheorem{lemm}{Lemme}[section]
\newtheorem{coro}[lemm]{Corollaire}
\newtheorem{defi}[lemm]{D\'efinition}
\newtheorem{prop}[lemm]{Proposition}
\newtheorem{conj}[lemm]{Conjecture}
\newtheorem{prob}{Probl\`eme}

\newtheorem{ques}{Question}

\newtheorem{rema}[lemm]{Remarque}

\newenvironment{demo}[1][D\'emonstration]{\noindent {\bf #1~: }}{\hfill$\Box$\medskip}
\title{R\'ecurrence et g\'en\'ericit\'e}
\author{Christian Bonatti et Sylvain Crovisier}
\date{\today}

\begin{document}
\maketitle

\begin{abstract}
Nous montrons un lemme de  connexion $C^1$ pour les pseudo-orbites des diff\'eomorphismes des vari\'et\'es compactes. Nous explorons alors les cons\'equences pour les diff\'eomorphismes $C^1$-g\'en\'eriques. Par exemple, les diff\'eomorphismes conservatifs $C^1$-g\'en\'eriques sont transitifs.
\end{abstract}
{\selectlanguage{english}
\begin{abstract}{We prove  a $C^1$-connecting lemma for pseudo-orbits of diffeomorphisms on compact manifolds. We explore some consequences for $C^1$-generic diffeomorphisms. For instance, $C^1$-generic conservative diffeomorphisms are transitive\footnote{An announcement of these results, written in english, \cite{BC}, is available at the adresse \\
{\bf \tt http://math.u-bourgogne.fr/topolog/bonatti/preprints/chaincras.ps}}.

}\end{abstract}}

\tableofcontents

\section{Introduction}
\subsection{R\'ecurrence des diff\'eomorphismes g\'en\'eriques}
\subsubsection{Les divers types de r\'ecurrences} 
On sait depuis le d\'ebut du 20$^{\grave eme}$ si\`ecle que  la dynamique globale des diff\'eomorphismes ou des \'equations diff\'erentielles peut \^etre tr\`es complexe m\^eme en basse dimension~: pour d\'ecrire cette dynamique globale il ne suffit pas de chercher des \'etats d'\'equilibre (des orbites p\'eriodiques) qui seraient la limite de toute orbite. L'existence d'intersections homoclines transverses, d\'ecouverte par H. Poincar\'e \cite{Po}, est un ph\'enom\`ene simple et r\'esistant aux perturbations   et G. D. Birkhoff  \cite{Bi2} a montr\'e qu'il entra\^\i ne l'existence d'une infinit\'e d'orbites p\'eriodiques dont la p\'eriode tend vers l'infini. 

Cette complexit\'e  des syst\`emes a amen\'e les dynamiciens \`a chercher la meilleure notion  d'ensemble qui concentrerait la dynamique non-triviale. De nombreuses d\'efinitions ont \'et\'e uti\-li\-s\'ees. En voici quelques unes parmi les plus fr\'equentes, exprim\'ees ici pour un diff\'eomorphisme $f$ d'une vari\'et\'e compacte $M$~: 
\begin{itemize}
\item l'ensemble ${Per}(f)$ des points p\'eriodiques,
\item les ensembles des points positivement ou n\'egativement r\'ecurrents de $f$, que l'on note $Rec_+(f)$ et $Rec_-(f)$ ,
\item les ensembles $L_+(f)$ et $L_-(f)$ qui sont l'union des ensembles $\omega$-limite et $\alpha$-limite des points de la vari\'et\'e, et leur union  $L(f)= L_+(f)\cup L_-(f)$,
\item l'ensemble $\Om(f)$ des {\em points non-errants} de $f$ (appel\'e \'egalement {\em ensemble non-errant}), introduit par Birkhoff dans \cite{Bi1}~: un point est dit {\em errant} s'il admet un voisinage disjoint de tous ses it\'er\'es. Cet ensemble est compact, invariant, et contient l'adh\'erence de tous les ensembles cit\'es ci-dessus. 
\end{itemize}

En g\'en\'eral ces ensembles ne co\"\i ncident pas : pour deux quelconques d'entre eux il existe des diff\'eomorphismes pour lesquels ils sont diff\'erents (voir \cite{Sh} pour une exposition d\'etaill\'ee sur les diff\'erents types de r\'ecurrence).

Plus r\'ecemment C. Conley \cite{Co}  et R. Bowen \cite{Bo} ont introduit les notions de {\em pseudo-orbites} et  de {\em points r\'ecurrents par cha\^\i nes}~: ce sont les points qui admettent des pseudo-orbites ferm\'ees, c'est-\`a-dire que l'on peut fermer leur  orbite si l'on admet une erreur abitrairement petite faite \`a chaque it\'eration.  Cette notion provient d'une vision perturbative de la dynamique~: les pseudo-orbites sont les orbites que l'on voit si l'on ajoute un bruit \`a la dynamique. 

La diff\'erence entre {\em l'ensemble  r\'ecurrent par cha\^\i nes}  $\cR(f)$ (i.e. l'ensemble des points r\'e\-cur\-rents par cha\^\i nes) et l'ensemble non-errant se comprend bien dans l'exemple des dynamiques hyperboliques, i.e. v\'erifiant l'Axiome A. Dans ce cas le th\'eor\`eme de d\'ecomposition spectrale  de Smale (voir \cite{Sm}) d\'ecompose l'ensemble non-errant $\Om(f)$ en union finie d'ensembles hyperboliques transitifs $\La_i$ disjoints, appel\'es {\em pi\`eces basiques}. On appelle alors {\em cycle} toute suite p\'eriodique de pi\`eces basiques telle que la vari\'et\'e instable de chacune rencontre la vari\'et\'e stable de la suivante. Les points errants mais r\'ecurrents  par cha\^\i nes sont alors les points d'intersections des vari\'et\'es invariantes des pi\`eces basiques d'un m\^eme cycle. J. Palis \cite{Pa} a  montr\'e que  ces intersections de vari\'et\'es invariantes formant un cycle sont responsables des $\Om$-explosions. Dans ce cadre, $\cR(f)$ appara\^\i t comme l'accroissement potentiel de $\Om(f)$, ou plus pr\'ecis\'ement comme la limite sup\'erieure des $\Om(g)$ pour $g$ convergeant vers $f$. 
C'est donc naturellement l'hyperbolicit\'e de l'ensemble r\'ecurrent par cha\^\i ne, plut\^ot que celle de l'ensemble non-errant qui est la clef de la stabilit\'e structurelle, et les hypoth\`eses de la th\'eorie de Smale s'expriment plus simplement en termes de $\cR(f)$:
\begin{equation*}
\begin{split}
f \mbox{ v\'erifie l'Axiome A et la condition ``pas de cycles"}&\Longleftrightarrow \\
\Om(f)=\overline{Per(f)} \mbox{, est hyperbolique et ne pr\'esente pas de cycles}&\Longleftrightarrow \\ 
\cR(f) \mbox{ est hyperbolique}.&\\
\end{split}
\end{equation*}
Cette condition est \'equivalente \`a l'$\Om$-stabilit\'e (stabilit\'e structurelle restreinte \`a $\Om(f)$), et \'e\-qui\-va\-len\-te de fait \`a la $\cR(f)$-stabilit\'e, {\em a priori} plus forte. La preuve de l'$\Om$-stabilit\'e passe par deux \'etapes:
\begin{itemize}
\item la notion de filtration, qui structure de fa\c con robuste la dynamique globale autour des pi\`eces basiques~: la vari\'et\'e est coup\'ee en tranches ordonn\'ees contenant chacune une unique pi\`ece basique ; les orbites respectent l'ordre sur les tranches~: si une orbite sort d'une tranche c'est pour entrer dans la suivante, et de ce fait elle ne reviendra plus visiter les tranches ant\'erieures (chaque d\'epart est sans retour)~; une orbite qui entre dans une tranche et qui n'en sort plus appartient \`a la vari\'et\'e stable de la pi\`ece basique correspondante. En particulier, toute  orbite s'accumule (positivement et n\'egativement) sur une unique pi\`ece basique.
\item l'\'etude locale de la dynamique au voisinage de chaque pi\`ece basique. Par exemple le lemme de pistage (shadowing lemma, voir \cite{Bo}) montre que toute pseudo-orbite dont les points appartiennent \`a une pi\`ece basique reste voisine de (``est pist\'ee par") une unique orbite de cette pi\`ece basique.  En cons\'equence, tout point dont l'orbite s'accumule  positivement sur une pi\`ece basique appartient \`a la vari\'et\'e stable d'un point de cette pi\`ece basique. 
\end{itemize}

Des filtrations structurant la dynamique globale existent en toute g\'en\'eralit\'e, associ\'ees \`a l'ensemble r\'ecurrent par cha\^\i nes~:

On d\'efinit une relation entre les points de $\cR(f)$ : deux points sont \'equivalents s'il existe des pseudo-orbites (de sauts arbitrairement petits) joignant l'un \`a l'autre et r\'eciproquement. L'ensemble $\cR(f)$ est ainsi l'union des classes d'\'equivalences qui sont des compacts invariants (nous les appellerons les {\em classes de r\'ecurrence par cha\^\i nes}). C. Conley (voir \cite{Co} et \cite[Fundamental Theorem of Dynamical Systems]{Ro})  montre l'existence de fonctions $\psi\colon M\to \RR$, dites {\em fonctions de Lyapunov}, strictement croissantes le long de toute orbite de $M\setminus \cR(f)$, et constantes le long de celles de $\cR(f)$. On peut de plus choisir $\psi$ telle que deux points de $\cR(f)$ soient dans un m\^eme niveau si et seulement s'ils sont \'equivalents et telle que l'image $\psi(\cR(f))$ soit un compact totalement discontinu de $\RR$. Les niveaux d'une telle fonction de Lyapunov $\psi$ permettent alors de construire des filtrations s\'eparant les classes de r\'ecurrence par cha\^\i nes.    

Cependant une classe de r\'ecurrence par cha\^\i nes n'est pas {\em a priori} une ``pi\`ece \'el\'ementaire de la dynamique" puisqu'elle peut parfois se d\'ecomposer en union de  compacts invariants plus petits.

\subsubsection{Dynamiques g\'en\'eriques}
La vari\'et\'e des comportements possibles d'un syst\`eme dynamique a amen\'e l'id\'ee d'\'eviter les pathologies fragiles (ne r\'esistant pas \`a d'infimes perturbations du syst\`eme) et de concentrer l'\'etude sur les comportements g\'en\'eriques.

Cette \'etude est intimement li\'ee aux r\'esultats perturbatifs~: ainsi le th\'eor\`eme de transversalit\'e de R. Thom se traduit pour les dynamiques g\'en\'eriques, en toute topologie $C^r$, par le th\'eor\`eme de Kupka-Smale~: $C^r$-g\'en\'eriquement, toute orbite p\'eriodique est hyperbolique et toute intersection entre vari\'et\'es invariantes d'orbites p\'eriodiques est une intersection transverse.  

Dans cet esprit, un r\'esultat essentiel de l'\'etude des dynamiques $C^1$-g\'en\'eriques est le c\'el\`ebre  lemme de fermeture de C. Pugh (Pugh's Closing Lemma \cite{Pu,PuRo}). En permettant de cr\'eer par ($C^1$-petite) perturbation une orbite p\'eriodique passant pr\`es d'un point non-errant arbitraire, il montre que, $C^1$-g\'en\'eriquement, l'ensemble non-errant $\Om(f)$ est l'adh\'erence de l'ensemble des points p\'eriodiques. Ceci entra\^\i ne que, $C^1$-g\'en\'eriquement, tous les candidats pour repr\'esenter la dynamique non-triviale se valent, sauf \'eventuellement l'ensemble r\'ecurrent par cha\^\i nes pour lequel la question restait ouverte~; en formule~:

$$
\mbox{$f$ est $C^1$-g\'en\'erique}\Longrightarrow \overline{Per(f)}=\overline{Rec_+(f)}=\overline{Rec_-(f)}=\overline{L(f)}=\Om(f)\subset \cR(f).
$$

Par analogie \`a la description donn\'ee par Smale de la dynamique globale des diff\'eomorphismes Axiome A sans cycles, on aimerait obtenir une ``d\'ecomposition en pi\`eces \'el\'ementaires de la dynamique non-triviale" d'un diff\'eomorphisme g\'en\'erique. Un candidat naturel est bien s\^ur l'ensemble r\'ecurrent par cha\^\i nes $\cR(f)$ et sa d\'ecomposition en classes d'\'equivalence de la th\'eorie de Conley (voir ci-dessus).

L'\'etude des dynamiques $C^1$-g\'en\'eriques a connu r\'ecemment de nombreux d\'eveloppements (voir \cite{BD2,Ar,CMP,CM,MP,Ab,We2}), \`a la suite d'un nouveau r\'esultat de $C^1$-perturbation, le lemme de connexion de S. Hayashi \cite{Ha} et ses variantes par L. Wen et Z. Xia \cite{WeXi} et par M.-C. Arnaud \cite{Ar}. 
Ce lemme de connexion permet de faire passer (par une $C^1$-perturbation du diff\'eomorphisme) l'orbite d'un point $x$ par un point $y$, si, pour la dynamique initiale, les orbites positive de $x$ et n\'egative de $y$ s'accumulent sur un m\^eme point non-p\'eriodique.

Il existe, pour les diff\'eomorphismes g\'en\'eriques, un second candidat pour la d\'ecomposition de la dynamique non-triviale en pi\`eces \'el\'ementaires.
Gr\^ace au lemme de connexion, l'ensemble non-errant d'un diff\'eomorphisme g\'en\'erique se d\'ecompose en union disjointe de compacts invariants {\em faiblement transitifs maximaux}~: pour tout couple de points d'un tel ensemble, tout voisinage du premier point poss\`ede des it\'er\'es positifs rencontrant un voisinage arbitraire du second de ces points, et r\'eciproquement. En effet \cite{Ar,GW}
montrent que cette relation sur les couples de points de $\Om(f)$ est, pour un diff\'eomorphisme g\'en\'erique, une relation d'\'equivalence  ferm\'ee dont les classes d'\'equivalence sont les ensembles faiblement transitifs maximaux (ce n'est pas vrai pour tous les diff\'eomorphismes). 
Cette nouvelle d\'ecomposition en pi\`eces \'el\'ementaires n'est pas ind\'ependante de celle donn\'ee par Conley~: en effet, un faiblement transitif maximal est toujours inclus dans une classe de r\'ecurrence par cha\^\i nes. Ces r\'esultats ne donnaient cependant pas de filtrations associ\'ees \`a cette partition de $\Om(f)$ en classes. 

Suivant une approche l\'eg\`erement diff\'erente, \cite{BD2} avait propos\'e les {\em classes homoclines } (a\-dh\'e\-ren\-ce des points homoclines transverses) des orbites p\'eriodiques comme \'etant les candidats naturels \`a \^etre les pi\`eces \'el\'ementaires de la dynamique~: les classe homoclines sont des ensembles transitifs canoniquement associ\'es aux orbites p\'eriodiques et co\"\i ncident avec les pi\`eces basiques dans le cas des dynamiques hyperboliques. De fait, \cite{Ar,CMP} montrent que les classes homoclines des diff\'eomorphismes g\'en\'eriques sont des {\em ensemble transitifs satur\'es}~: elles contiennent tout ensemble transitif qu'elles intersectent (en particulier elles sont disjointes ou confondues)~; leurs preuves montrent \'egalement que les classes homoclines sont (g\'en\'eriquement) des ensembles faiblement transitifs maximaux\footnote{Une classe homocline $H(P,f)$ d'un diff\'eomorphisme g\'en\'erique $f$ est l'intersection de l'adh\'erence des vari\'et\'es instable et stable du point $P$, d'apr\`es \cite{Ar}. Ces adh\'erences sont {\em stables au sens de Lyapunov} pour $f$ et $f^{-1}$, respectivement, d'apr\`es \cite{CMP}, c'est-\`a-dire qu'elles admettent  une base de voisinages positivement invariants pour $f$ et $f^{-1}$, respectivement. On en d\'eduit que tout ensemble faiblement transitif intersectant $H(P,f)$ est inclus dans des voisinages arbitrairement petits des adh\'erences des vari\'et\'es invariantes, et donc dans $H(P,f)$. }. F. Abdenur \cite{Ab} montre alors l'existence de deux $C^1$-ouverts disjoints $\cS$ et $\cA$ dont l'union est $C^1$-dense dans l'ensemble des diff\'eomorphismes et tels que~:
\begin{itemize}
\item les diff\'eomorphismes g\'en\'eriques de $\cS$ poss\`edent une infinit\'e de classes homoclines (on parle alors de {\em dynamique sauvage})~; 
\item  pour un diff\'eomorphisme g\'en\'erique $f$ de $\cA$,  l'ensemble non-errant $\Om(f)$ est l'union  d'un nombre fini de classes homoclines  disjointes, ce nombre \'etant localement constant  parmi les diff\'eomorphismes g\'en\'eriques~; de plus il existe une filtration de la vari\'et\'e associ\'ee \`a cette partition de $\Om(f)$ en  classes homoclines. On parle alors de {\em dynamique apprivois\'ee}.  Pour ces dynamiques apprivois\'ees, les ensembles $\cR(f)$ et $\Om(f)$ co\"\i ncident et de plus la partition de $\cR(f)$ en classes de r\'ecurrence par cha\^\i nes co\"\i ncide avec la partition  de $\Om(f)$ en ensembles transitifs qui sont les classes homoclines. 
\end{itemize}

Les syst\`emes apprivois\'es admettent donc une description de leur dynamique globale tr\`es proche de celles des diff\'eomorphismes de type Axiome A. La dynamique en restriction aux classes homoclines n'est cependant pas n\'ecessairement hyperbolique~: il s'agit bien d'une classe de diff\'eomorphismes plus vaste que celle des Axiome A. Des exemples de dynamiques apprivois\'ees non-hyperboliques peuvent \^etre construits \`a l'aide d'ensembles {\em robustement transitifs}, pour lesquels il existe de nombreux exemples locaux ou globaux.  Il n'est \`a ce jour  pas connu si les classes homoclines des dynamiques apprivois\'ees sont n\'ecessairement robustement transitives (dans ce cas l'ensemble des dynamiques apprivois\'ees serait un ouvert dense de $\cA$). Nous verrons cependant, en cons\'equence de nos r\'esultats, qu'elles v\'erifient une propri\'et\'e proche de la robuste transitivit\'e  que nous appellerons {\em robuste r\'ecurrence par cha\^\i nes}.

Le fait que le $C^1$-ouvert $\cS$ est non-vide pour les vari\'et\'es de dimension sup\'erieure ou \'egale \`a $3$ (autrement dit, l'existence de dynamiques sauvages) \`a \'et\'e montr\'e dans \cite{BD2}. Ce n'est pas connu en dimension $2$ pour la topologie $C^1$ consid\'er\'ee ici~; rappelons qu'en topologie $C^2$, la coexistence g\'en\'eriques sur les surfaces d'une infinit\'e de classes homoclines (en fait, des puits ou des sources) est connue depuis \cite{N,Ne}. De plus, ruinant les espoirs de \cite{BD2},  \cite{BD3} montre l'existence de diff\'eomorphismes localement g\'en\'eriques poss\`edant une famille non-d\'enombrable  d'ensembles transitifs satur\'es, stables au sens de Lyapunov pour les temps positifs et n\'egatifs, et qui ne contiennent aucune orbite p\'eriodique~: de fait, pour les exemples connus, la dynamique en restriction \`a ces ensembles transitifs est minimale.

\vskip 2mm
Dans cet article, nous montrons un nouveau lemme de perturbation (voir le th\'eor\`eme~\ref{t.connect})~: si deux points sont joints par des pseudo-orbites de sauts arbitrairement petits, alors il existe des $C^1$-petites perturbations de la dynamique telles que l'orbite positive du premier point passe par le second de ces points. De plus, dans le cas conservatif, ces perturbations peuvent \^etre choisies pr\'eservant le volume. En particulier, comme dans le cas des $\Om$-explosions associ\'ees aux cycles des diff\'eomorphismes Axiome A, l'ensemble $\cR(f)$ est l'accroissement potentiel de $\Om(f)$.

Ceci nous permet de montrer que, g\'en\'eriquement, les ensembles $\cR(f)$ et $\Om(f)$ co\"\i ncident et que de plus la partition de $\cR(f)$ en classes de r\'ecurrence par cha\^\i nes co\"\i ncide (g\'en\'eriquement) avec la partition (d\'efinie g\'en\'eriquement)  de $\Om(f)$ en ensembles faiblement transitifs. Les deux candidats que nous proposions comme pi\`eces \'el\'ementaires de la dynamique sont g\'en\'eriquement les m\^emes. D'apr\`es la th\'eorie de Conley, il existe donc une filtration s\'eparant les pi\`eces \'e\-l\'e\-men\-tai\-res, ce qui montre que toute orbite positive ou n\'egative s'accumule sur une unique de ces pi\`eces (en fait, sur une partie de cette pi\`ece). Chaque pi\`ece \'el\'ementaire est faiblement transitive\footnote{Ceci sugg\`ere d'appeler {\em pi\`ece \'el\'ementaire de la dynamique} d'un diff\'eomorphisme $f$ quelconque,  toute classe de r\'ecurrence par cha\^\i nes qui est faiblement transitive.}. De plus \cite{CM} montrent que pour un diff\'eomorphisme g\'en\'erique, les points g\'en\'eriques de la vari\'et\'e ont leur orbite positive (ou n\'egative) qui s'accumule sur toute une pi\`ece \'el\'ementaire, qui de plus est stable au sens de Lyapunov.

Nous donnons alors une liste de cons\'equences plus ou moins directes, dont la plus spectaculaire est sans doute dans le cas conservatif~: les diff\'eomorphismes $C^1$-g\'en\'eriques d'une vari\'et\'e compacte connexe qui pr\'eservent le volume sont transitifs. Ceci est particuli\`erement remarquable pour les diff\'eomorphismes pr\'eservant l'aire des surfaces, pour lesquels on sait que ce r\'esultat est faux en topologie $C^4$. En effet, la th\'eorie KAM assure l'existence de $C^4$-ouverts de diff\'eomorphismes poss\'edant des disques p\'eriodiques (qui emp\^echent la transitivit\'e) (voir par exemple \cite[section II.4.c]{M} ou \cite[chapitre IV]{He0}).

\subsection{\'Enonc\'e pr\'ecis des r\'esultats}
\subsubsection{\'Enonc\'e du lemme de connexion pour les pseudo-orbites}

Dans tout ce travail nous consid\'erons une vari\'et\'e compacte $M$, munie d'une m\'etrique riemannienne arbitraire, et parfois munie d'une forme volume $\omega$ (sans relation avec la m\'etrique). Nous notons $\diff^1(M)$ l'ensemble des diff\'eomorphismes de classe $C^1$ de $M$ muni de la topologie $C^1$ et $\diff_\omega^1(M)\subset \diff^1(M)$ le sous-ensemble de ceux qui pr\'eservent le volume $\omega$.

Rappelons que, pour un espace m\'etrique complet, une partie est dite {\em r\'esiduelle} si elle contient une intersection d\'enombrable d'ouverts denses. Une propri\'et\'e est dite {\em g\'en\'erique} si elle est v\'erifi\'ee sur un ensemble r\'esiduel. Nous parlerons ici, par abus de langage, de diff\'eomorphismes g\'en\'eriques~: la phrase {\em ``un diff\'eomorphisme g\'en\'erique v\'erifie la propri\'et\'e $P$"} signifie que la propri\'et\'e $P$ est g\'en\'erique.

Soit $f\in \diff^1(M)$ un diff\'eomorphisme de $M$.  Pour tout $\varepsilon >0$, une {\em $\varepsilon$-pseudo-orbite} de $f$ est une suite (finie ou infinie) de points  $(x_i)$ v\'erifiant $d(x_{i+1},f(x_i))<\varepsilon$.
On d\'efinit les relations binaires suivantes pour les paires $(x,y)$ de points de $M$~:

\begin{itemize}
\item Pour tout $\varepsilon>0$, on note $x\dashv_\varepsilon y$ s'il existe une $\varepsilon$-pseudo-orbite $x_0=x,\dots, x_k=y$ avec $k\geq 1$.
\item On note $x\dashv y$ si, pour tout $\varepsilon>0$, on a $x\dashv_\varepsilon y$. On notera parfois $x\dashv_f y$ pour pr\'eciser le syst\`eme dynamique consid\'er\'e.
\item On note $x\prec y$ (ou $x\prec_f y$) si, pour tous voisinages $U,V$ de $x$ et $y$, respectivement, il existe $n\geq 1$ tel que $f^n(U)$ rencontre $V$. 
\end{itemize}

Voici quelques propri\'et\'es \'el\'ementaires de ces relations.

\begin{enumerate}
\item
Les relations $\dashv$ et $\dashv_\varepsilon$ sont, par construction, transitives. L'ensemble  r\'ecurrent par cha\^\i nes $\cR(f)$ est l'ensemble des points $x$ de $M$ tels que $x\dashv x$. 
\item  La relation $x\prec y$ n'est pas {\em a priori} transitive. L'ensemble non-errant $\Om(f)$ est l'ensemble des points $x$ de $M$ tels que $x\prec x$. 
\end{enumerate}

Marie-Claude Arnaud montre dans \cite{Ar} que la relation $\prec$ est transitive pour les dif\-f\'eo\-mor\-phis\-mes g\'en\'eriques.
Par des m\'ethodes similaires nous montrons~:

\begin{theo}\label{t.equivalence} Il existe une partie r\'esiduelle $\cG$ de $\diff^1(M)$ (ou de $\diff_\omega^1(M)$) telle que pour tout diff\'eomorphisme $f$ de $\cG$ et tout couple $(x,y)$ de points de $M$ on a~:
$$x\dashv_f y\Longleftrightarrow x\prec_fy .$$
\end{theo}
Ce th\'eor\`eme est la cons\'equence pour les dynamiques g\'en\'eriques du r\'esultat perturbatif suivant~:

\begin{theo}\label{t.connect}
Soit $f$ un diff\'eomorphisme d'une vari\'et\'e compacte  $M$ v\'erifiant l'une des deux hypoth\`eses suivantes~:
\begin{enumerate}
\item toutes les orbites p\'eriodiques de $f$ sont hyperboliques,
\item $M$ est une surface compacte et toute orbite p\'eriodique est ou bien hyperbolique ou bien elliptique irrationnelle
(ses valeurs propres sont complexes, de module $1$, mais ne sont pas des racines de $1$).
\end{enumerate}
Soit $\cU$ un $C^1$-voisinage de $f$ dans $\diff^1(M)$ (ou dans $\diff^1_\omega(M)$, si $f$ pr\'eserve une forme volume $\omega$). Alors, pour toute paire $(x,y)$ de points de $M$ telle que $x\dashv y$, il existe un diff\'eomorphisme $g$ dans $\cU$ et un entier $n>0$ tel que $g^n(x)=y$.
\end{theo}
\begin{rema}\label{r.connect} Dans le th\'eor\`eme~\ref{t.connect} ci-dessus, si le diff\'eomorphisme $f$ est de classe $C^r$, avec $r\in (\NN\setminus\{0\})\cup\{\infty\}$, la $C^1$-perturbation $g$ peut \^etre choisie de classe $C^r$.  En effet le diff\'eomorphisme $g$ est obtenu gr\^ace \`a un nombre fini de $C^1$-perturbations donn\'ees par le connecting lemma (th\'eor\`eme~\ref{t.connecting}), chacune de ces perturbations \'etant elle-m\^eme de classe $C^r$. 
\end{rema}

Voici quelques cons\'equences de ces r\'esultats~:

\begin{coro}\label{c.generic} Il existe une partie r\'esiduelle $\cG$ de $\diff^1(M)$ telle que pour tout diff\'eomorphisme $f$ de $\cG$, son ensemble $\cR(f)$ de points r\'ecurrents par cha\^\i nes co\"\i ncide avec son ensemble $\Om(f)$ de points non-errants.
\end{coro}

\begin{coro}\label{c.homocline} Supposons $M$ connexe. Il existe une partie r\'esiduelle $\cG$ de $\diff^1(M)$ telle que si $f\in\cG$ v\'erifie $\Om(f)=M$ alors il est transitif. De plus $M$ est alors une unique classe homocline.
\end{coro}

Un argument classique (voir~\cite{T}) permet de renforcer l\'eg\`erement l'\'enonc\'e du corollaire~\ref{c.homocline}~: un diff\'eomorphisme g\'en\'erique $f$ tel que $\Om(f)=M$ est {\it topologiquement m\'elangeant}, i.e. pour tous ouverts $U$ et $V$ non vides de $M$, tout it\'er\'e positif assez grand de $U$ rencontre $V$.

\subsubsection{D\'ecomposition de la dynamique des diff\'eomorphismes g\'en\'eriques en pi\`eces \'el\'ementaires}\label{s.generique}
Consid\'erons la relation sym\'etris\'ee $\vdashv$ de $\dashv$ d\'efinie par $x\vdashv y$ si $x\dashv y$ et $y\dashv x$. Cette relation induit une relation d'\'equivalence sur $\cR(f)$, dont les classes d'\'equivalence sont appel\'ees {\em classes de r\'ecurrence par cha\^\i nes}. 

On dit qu'un compact $\Lambda$ invariant par $f$ est {\em faiblement transitif} si, pour tous $x,y\in\La$,  on a $x\prec y$. Un ensemble $\La$ est {\em faiblement transitif maximal} s'il est maximal pour $\subset$ parmi les ensembles faiblement transitifs. L'adh\'erence d'une union croissante de compacts faiblement transitifs \'etant faiblement transitive, le lemme de Zorn implique que tout faiblement transitif est inclus dans un faiblement transitif maximal. Dans le cas o\`u la relation $\prec_f$ est transitive (c'est \`a dire pour $f$ g\'en\'erique) les faiblement transitifs maximaux sont les classes d'\'equivalence de la relation sym\'etris\'ee de $\prec$ induite sur l'ensemble $\Om(f)$. Pour un diff\'eomorphisme g\'en\'erique on obtient alors~:

\begin{coro} Il existe une partie r\'esiduelle $\cG$ de $\diff^1(M)$ telle que pour tout $f\in \cG$ les classes de r\'ecurrence par cha\^\i nes sont exactement les faiblement transitifs maximaux de $f$. 
\label{c.coincide}
\end{coro}

Les r\'esultats de Conley sur la d\'ecomposition de $\cR(f)$ en classes de r\'ecurrence par cha\^\i nes vont donc s'appliquer (pour les diff\'eomorphismes g\'en\'eriques) \`a la d\'ecomposition de $\Om(f)$ en ensembles faiblement transitifs maximaux.
Rappelons donc ces r\'esultats (voir \cite{Co}).

Pour tout espace compact m\'etrique $X$ et tout hom\'eomorphisme $h\colon X\to X$, il existe une fonction continue $\phi\colon X\to \RR$, appel\'ee {\em fonction de Lyapunov de $h$},  qui est strictement croissante le long des orbites de $X\setminus \cR(h)$, constante sur chaque classe de r\'ecurrence par cha\^\i nes  et injective sur l'ensemble de ces classes. De plus l'image $\phi(\cR(h))$ est un compact totalement discontinu de $\RR$. Une telle fonction de Lyapunov permet de construire une filtration, adapt\'ee \`a $h$ et s\'eparant les classes de r\'ecurrence par cha\^\i nes~:  

On appelle {\em filtration adapt\'ee \`a $h$} une famille $\{X_i\}_{i\in I}$, index\'ee par une partie $I\subset \RR$, de  compacts  de $X$ v\'erifiant les deux propri\'et\'es suivantes~:
\begin{itemize}
\item $X_j$ est contenu dans l'int\'erieur de $X_i$ si $j>i$~;
\item pour tout $i$, $h(X_i)$ est contenu dans l'int\'erieur de $X_i$.
\end{itemize} 
On dira que la filtration $\{X_i\}$ s\'epare deux ensembles invariants $K_1$ et $K_2$ s'il existe $i$ tel que l'int\'erieur de $X_i$ et $X\setminus X_i$ contiennent chacun un et un seul de ces deux ensembles. 
L'ensemble maximal invariant $\bigcap_{n\in\ZZ} h^n(X_i)$ dans $X_i$ est un {attracteur topologique}\footnote{On appelle {\em attracteur topologique}  l'ensemble maximal invariant d'un compact $U$ dont l'image est contenue dans  l'int\'erieur de $U$~; un {\em r\'epulseur topologigue} est un attracteur topologique pour la dynamique inverse.} (en g\'en\'eral non-transitif). De m\^eme l'ensemble maximal invariant du compl\'ementaire de $X_j$ est un {r\'epulseur topologique}. L'ensemble maximal invariant dans $X_i\setminus X_j$ est donc l'intersection d'un attracteur et d'un r\'epulseur. Plus g\'en\'eralement, toute classe de r\'ecurrence par cha\^\i nes est l'intersection d'un ensemble Lyapunov stable pour $h$ et d'un ensemble Lyapunov stable pour $h^{-1}$ (les ensembles poss\'edant cette propri\'et\'es sont appel\'es {\em ensembles neutres} par \cite{CMP}).  

Si $\phi$ est une fonction de Lyapunov de $h$, on peut obtenir une filtration de $X$, adapt\'ee \`a $h$ et s\'eparant les classes de r\'ecurrence par cha\^\i nes, de la fa\c con suivante~: on consid\`ere  une partie $I$ de $\RR\setminus\phi(\cR(h))$  telle que toute composante connexe de $\RR\setminus \phi(\cR(h))$ contient exactement un point de $I$. La famille $\{\phi^{-1}([i,+\infty[)\}_{i\in I}$ est la filtration annonc\'ee. 

Comme cons\'equence directe de nos r\'esultats et de la th\'eorie de Conley, nous obtenons donc~:
\begin{coro} Il existe une partie r\'esiduelle $\cG$ de $\diff^1(M)$ telle que tout $f\in \cG$ poss\`ede une filtration adapt\'ee qui s\'epare les ensembles faiblement transitifs maximaux.
\end{coro}

On aimerait que, pour tout diff\'eomorphisme g\'en\'erique, tout point g\'en\'erique appartienne au bassin d'un attracteur topologique transitif. Nos r\'esultats ne semblent pas r\'epondre directement \`a ce probl\`eme. Cependant une version affaiblie de cette question a \'et\'e propos\'ee par M. Hurley dans \cite{Hu}. Il d\'efinit un {\em quasi-attracteur} comme \'etant l'intersection d'une famille d'attracteurs topologiques. Il conjecture~:

\begin{conj}[{\bf Hurley}] Pour tout $r\geq 1$, il existe une partie r\'esiduelle $\cG$ de $\diff^r(M)$ telle que pour tout $f\in \cG$, l'union des bassins des quasi-attracteurs r\'ecurrents par cha\^\i nes est une partie r\'esiduelle de $M$.
\end{conj}
Il a d\'emontr\'e cette conjecture en topologie $C^0$ (pour les champs de vecteurs). \cite{Ar} puis \cite{MP} ont donn\'e des r\'esultats partiels dans ce sens qui nous permettent de montrer la conjecture en topologie $C^1$, avec le th\'eor\`eme~\ref{t.equivalence} et la proposition suivante~:

\begin{prop}\label{p.quasiattracteur} 
\begin{enumerate}
\item Pour tout hom\'eomorphisme $f$, un quasi-attracteur r\'ecurrent par cha{\^\i}\-nes est une classe de r\'ecurrence par cha\^\i nes stable au sens de Lyapunov.  
\item Il existe une partie r\'esiduelle $\cG$ de $\diff^1(M)$ telle que pour tout $f\in \cG$, une classe de r\'ecurrence par cha\^\i nes qui est stable au sens de Lyapunov est un quasi-attracteur (bien s\^ur r\'ecurrent par cha\^\i nes).
\end{enumerate}
\end{prop}

Voici donc la r\'eponse positive \`a la conjecture de Hurley en topologie $C^1$~:
\begin{coro}\label{c.hurley} Il existe une partie r\'esiduelle $\cG$ de $\diff^1(M)$ telle que pour tout $f\in \cG$, l'union des bassins des quasi-attracteurs r\'ecurrents par cha\^\i nes est une partie r\'esiduelle de $M$~; plus pr\'ecis\'ement, pour tout $f\in\cG$, l'ensemble des points dont l'$\omega$-limite est une classe de r\'ecurrence par cha\^\i nes stable au sens de Lyapunov est r\'esiduel. 
\end{coro}

Commentons ce r\'esultat  par quelques  remarques~:
\begin{rema}
\begin{enumerate}
\item
Pour un compact invariant $\La$, on d\'efinit son ensemble stable comme l'ensemble des points dont l'ensemble $\omega$-limite est contenu dans $\La$. Les vari\'et\'es stables des classes de r\'ecurrence par cha\^\i nes forment une partition de $M$. 
\item Comme il a \'et\'e remarqu\'e dans \cite{MP}, l'application qui \`a un  point associe l'adh\'erence de son orbite positive est semi-continue inf\'erieurement, et donc continue sur une partie r\'esiduelle de $M$. On en d\'eduit qu'il existe une partie r\'esiduelle de $M$ telle que, en restriction \`a cette partie, les ensembles stables des classes de r\'ecurrence par cha\^\i nes sont ferm\'es. 
\item De l'item pr\'ec\'edent on d\'eduit que, si l'ensemble stable d'une classe de r\'ecurrence par cha\^\i nes est d'adh\'erence non-maigre, alors l'ensemble stable lui-m\^eme est localement r\'e\-si\-du\-el (r\'esiduel dans un ouvert non-vide).
Ceci provient de ce que tout bor\'elien non-maigre, est localement r\'esiduel dans un ouvert non-vide (voir \cite[Proposition 8.26]{Ke}), et que la diff\'erence entre l'ensemble stable et son adh\'erence est un ensemble maigre, d'apr\`es l'item pr\'ec\'edent.   
\end{enumerate}
\end{rema}

\subsubsection{Classes de r\'ecurrence par cha\^\i nes et orbites p\'eriodiques}
Rappelons que d'apr\`es le lemme de fermeture (closing lemma) de C. Pugh, l'ensemble des points p\'eriodiques d'un diff\'eomorphisme g\'en\'erique est dense dans  $\Om(f)$, et on aimerait utiliser ces orbites p\'eriodiques pour mieux comprendre la dynamique des classes de r\'ecurrence par cha\^\i nes. Rappelons que la {\em classe homocline} $H(p,f)$ de l'orbite d'un point p\'eriodique hyperbolique $p$ est l'adh\'erence de l'ensemble des points d'intersection transverse de ses vari\'et\'es stable et instable. C'est un ensemble transitif, et comme nous l'avons vu \`a la section~\ref{s.generique}, les r\'esultats de \cite{CMP} entra\^\i nent que, pour un diff\'eomorphisme g\'en\'erique, toute classe homocline est un ensemble faiblement transitif maximal. Ceci montre, d'apr\`es le corollaire~\ref{c.coincide}~:
\begin{rema}\label{r.homoclasse}
Les classes homoclines d'un diff\'eomorphisme g\'en\'erique sont des classes de r\'ecurence par cha\^\i nes.
\end{rema}

Une classe de r\'ecurrence par cha\^\i nes d'un diff\'eomorphisme g\'en\'erique qui n'est pas une classe homocline ne contient donc aucune orbite p\'eriodique~; on appellera {\em classe ap\'eriodique} toute classe de r\'ecurrence par cha\^\i nes sans orbite p\'eriodique.   

\begin{coro}\label{c.interieur} Il existe une partie r\'esiduelle $\cG$ de $\diff^1(M)$ telle que, pour tout $f\in \cG$, toute composante connexe d'int\'erieur non-vide de $\Om(f)=\cR(f)$ est p\'eriodique et son orbite est une classe homocline.
\end{coro}

Le closing lemma de Pugh et la remarque~\ref{r.homoclasse} montrent~:
\begin{rema}
Pour $f$ g\'en\'erique, toute classe de r\'ecurrence par cha\^\i nes qui est isol\'ee dans $\cR(f)$ est une classe homocline. C'est en particulier le cas des classes qui sont des attracteurs ou des r\'epulseurs topologiques.
\end{rema} 
Une classe isol\'ee  est l'ensemble maximal invariant dans un voisinage $M_i\setminus \overline{M_j}$, appel\'e voisinage filtrant\footnote{Plus g\'en\'eralement, on appelle {\em voisinage isolant} d'un compact invariant $K$ tout voisinage ouvert $U$ de $K$ tel que  $K$ est l'ensemble maximal invariant de l'adh\'erence de $U$. On dira que $U$ est un {\em voisinage filtrant} s'il peut s'\'ecrire $U= U_1\setminus \overline{U_2}$ o\`u  $U_1$ et $U_2$ sont des ouverts v\'erifiant $\overline{U_2}\subset U_1$, $f(\overline U_1)\subset U_1$ et $f(\overline U_2)\subset U_2$.}.

On dit qu'un ensemble compact invariant $\La$ est {\em robustement transitif} s'il admet un voisinage isolant $U$ tel que, pour tout diff\'eomorphisme $g$ suffisamment $C^1$-proche de $f$, le maximal invariant $\La_g$ de $g$ dans $\overline{U}$ est encore transitif. On aimerait montrer que, pour un diff\'eomorphisme g\'en\'erique, les classes homoclines isol\'ees sont robustement transitives. Notre travail donne un r\'esultat dans ce sens~:

\begin{coro}\label{c.isolee} Il existe une partie r\'esiduelle $\cG$ de $\diff^1(M)$ telle que, pour tout $f\in \cG$, toute classe homocline $\La$   de $f$ isol\'ee dans $\cR(f)$ est {\em robustement r\'ecurrente par cha\^\i nes}~: pour tout voisinage isolant $U$ de $\La$, pour tout diff\'eomorphisme $g$ suffisamment $C^1$-proche de $f$, l'ensemble maximal invariant de $g$ dans $\overline U$ est r\'ecurrent par cha\^\i nes.
\end{coro}

Pour les classes non-isol\'ees, un travail r\'ecent (voir \cite{Cr}) pr\'ecise la fa\c con dont une classe de r\'ecurrence par cha\^\i nes est approch\'ee par les orbites p\'eriodiques~:
\begin{theor} Il existe une partie r\'esiduelle $\cG$ de $\diff^1(M)$ (ou de $\diff^1_\omega(M)$) telle que pour $f\in \cG$ tout ensemble faiblement transitif maximal est la limite, pour la topologie de Hausdorff, d'une suite d'orbites p\'eriodiques de $f$.
\end{theor}

De plus les classes de r\'ecurrence par cha\^\i nes ont une propri\'et\'e de semi-continuit\'e sup\'erieure~: si $x_n\in \cR(f)$ est une suite de points convergeant vers un point $x$ alors pour $n$ assez grand, la classe de $x_n$ est contenue dans un voisinage arbitrairement petit de la classe de $x$.

\subsubsection{Dynamiques conservatives}
Dans le cas conservatif, l'ensemble non-errant co\"\i ncide toujours avec la vari\'et\'e $M$. On obtient donc~:

\begin{theo}\label{t.conservatif} Supposons $M$ connexe.  Il existe une partie r\'esiduelle $\cG_\omega$ dans l'ensemble $\diff_\omega^1(M)$ des diff\'eomorphismes pr\'eservant $\omega$ pour laquelle tout diff\'eomorphisme $f\in\cG_\omega$ est transitif. 

De plus  $M$ est une unique classe homocline. 
\end{theo}
Comme pour le corollaire~\ref{c.homocline}, un diff\'eomorphisme conservatif g\'en\'erique est donc topologiquement m\'elangeant.
\begin{rema}
Ce r\'esultat \'etait d\'ej\`a connu en topologie $C^0$ pour les hom\'eomorphismes g\'en\'eriques.
C'est en effet une cons\'equence directe du th\'eor\`eme d'Oxtoby et Ulam~\cite{OU} qui montre que pour des hom\'emorphismes conservatifs g\'en\'eriques, la mesure de Lebesgue est ergodique.
\end{rema}

Nos r\'esultats permettent de r\'epondre \`a une question pos\'ee par M.~Herman dans \cite{He}~:

\begin{ques}[Herman]
Est-ce que tout diff\'eomorphisme $f$ de classe $C^\infty$ d'une vari\'et\'e compacte pr\'eservant une forme volume $\omega$ v\'erifie l'une des deux propri\'et\'es suivantes ?
\begin{enumerate}
\item ou bien $f$ poss\`ede une d\'ecomposition domin\'ee\footnote{On dit qu'un ensemble $f$-invariant $K$  poss\`ede une {\em d\'ecomposition domin\'ee} si le fibr\'e tangent $TM|_K$ \`a $M$ au-dessus de $K$ admet un scindement $TM(x)=E(x)\oplus F(x), x\in K$ en somme directe de deux sous-fibr\'es $Df$-invariants continus tels que l'expansion des vecteurs dans $E$ est uniform\'ement plus petite que dans $F$~: il existe un entier $\ell$ tel que, pour tout point $x\in K$ et tout couple de vecteurs non-nuls $u\in E(x)$ et $v\in F(x)$, on ait~:
$$ \frac{\|Df^\ell(u)\|}{\|u\|}\leq \frac12\frac{\|Df^\ell(v)\|}{\|v\|}
.$$
On dira que cette d\'ecomposition est {\em $\ell$-domin\'ee.}}~;
\item ou bien $f$ n'a pas d'{\em exposants stables}~:  $f$ est approch\'e en topologie $C^1$ par une suite de diff\'eomorphismes $(g_k)$ de classe $C^\infty$ poss\'edant une orbite $p_k$ de p\'eriode $n_k$ telle que tous les exposants   de $D g_k^{n_k}(p_k)$ (logarithme des valeurs propres divis\'e par la p\'eriode) sont plus petits que $1/k$ en module.
\end{enumerate}
\end{ques}

A.~Arbieto et C.~Matheus nous ont signal\'e le raisonnement suivant qui, comme cons\'equence directe de~\cite{BDP} et du th\'eor\`eme~\ref{t.conservatif}, r\'epond partiellement \`a cette question~:

\cite{BDP} montre que si $f$ ne peut pas \^etre approch\'ee en topologie $C^1$ par un diff\'eomorphisme ayant un point p\'eriodique dont la diff\'erentielle, \`a la p\'eriode, est l'identit\'e, alors il existe $\ell$ tel que, pour tout diff\'eomorphisme $g$ dans un $C^1$-voisinage de $f$, toute classe homocline poss\`ede une d\'ecomposition $\ell$-domin\'ee.
D'apr\`es le th\'eor\`eme~\ref{t.conservatif}, il existe une suite $h_k$ de diff\'eomorphismes de classe $C^1$ convergeant vers $f$ en topologie $C^1$, pour lesquels $M$ est une classe homocline, et poss\`ede donc une d\'ecomposition $\ell$-domin\'ee. Les d\'ecompositions $\ell$-domin\'ees passant \`a la limite (voir \cite[Corollary 1.5]{BDP}), la vari\'et\'e $M$ h\'erite d'une d\'ecomposition domin\'ee pour $f$. 

Ceci donne presque une r\'eponse positive \`a  la question d'Herman~: les diff\'eomorphismes $g_k$ de l'item 2 obtenus par ce raisonnement ne sont {\em a priori} que de classe $C^1$. Le ``lissage" des diff\'eomorphismes $g_k$ pose une difficult\'e~: \`a notre connaissance on ne sait toujours pas\footnote{Cette question est le probl\`eme 44 de la liste de \cite{PP}. E. Zehnder dans \cite{Z} donne une r\'eponse positive dans le cas o\`u $M$ est une surface~; plus g\'en\'eralement il montre la densit\'e des diff\'eomorphismes symplectiques $C^\infty$ dans l'espace des diff\'eomorphismes symplectiques $C^r$, pour tout $r\geq 1$, en toute dimension. Il donne aussi une r\'eponse partielle dans le cas des diff\'eomorphismes conservatifs d'une vari\'et\'e de dimension $\geq 3$~: les diff\'eomorphismes de classe $C^\infty$ pr\'eservant $\omega$ sont denses en topologie $C^r$ dans $\diff^{r+\alpha}_\omega(M)$, o\`u $\alpha>0$ et $r\geq 1$ est un entier. } si les diff\'eomorphismes de classe $C^r$, $r\in[2,+\infty]$, pr\'eservant la forme volume $\omega$ sont denses dans $\diff^1_\omega(M)$. En raffinant le raisonnement de Arbieto et Matheus, nous pouvons cependant donner une r\'eponse compl\`ete \`a la question d'Herman~:

\begin{theo}\label{t.herman}
Soit $r\in (\NN\setminus \{0\})\cup \{\infty\}$ et $f$ un diff\'eomorphisme de classe $C^r$ d'une vari\'et\'e compacte pr\'eservant une forme volume $\omega$. Alors, $f$ v\'erifie une des deux propri\'et\'es suivantes~:
\begin{enumerate}
\item ou bien $f$ poss\`ede une d\'ecomposition domin\'ee,
\item ou bien $f$ est approch\'e en topologie $C^1$ par une suite de diff\'eomorphismes $g$ de classe $C^r$ poss\'edant une orbite $p$ de p\'eriode $n$ telle que $D g^n(p)= Id$.
\end{enumerate}
\end{theo}

Ce r\'esultat est \`a rapprocher des travaux de J. Bochi (voir \cite{Boc}) qui montre que tout diff\'eomorphisme conservatif g\'en\'erique  d'une surface compacte,  ou bien est un diff\'eomorphisme d'Anosov, ou bien poss\`ede un ensemble de mesure totale de points dont les exposants de Lyapunov sont nuls.  J. Bochi et M. Viana (voir \cite{BV}) ont donn\'e une g\'en\'eralisation  en dimension plus grande~: pour Lebesgue-presque tout point $x$, ou bien les exposants de Lyapunov de  $x$ sont  tous nuls, ou bien la d\'ecomposition d'Oseledec (en espaces de Lyapunov) au-dessus de l'orbite de $x$ est domin\'ee.

\subsection{Probl\`emes ouverts}
Voici d'abord une liste de questions concernant la dynamique g\'en\'erique des diff\'eomorphismes $C^1$ (non-conservatifs) des vari\'et\'es compactes~:

\begin{prob}[Cas non-conservatif] Soit $f$ un diff\'eomorphisme $C^1$-g\'en\'erique de $\diff^1(M)$.
\begin{enumerate}
\item Une classe homocline d'int\'erieur non-vide de $f$ co\"\i ncide-t-elle avec la vari\'et\'e ?

{
\em
On a vu que, si  $\cR(f)$ est d'int\'erieur non-vide, il existe une classe homocline d'int\'erieur non-vide.
Une telle classe est stable au sens de Lyapunov pour $f$ et $f^{-1}$. On se demande donc plus g\'en\'eralement~:} 

Une classe homocline stable au sens de Lyapunov pour $f$ et $f^{-1}$ est-elle \'egale \`a toute la vari\'et\'e~?
\item Une classe homocline  stable au sens de Lyapunov  de $f$ est-elle un attracteur topologique~? 
{\em Dans le cas o\`u la r\'eponse serait n\'egative, on peut att\'enuer cette question de la fa\c con suivante~:}
Sa vari\'et\'e stable contient-elle un sous-ensemble localement r\'esiduel~?
\item Une classe homocline de $f$ qui est robustement r\'ecurrente par cha\^\i nes est-elle un ensemble robustement transitive~?

\item Est-ce que $f$ peut poss\`eder des ensembles faiblement transitifs maximaux qui ne sont pas transitifs~?
{\em On sait (voir \cite{BD3}) qu'il existe g\'en\'eriquement des ensembles faiblement transitifs maximaux qui ne sont pas des classes homoclines, mais les exemples connus sont transitifs.}
\item Est-ce que $f$ peut avoir un ensemble (infini mais) d\'enombrable d'ensembles faiblement transitifs maximaux~? Si un tel diff\'eomorphisme existe, poss\`ede-t-il des ensembles faiblement transitifs maximaux qui ne soient pas des classes homoclines~? 
\item Pour tout diff\'eomorphisme de $\diff^1(M)$, et pour tout point $x$ de $M$, existe-t-il une $C^1$-perturbation $g$ pour laquelle $x$ est sur une vari\'et\'e stable d'orbite p\'eriodique~? {\em Ceci aurait pour cons\'equence~: g\'en\'eriquement, les vari\'et\'es stables et instables d'orbites p\'eriodiques sont denses dans la vari\'et\'e.}
\item Est-ce que l'union des bassins des attracteurs topologiques faiblement transitifs de $f$ est dense dans la vari\'et\'e~?
{\em Un tel attracteur est une classe homocline et est transitif. Remarquons qu'une r\'eponse positive impliquerait une r\'eponse positive au probl\`eme pr\'ec\'edent~: en effet, d'apr\`es \cite{CMP} l'adh\'erence de la vari\'et\'e stable d'une orbite p\'eriodique est stable au sens de Lyapunov pour $f^{-1}$~; en cons\'equence, la vari\'et\'e stable de cette orbite est dense dans le bassin d'attraction de sa classe homocline. 

Voici deux remarques allant dans le sens d'une r\'eponse positive \`a cette question~:
\begin{itemize}
\item L'union des classes de r\'ecurrence par cha\^\i nes ap\'eriodiques est maigre~: d'apr\`es le corollaire~\ref{c.interieur}, ces classes sont contenues dans le bord de l'ensemble $\cR(f)$, qui est un ferm\'e d'int\'erieur vide.
\item Le corollaire~\ref{c.hurley} montre que pour un point $x$ g\'en\'erique d'un diff\'eomorphisme g\'en\'erique, l'ensemble $\omega(x)$ est un quasi-attracteur r\'ecurrent par cha\^\i nes. 
\end{itemize}
}
\end{enumerate}
\end{prob}
Passons au cas conservatif. Nous avons vu qu'un diff\'eomorphisme g\'en\'erique $f\in\diff^1_\omega(M)$ est transitif. Par cons\'equent, l'ensemble des points de $M$ d'orbite (positive et n\'egative) dense dans $M$ est r\'esiduel dans $M$. Il est naturel de chercher \`a comprendre la mesure de cet ensemble~:

\begin{prob}[Cas conservatif] Soit $f\in\diff^1_\omega(M)$ un diff\'eomorphisme conservatif g\'en\'erique.
\begin{enumerate}
\item L'ensemble des points de $M$, d'orbites positive et n\'egative par $f$ denses dans $M$, est-il de mesure de Lebesgue totale dans $M$~?
\item La mesure de Lebesgue est-elle ergodique pour $f$~? {\em Une r\'eponse positive \`a cette question impliquerait une r\'eponse positive \`a la pr\'ec\'edente. En topologie $C^0$, c'est le th\'eor\`eme d'Oxtoby-Ulam \cite{OU}.}
\end{enumerate}
\end{prob}

\subsection{Pr\'esentation de l'article}
\subsubsection{Id\'ee  de la preuve du th\'eor\`eme~\ref{t.connect}} 

\`A la base de tous les r\'esultats perturbatifs en topologie $C^1$, il y a la remarque que la taille du support d'une $C^1$-petite perturbation bougeant un point sur une distance $\delta$ est proportionnelle \`a $\delta$. \`A partir d'un r\'esultat d'alg\`ebre lin\'eaire, Pugh a montr\'e, pour son ``closing lemma", que, quitte \`a r\'epartir dans le temps cette perturbation, on peut rendre la constante de proportionnalit\'e arbitrairement proche de $1$. 

Pour fermer un segment d'orbite dont les extr\'emit\'es sont suffisamment proches, il suffit de s\'electionner un sous-segment dont les extr\'emit\'es sont proches l'une de l'autre (\`a distance $\delta$) mais sont ``\'eloign\'ees" (plus que la constante de proportionalit\'e fois $\delta$) des autres points interm\'ediaires de l'orbite~: c'est la strat\'egie suivie par Pugh pour le closing lemma. 

\vskip 2mm

Pour le connecting lemma, Hayashi a renonc\'e \`a r\'ealiser la connexion en une seule perturbation qui ne modifierait pas les orbites sur les points interm\'ediaires. Nous consid\'erons ici la version donn\'ee par M.-C. Arnaud dans \cite{Ar} de ce connecting lemma. \'Etant donn\'es deux segments d'orbites tels que l'extr\'emit\'e de l'un est proche de l'origine de l'autre,  ces segments d'orbites seront perturb\'es en plusieurs endroits, la strat\'egie \'etant de faire de petites perturbations permettant de raccourcir ces segments chaque fois qu'ils reviennent tr\`es pr\`es d'eux-m\^emes o\`u de l'autre segment~: lors du premier passage pr\`es d'un point, plut\^ot que de poursuivre l'orbite initiale, on saute directement sur le dernier passage pr\`es de ce point, abandonnant le segment d'orbite interm\'ediaire.

La difficult\'e est de s\'electionner ces sauts, afin de pouvoir les r\'ealiser par des perturbations \`a supports disjoints.  
Une premi\`ere notion de proximit\'e utilise un quadrillage, dans une carte au voisinage d'un point. Deux passages des orbites seront consid\'er\'es comme proches s'ils appartiennent \`a un m\^eme carreau de ce quadrillage. L'id\'ee est donc, \`a chaque passage dans un carreau, de sauter directement sur le dernier passage dans ce carreau. Une difficult\'e suppl\'ementaire provient de ce que les perturbations r\'ealisant ces sauts peuvent interf\'erer si elles concernent des carreaux adjacents. Comme nous l'avons dit \`a propos du closing lemma, ces perturbations sont r\'eparties dans le temps, i.e. le long d'it\'er\'es des carreaux. \`A chaque instant, les supports des perturbations sont tr\`es petits (bien plus petits que la taille des carreaux). Les interf\'erences ont lieu quand ces supports se rencontrent~: cela signife que les points concern\'es sont tous concentr\'es dans une boule tr\`es petite devant la taille des carreaux, ce qui d\'efinit une nouvelle notion de proximit\'e. On applique alors de nouveau l'id\'ee de sauter directement sur le point proche correspondant au dernier passage proche.

\vskip 2mm
Pour r\'esumer, le connecting lemma r\'ealise  par une vraie orbite une pseudo-orbite qui ne pr\'esentait qu'un seul saut entre deux orbites initiales. Pour ce faire, il utilise un quadrillage tel que ce saut se fasse dans un de ses carreaux. Cependant, la preuve utilise, et permet donc, un nombre arbitraire de sauts pourvu que chacun d'entre eux ait lieu dans un carreau. C'est cette version du connecting lemma dont nous nous servirons ici. Elle est \'enonc\'ee pr\'ecis\'ement \`a la section~\ref{s.boites} (th\'eor\`eme~\ref{t.connecting}). 

On aimerait utiliser cette id\'ee pour traiter tous les sauts d'une pseudo-orbite. La difficult\'e est que les sauts d'une pseudo-orbite n'ont aucune raison de se trouver dans les carreaux d'un m\^eme quadrillage. L'id\'ee que nous avons suivie est de couvrir l'espace des orbites de la dynamique initiale par des quadrillages disjoints (c'est le corollaire~\ref{c.coloriage}). Toute orbite passe d'un carreau \`a un autre en temps fini, ce qui permet, pour toute pseudo-orbite de regrouper ses sauts dans les diff\'erents carreaux des quadrillages. Quadrillage par quadrillage, le connecting lemma permet alors de supprimer tous les sauts de la pseudo-orbite. Comme dans le connecting lemma, la nouvelle orbite est en g\'en\'eral plus courte que la pseudo-orbite initiale, mais ses extr\'emit\'es sont les m\^emes. Nous obtenons ainsi un lemme de connexion pour les pseudo-orbites (c'est le th\'eor\`eme~\ref{t.connect}).

Les perturbations associ\'ees \`a chaque quadrillage sont r\'eparties sur les $N$ premiers it\'er\'es de ce quadrillage. On cherche de ce fait une famille de quadrillages dont les $N$ premiers it\'er\'es sont deux \`a deux disjoints et dont l'int\'erieur des carreaux rencontre toute orbite. La premi\`ere \'etape de cette construction consiste \`a construire un ouvert $U$ disjoint de ses $N$ premiers it\'er\'es et rencontrant toute orbite (bien s\^ur les orbites p\'eriodiques de petite p\'eriode devront suivre un traitement sp\'ecial). On a appel\'e ``tour topologique" un tel ouvert par analogie aux tours utilis\'ees en th\'eorie ergodique (c'est le th\'eor\`eme~\ref{t.coloriage}).

Voici le plan de construction de cette tour topologique en supposant ici pour simplifier qu'il n'y a pas d'orbite p\'eriodique~: on choisit un entier $K$ bien plus grand que $N$~; comme il n'y a pas d'orbite p\'eriodique, on peut recouvrir la vari\'et\'e par une famille finie d'ouverts $(V_i)$ chacun disjoint de ses $K$ premiers it\'er\'es. Pour construire $U$, on commence par garder $V_0$ et on l'appelle $U_0$. On \^ote de $V_1$ les points qui rencontrent $U_0$ en temps inf\'erieur \`a $K$ it\'er\'es. Le reste des points de $V_1$ est recouvert par une famille finie d'ouverts $(W_j)$ et un ``lemme de coloriage" permettra de choisir des entiers $k_j\leq K$ tels que $U_1=U_0\cup \bigcup_j f^{k_j}(W_j)$ soit un ouvert disjoint de ses $N$ premiers it\'er\'es. De la m\^eme fa\c con on \^ote de $V_2$ les points rencontrant $U_1$ en moins de $K$ it\'er\'es et on applique le m\^eme proc\'ed\'e pour construire un ouvert $U_2$ disjoint de ses $N$ premiers it\'er\'es et rencontrant toutes les orbites des points de $V_0\cup V_1\cup V_2$. L'ouvert $U$ annonc\'e est construit en incorporant successivement \`a l'ouvert $U_i$ les orbites de $V_{i+1}$ par la m\^eme construction.

\subsubsection{Organisation de cet article}

L'outil technique essentiel est le lemme de connexion de Hayashi. La version que nous utilisons est pr\'esent\'ee en d\'etail en section~\ref{s.boites}~: nous d\'efinissons tout d'abord la notion de bo\^\i te de perturbation (d\'efinition~\ref{d.boites})~; nous pr\'esentons alors le lemme de connexion comme un lemme d'existence de bo\^\i tes de perturbation (th\'eor\`eme~\ref{t.connecting}). La preuve, tr\`es proche de celle de M.-C. Arnaud est donn\'ee dans l'appendice~\ref{a.A}. Elle se fait en deux temps~: le premier est un r\'esultat perturbatif \'enonc\'e dans~\cite{Ar} et d\^u \`a C.~Pugh et C.~Robinson (section~\ref{a.cube})~; le second s'appuie sur deux lemmes de s\'election de pseudo-orbites (sections~\ref{s.connecting} et~\ref{s.suite2}) et utilise un r\'esultat de regroupement de points proximaux (section~\ref{s.boules}).

L'autre clef de la d\'emonstration, pr\'esent\'ee en section~\ref{s.coloriage} est un r\'esultat d'existence de tours topologiques (th\'eor\`eme~\ref{t.coloriage}, section~\ref{s.coloriage1}). La d\'emonstration (section~\ref{s.coloriage3}) utilise un lemme de coloriage (section~\ref{ss.coloriage}).

Ceci permet alors (section~\ref{s.connect}) de montrer le th\'eor\`eme~\ref{t.connect}. \`A partir de l'existence des tours topologiques, on construit un ensemble fini de bo\^\i tes de perturbations disjointes intersectant chaque orbite non-p\'eriodique (section~\ref{s.recouvre}). Toute pseudo-orbite ayant des sauts assez petits peut \^etre l\'eg\`erement modifi\'ee de fa\c con \`a concentrer ses sauts dans les carreaux des quadrillages des bo\^\i tes de perturbation. La propri\'et\'e des bo\^\i tes de perturbation permet ensuite de refermer ces pseudo-orbites en vraies orbites (sections~\ref{s.regroupe} et~\ref{s.fin}).

Nous avons rassembl\'e les cons\'equences du lemme de connexion en section~\ref{s.consequences}. Un argument de cat\'egorie de Baire (section~\ref{s.equivalence}) donne le th\'eor\`eme~\ref{t.equivalence}. Les autres r\'esultats se trouvent en section~\ref{s.consequences2}.

Le cas conservatif fait l'objet de la section~\ref{s.conservatif}. Lorsque la dimension de la vari\'et\'e est sup\'erieure ou \'egale \`a $3$ (section~\ref{s.conservatif3}), les d\'emonstrations sont presque inchang\'ees par rapport au cas g\'en\'eral : ceci provient du fait que g\'en\'eriquement les orbites p\'eriodiques sont toutes hyperboliques. La section~\ref{s.conservatif2} traite du cas de la dimension $2$ o\`u une \'etude au voisinage des orbites p\'eriodiques elliptiques doit \^etre r\'ealis\'ee.

\subsubsection{Remerciements}
Nous remercions Fran\c cois B\'eguin, Bassam Fayad, Marguerite Flexor, Fr\'ed\'eric Le~Roux, Enrique Pujals, Marcelo Viana et Jean-Christophe Yoccoz  pour leur \'ecoute attentive et plus particuli\`erement Marie-Claude Arnaud et Th\'er\`ese Vivier qui nous ont aid\'es \`a bien comprendre la preuve du Connecting Lemma ainsi que Flavio Abdenur pour les discussions sur la dynamique g\'en\'erique.
  
\section{Bo\^\i tes de perturbation et \'enonc\'e du ``connecting lemma"}\label{s.boites}
L'outil essentiel de la preuve est une version du connecting lemma que l'on peut extraire de la preuve donn\'ee dans \cite{Ar}, et que nous pr\'esentons dans cette section.

Pour tout $i\geq 0$, notons $\alpha_i =1+ \sum_{j=0}^{i} 2^{-j}$.
Appelons {\em carreau standard} de $\RR^d$ tout cube $C$  d\'efini de la fa\c con suivante~: 
\begin{itemize}
\item ou bien $C=[-1,1]^d$,
\item ou bien il existe des entiers $j_0\in\{1,\dots,d\}$ et $i\geq 0$ ainsi que, pour tout $j\in \{1,\dots,d\}$,  un entier $k(j)\in[-2^i\alpha_{i},2^i\alpha_{i}-1]$ tels que $k(j_0)\in\{-2^i\alpha_{i},2^i\alpha_{i}-1\}$ et tels que $C$ est le produit 
$$
C=\prod_{j=1}^d\left[\frac{k(j)}{2^i}, \frac{k(j)+1}{2^i}\right].
$$ 
\end{itemize}

\begin{figure}
\begin{center}
\epsfig{file=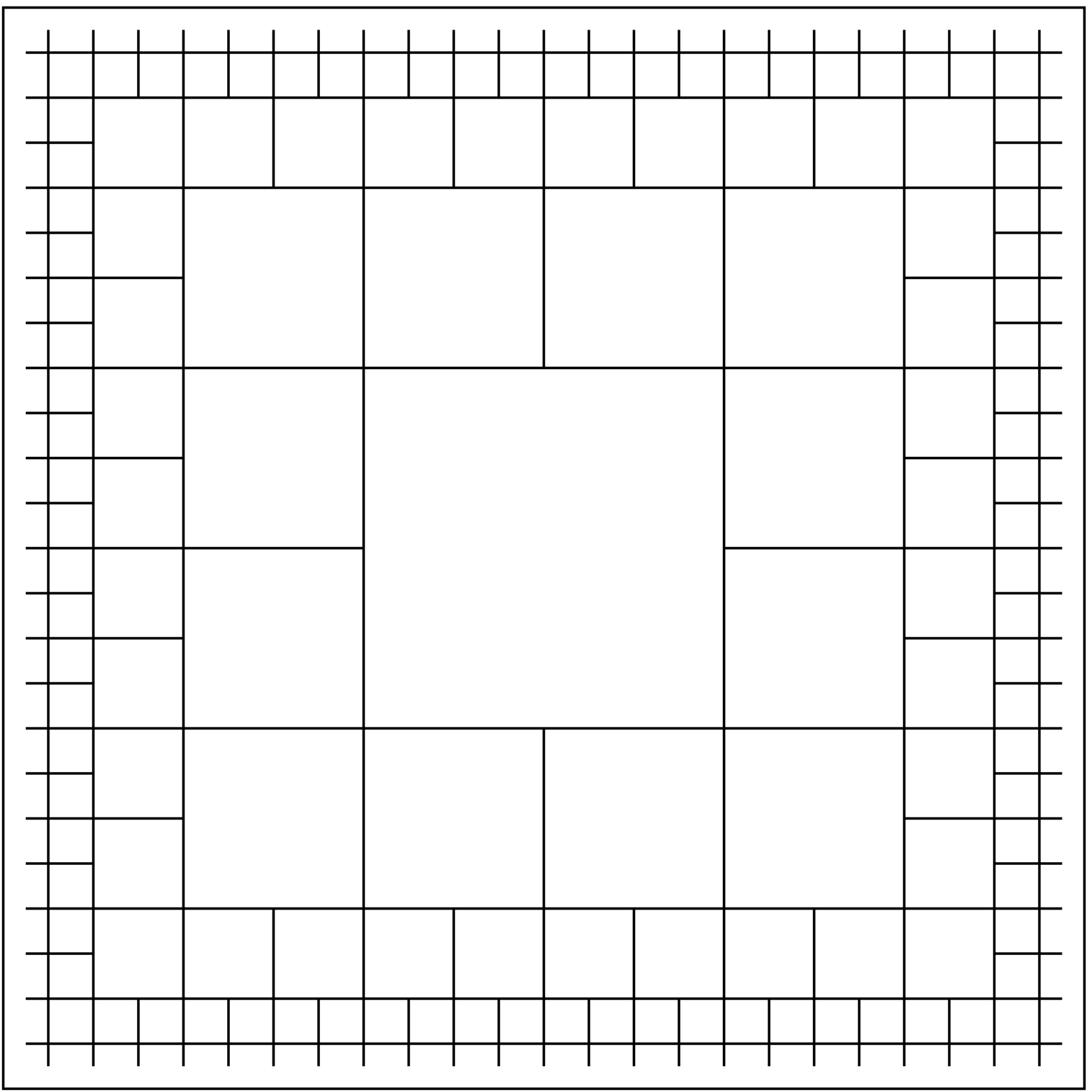, width=5cm}
\end{center}
\caption{Cube quadrill\'e de $\RR^2$. \label{cube}}
\end{figure}

Remarquons que les carreaux standards sont d'int\'erieurs deux \`a deux disjoints, que deux carreaux d'intersection non vide ont des longueurs de c\^ot\'es qui sont de rapport $2$, $1$, ou $\frac 12$, et finalement que l'union de tous les carreaux standards est le cube ouvert $\cC=]-3,3[^d$ (voir la figure~\ref{cube}). 

La collection de tous les carreaux standards est appel\'ee le {\em quadrillage standard}.  Nous appellerons {\em cube quadrill\'e} de $\RR^d$, l'image du cube ouvert $]-3,3[^d$ pav\'e de ses carreaux standards, par la compos\'ee d'une homoth\'etie et d'une translation de $\RR^d$.

Pour toute carte $\varphi\colon U\to\RR^d$ de la vari\'et\'e $M$, nous appellerons {\em cube quadrill\'e de $\varphi$} l'image inverse par $\varphi$ d'un cube quadrill\'e de $\RR^d$ contenu dans $\varphi(U)$. Ceci d\'efinit aussi les {\em carreaux} et le {\em quadrillage} de ce cube quadrill\'e.

\begin{defi}\label{d.boites} Soit $f$ un diff\'eomorphisme d'une vari\'et\'e compacte $M$, et $\cU$ un $C^1$-voisinage de $f$. Pour tout entier $N>0$, on appelle {\em bo\^\i te de perturbation $B$ d'ordre $N$ (pour $(f,\cU$))} tout cube quadrill\'e de $M$ disjoint de tous ses it\'er\'es $f^i(B)$, $i\in\{1,\dots,N\}$ et v\'erifiant la propri\'et\'e suivante.

Pour toute suite finie $\{(x_i,y_i)\}_{i\in\{1,\dots,\ell\}}$ de paires  de points de $B$ telle que pour tout $i\in\{1,\dots,\ell\}$ les points $x_i$ et $y_i$ appartiennent \`a un m\^eme carreau de $B$, il existe~:
\begin{itemize}
\item
un diff\'eomorphisme $g\in\cU$ co\"\i ncidant avec $f$ hors de l'union $\bigcup_{t=0}^{N-1} f^t(B)$, 
\item  une  suite strictement croissante $n_0=1<n_1<\cdots<n_s\leq \ell$,
\end{itemize}
tels que  $g^N(x_{n_k})=f^N(y_{n_{k+1}-1})$ pour tout $k\neq s$, et $g^N(x_{n_s})= f^N(y_\ell).$

On appelle {\em support de la bo\^\i te de perturbation $B$} l'union $supp(B)=\bigcup_0^N f^t(B)$
\end{defi}

Cette d\'efinition peut sembler abstraite et difficile \`a interpr\'eter. La remarque ci-dessous explique la propri\'et\'e essentielle des bo\^\i tes de perturbations~:

Soient $x$ et $y$ deux points hors du support de la bo\^\i te de perturbation $B$ et soit $x_0=x, x_1,\dots, x_n=y$ une $\varepsilon$-pseudo-orbite joignant $x$ \`a $y$. 

On dira que la pseudo-orbite {\em pr\'eserve le quadrillage} (de la bo\^\i te $B$) si les intersections de la pseudo-orbite avec le support de $B$ est une union de segments $x_i,x_{i+1}, \dots,x_{i+k}\dots x_{i+N}$ de la forme $x_i\in B$, $x_{i+k}=f^k(y_i)$, $k\in\{1,\dots,N\}$, o\`u $y_i$ est un point de $B$ appartenant au m\^eme carreau que $x_i$. 

On dira que cette pseudo-orbite {\em n'a pas de sauts dans $supp(B)$}  si les intersections de la pseudo-orbite avec le support de $B$ est une union de segment $x_i,x_{i+1}, \dots,x_{i+k}\dots x_{i+N}$ de la forme $x_i\in B$, $x_{i+k}=f^k(x_i)$, $k\in\{1,\dots,N\}$.

\begin{figure}
\begin{center}
\input{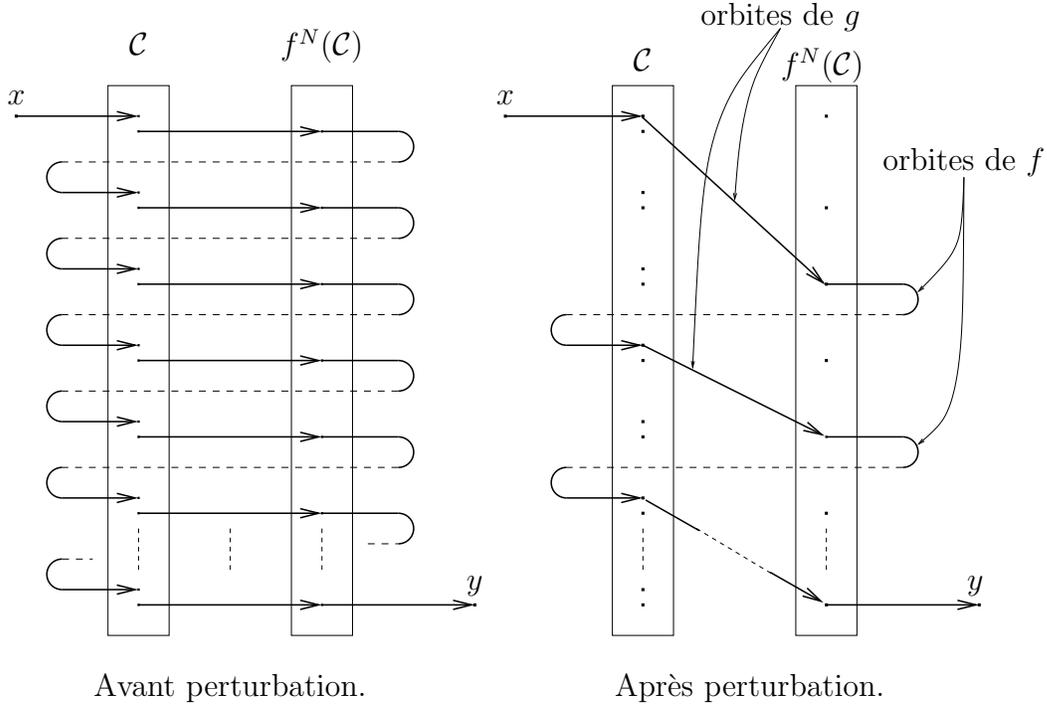}
\caption{Perturbation associ\'ee \`a une pseudo-orbite ayant tous ses sauts dans une m\^eme $B$.\label{f.raccourcit}}
\end{center}
\end{figure}

\begin{rema}\label{r.tuboite} Soit $B$ une bo\^\i te de perturbation d'ordre $N$ associ\'ee \`a $(f,\cU)$. 

Soient $x$ et $y$ deux points hors du support de la bo\^\i te $B$ et soit $x_0=x, x_1,\dots, x_n=y$ une $\varepsilon$-pseudo-orbite joignant $x$ \`a $y$ qui pr\'eserve le quadrillage de $B$. Alors,
 par d\'efinition d'une bo\^\i te de perturbation, il existe un diff\'eomorphisme $g\in\cU$, co\"\i ncidant avec $f$ hors de $\bigcup_{t=0}^{N-1} f^t(B)$, et  une $\varepsilon$-pseudo-orbite $z_0=x,z_1,\dots,z_m=y$ de $g$ qui n'a pas de saut dans le support de $B$. 

De plus $\{z_0,\dots,z_m\}\setminus \bigcup_1^{N-1} f^i(B)$ est form\'e de points de $\{x_0,\dots ,x_n\}$ (et plus pr\'ecis\'ement, de segments de la pseudo-orbite $(x_i)$ dont l'origine est ou bien $x$ ou bien appartient \`a $f^N(B)$ et l'extr\'emit\'e est ou bien $y$ ou bien appartient \`a $B$).
(La figure~\ref{f.raccourcit} pr\'esente un exemple de combinatoire lorsque tous les sauts ont lieu dans la m\^eme bo\^\i te $B$.)
\end{rema}
 
Si $\cB$ est une famille de bo\^\i tes de perturbation de supports deux \`a deux disjoints, on dira qu'une pseudo-orbite {\em n'a pas de saut hors des bo\^\i tes de $\cB$} si elle pr\'eserve les quadrillages  et si pour tout $x_i$ n'appartenant pas \`a l'une de ces bo\^\i tes, on a $x_{i+1}=f(x_i)$. Ceci s'applique aux point $x_i$ appartenant \`a l'image $f^j(B)$, avec $B\in\cB$ et $j\in\{1,\dots,N\}$.

De la remarque ci-dessus, on d\'eduit imm\'ediatement~:
\begin{lemm}\label{l.tuboite}
Soit $\cB$ une famille de bo\^\i tes de perturbation de supports deux \`a deux disjoints, soient $x$ et $y$ deux points hors du support de ces bo\^\i tes et soit $x_0=x, x_1,\dots, x_n=y$ une $\varepsilon$-pseudo-orbite joignant $x$ \`a $y$.
On suppose que cette pseudo-orbite pr\'eserve les quadrillages et n'a aucun saut hors des bo\^\i tes de $\cB$.

Consid\'erons une bo\^\i te $B\in \cB$ d'ordre $N$. Il existe alors un diff\'eomorphisme $g\in\cU$, co\"\i ncidant avec $f$ hors de $\bigcup_{t=0}^{N-1} f^t(B)$, et  une $\varepsilon$-pseudo-orbite $z_0=x,z_1,\dots,z_m=y$ de $g$ qui pr\'eserve les quadrillages et n'a aucun saut hors des bo\^\i tes de $\cB\setminus\{B\}$.
\end{lemm}

De la preuve du connecting lemma dans \cite{Ar} on extrait le r\'esultat suivant~:

\begin{theo}\label{t.connecting} Soit $f$ un diffeomorphisme d'une vari\'et\'e compacte $M$ de dimension $d$. Pour tout voisinage $\cU$ de $f$ il existe $N>0$ v\'erifiant~: pour tout point  $x\in M$, il existe  une carte locale $\varphi:U_x\to\RR^d$ en $x$  telle que tout  cube quadrill\'e disjoint de ses $N$ premiers it\'er\'es est une bo\^\i te de perturbation d'ordre $N$ pour $(f, \cU)$.

De plus, cette propri\'et\'e est encore v\'erifi\'ee par les cartes proches de $\varphi$ en topologie $C^1$.
\end{theo}

Une carte $\varphi\colon U\to \RR^d$ dont les cubes quadrill\'es sont des bo\^\i tes de perturbation sera appel\'ee une {\em carte de perturbation}.

Cet \'enonc\'e n'est pas \'ecrit explicitement dans \cite{Ar}~; de plus, le fait que $N$ est uniforme sur toute la vari\'et\'e n'est pas explicite dans \cite{Ar}, mais a \'et\'e remarqu\'e dans \cite{We}. C'est pourquoi nous donnerons dans l'Appendice~\ref{a.A} une preuve compl\`ete  du th\'eor\`eme~\ref{t.connecting}.

\begin{rema}
{\em On peut retrouver \`a partir du th\'eor\`eme~\ref{t.connecting} l'\'enonc\'e du connecting lemma de~\cite{Ar} pour les diff\'eomorphismes $C^1$ des vari\'et\'es compactes~:}

Soient $f$ un diff\'eo\-mor\-phisme $C^1$ d'une vari\'et\'e compacte $M$, $\cU$ un voisinage de $f$ dans $\diff^1(M)$, $p_0$ un point de $M$ qui n'est pas p\'eriodique pour $f$, et $U$ un voisinage de $p_0$ dans $M$. Il existe, d'apr\`es le th\'eor\`eme~\ref{t.connecting}, un entier $N>0$ et une bo\^\i te de perturbation $B$ d'ordre $N$ pour $(f,\cU$) contenue dans $U$ dont l'un des carreaux $C$ contient $p_0$ dans son int\'erieur. Notons $W$ l'int\'erieur de $C$ (pour reprendre les notations de \cite{Ar}). Soient $p,q$ deux points de $M$ hors du support de la bo\^\i te $B$ tels que l'orbite positive de $p$ et l'orbite n\'egative de $q$ rencontrent $W$. Alors, il existe $g\in \cU$ tel que $q$ appartiennent \`a l'orbite positive de $p$.

{\em En effet, si $f^{n_p}(p)$ et $f^{-n_q}(q)$ appartiennent \`a $W$, pour $n_p,n_q>0$, et si $p_1,\dots,p_n=f^{n_p}(p)$ et $q_m=f^{-n_q}(q),\dots q_1$ sont les passages successifs dans $B$ des orbites respectives de $p$ sur l'intervalle de temps $\{0,\dots,n_p\}$ et de $q$ sur l'intervalle de temps $\{-n_q,\dots,0\}$ on applique la propri\'et\'e de la bo\^\i te de perturbation aux couples de points $(p_1,p_1),(p_2,p_2),\dots,(p_{n-1},p_{n-1}),(p_n,q_m)$, $(q_{m-1},q_{m-1}),\dots,(q_1,q_1)$.}
\end{rema}
\begin{rema}\label{r.connecting}\begin{itemize}
\item Si le diff\'eomorphisme $f$ du th\'eor\`eme~\ref{t.connecting} pr\'eserve une forme volume $\omega$, alors les bo\^\i tes de perturbation construites par ce th\'eor\`eme permettent des $C^1$-perturbations qui pr\'eservent $\omega$.

\item Si le diff\'eomorphisme $f$ du th\'eor\`eme~\ref{t.connecting} est de classe $C^r$, $r\geq 1$ les bo\^\i tes de perturbation construites par ce th\'eor\`eme permettent des $C^1$-perturbations qui sont de classe $C^r$. Ceci est imm\'ediat dans le cas non-conservatif, puisque les diff\'eomorphismes $C^r$ sont denses parmi les diff\'eomorphismes $C^1$. Dans le cas conservatif, cela provient de la nature des perturbations \'el\'ementaires qui apparaissent dans la d\'emonstration (voir le lemme~\ref{l.perturb} et la remarque qui suit).  
\end{itemize}
\end{rema}
\section{Existence d'une tour topologique}\label{s.coloriage}
L'une des clefs de la preuve du th\'eor\`eme~\ref{t.connect} est  l'existence d'une famille finie de bo\^\i tes de perturbation d'ordre $N$, de supports disjoints, et dans lesquelles  toute orbite entre en temps fini. Remarquons que l'union de telles bo\^\i tes de perturbation forme un ouvert $U$, disjoint de ses $N$ premiers it\'er\'es, et rencontrant tout segment d'orbite suffisamment long. De fait, la construction d'un tel ouvert (appel\'e ici {\em tour topologique}) est l'\'etape essentielle qui, avec le th\'eor\`eme~\ref{t.connecting}, nous permettra de construire une telle famille de bo\^\i tes de perturbation. Le but de cette partie et de prouver l'existence d'une telle tour topologique. 

 Clairement, l'existence d'orbites p\'eriodiques de p\'eriode inf\'erieure \`a $N$ est une obstruction incontournable~: l'ouvert $U$ devant rencontrer ces orbites, ne peut \^etre disjoint de ses $N$ premiers it\'er\'es. C'est pourquoi, dans les \'enonc\'es suivants, nous particularisons les orbites p\'eriodiques (de p\'eriode $\leq N$) et leurs vari\'et\'es invariantes locales\footnote{Pour tout point p\'eriodique hyperbolique $x$ de $f$, on rappelle que sa vari\'et\'e stable
locale de taille $\delta>0$ est l'ensemble~:
$$
W^s_\delta(x)=\{y\in M, \quad \forall n\in\mathbb{N}, d(f^n(x),f^n(y))<\delta\}.
$$
On d\'efinit de m\^eme la vari\'et\'e instable locale $W^u_\delta(x)$.
}.

Dans le cas d'un compact $K$ sans orbites p\'eriodiques, nous montrons l'existence pour tout $N_0\in\NN$ d'un ouvert $U$ disjoint de ses $N_0$ premiers it\'er\'es, par lequel passe toute orbite de $K$~:  la compacit\'e de $K$ implique alors que les orbites de $K$ ont un temps de retour uniform\'ement born\'e dans $U$. Voici maintenant l'\'enonc\'e g\'en\'eral, tenant compte des orbites p\'eriodiques~:

\subsection{\'Enonc\'e du r\'esultat}\label{s.coloriage1}

\begin{theo}\label{t.coloriage}
Pour tout entier naturel $d$, il existe $\kappa_d\in \NN$ tel que pour tout diff\'eomorphisme $f$ d'une vari\'et\'e de dimension $d$, pour tout $N_0\in\NN$, pour toute constante $\delta>0$ et tout compact invariant $K$ ne contenant pas de point p\'eriodique non-hyperbolique de p\'eriode inf\'erieure \`a $\kappa_d N_0$, il existe un ouvert $U$ ayant les propri\'et\'es suivantes~:

\begin{enumerate}
\item Notons $Per_{N_0}(f)$ l'ensemble des orbites p\'eriodiques de p\'eriode inf\'erieure \`a $N_0$ de $K$. Pour tout point $x\in K$ tel que  $x\notin \bigcup_{p\in Per_{N_0}}W^s_{\delta}(p)$, il existe $n>0$ tel que $f^n(x)\in U$.
\item Pour tout point $x\in K\setminus \bigcup_{p\in Per_{N_0}}W^u_{\delta}(p)$ il existe $n>0$ tel que $f^{-n}(x)\in U$.
\item Les ferm\'es $\bar U,f(\bar U),\dots,f^{N_0}(\bar U)$ sont deux \`a deux disjoints.
\end{enumerate}
De plus, les composantes connexes de $\bar U$ peuvent \^etre choisies de diam\`etre arbitrairement petit.
\end{theo}

La d\'emonstration de ce r\'esultat fera l'objet de la section~\ref{s.coloriage}.

\begin{coro}\label{c.coloriage2}
Sous les hypoth\`eses et conclusions du th\'eor\`eme~\ref{t.coloriage}, il existe un compact $V\subset U$ tel que~:
\begin{itemize}
\renewcommand{\labelitemi}{-}
\item Pour tout point $x\in K$ tel que  $x\notin \bigcup_{p\in Per_{N_0}}W^s_{2\delta}(p)$, il existe $n>0$ tel que $f^n(x)\in V$.
\item Pour tout point $x\in K\setminus \bigcup_{p\in Per_{N_0}}W^u_{2\delta}(p)$ il existe $n>0$ tel que $f^{-n}(x)\in V$.
\end{itemize}
\end{coro}
\begin{demo} \`A partir du th\'eor\`eme~\ref{t.coloriage}, l'existence d'un tel compact peut se montrer de la fa\c con suivante~: on consid\`ere pour chaque point p\'eriodique de $Per_{N_0}$ un domaine fondamental compact de sa vari\'et\'e instable. On peut le choisir contenu dans une vari\'et\'e instable locale suffisamment petite pour qu'il ne rencontre pas les vari\'et\'es $W^s_\delta(Per_{N_0})$. Par cons\'equent, chaque point de ces domaines fondamentaux a un it\'er\'e positif rencontrant $U$. Par compacit\'e, il existe un voisinage $O^u$ dans $M$ de ces domaines fondamentaux et un compact $V^u_1$ de $U$ tel que tout point de $O^u$ a un it\'er\'e positif dans l'int\'erieur de $V^u_1$.

L'union des vari\'et\'es stables locales $W^s_{2\delta}(x)$, $x\in Per_{dN}$ et des it\'er\'es n\'egatifs de $O^u$  contient un voisinage ouvert $\cO^u$ de $W^s_\delta(Per_{N_0})$.  On remarque que $K\setminus \cO^u$ est compact. Par d\'efinition de $U$, tout point de ce compact poss\`ede un it\'er\'e positif dans $U$ et donc dans l'int\'erieur d'un compact $V^u_2\subset U$. Tout point de $\cO^u\setminus W^s_{2\delta}(Per_{N_0})$ a par construction un it\'er\'e positif dans $O^u$ et donc un autre it\'er\'e positif dans l'int\'erieur du compact $V^u_1\subset U$.

On construit de m\^eme des compacts $V^s_1$ et $V^s_2$ de $U$ tels que tout point de $K\setminus W^u_{2\delta}(Per_{N_0})$ poss\`ede un it\'er\'e n\'egatif dans l'int\'erieur de $V^s_1\cup V^s_2$. Le compact annonc\'e est l'union $V^u_1\cup V^u_2\cup V^s_1\cup V^s_2$.
\end{demo}

\subsection{Un lemme de coloriage}\label{ss.coloriage}

Soit $\{V,U_1,\cdots,U_\ell\}$ un ensemble fini de sous-vari\'et\'es compactes \`a bord d'int\'erieurs non vides contenues dans une vari\'et\'e $M$ de dimension $d$. Un {\em coloriage de $V\setminus \cup_i int(U_i)$ par $k$ couleurs} est la donn\'ee d'une famille finie $\cV=\{V_j\}_{j\in J}$ de sous-vari\'et\'es compactes \`a bord, d'int\'erieurs non vides, de $M$ et d'une fonction $c\colon j\mapsto c(j)\in \{1,\dots,k\}$ telles que

\begin{itemize}
\item la r\'eunion des int\'erieurs des $V_i$ recouvre le compact $V\setminus \cup_i int(U_i)$~;
\item pour toute paire $(j,j')$, $j\neq j'$, les {\em couleurs} $c(j)$ et $c(j')$ sont diff\'erentes d\`es que $V_j$ et $V_{j'}$ se rencontrent.
\end{itemize}

Nous appellerons {\em  $\gamma$-peinture de $\{U_i\}$ \`a $k$ couleurs} une fonction $c_0$ qui, \`a toute vari\'et\'e $U_i$, associe une partie de $\{1,\dots,k\}$, de cardinal born\'e par $\gamma$. On dira qu'un coloriage  $(\cV, c)$ de $V\setminus\cup_i int(U_i)$  {\em respecte $c_0$}, si pour tout $j\in J$ et toute vari\'et\'e $U_i$ rencontrant  $V_j$, la couleur $c(j)$ n'appartient par \`a la peinture $c_0(i)$.

\medskip
Une suite finie de sous-espace vectoriels d'un espace vectoriel de dimension $d$ sera dite {\em en position g\'en\'erale} si la somme  de leur codimension est  \'egale \`a la codimension de leur intersection. Remarquons que, si une suite  de sous-espaces vectoriels est en position g\'en\'erale, toute sous-suite l'est \'egalement.  Une suite $E_1,\dots,E_i,E_{i+1}$ est en position g\'en\'erale si et seulement si la suite $E_1,\dots,E_i$ est en position g\'en\'erale et si $E_{i+1}$ est transverse  \`a l'intersection $\bigcap_1^i E_j$.

Une suite de sous-vari\'et\'es $\{X_i\}_{i\in I}$ d'une vari\'et\'e $M$  de dimension $d$ sera dite {\em en position g\'en\'erale}, si pour toute partie finie $J\subset I$ et tout point $x\in\bigcap_{j\in J} X_j$, la suite $\{T_x(X_j)\}_{j\in J}$ des espaces tangents est en position g\'en\'erale. Pour des sous-vari\'et\'es de codimensions non-nulles et en position g\'en\'erale, un point de $M$ appartient \`a au plus $d$ de ces  sous-vari\'et\'es. Nous dirons que les sous-vari\'et\'es $\{X_i\}_{i\in I}$  sont {\em en position g\'en\'erale au voisinage d'un compact $V\subset M$} s'il existe un ouvert contenant $V$ tel que les intersections des $X_i$ avec cet ouvert sont en position g\'en\'erale.

\vskip 2mm
Le but de cette section~\ref{ss.coloriage} est de montrer~:

\begin{prop}{\bf[Lemme de coloriage]} \label{p.coloriage}

Pour tout $d\in\NN$, il existe $k_d \in \NN$ tel que tout $k\geq k_d$ poss\`ede la propri\'et\'e suivante:

\vskip 2mm

Soient $V,U_1,\dots,U_\ell$ des sous-vari\'et\'es compactes \`a bord de dimension $d$ d'une vari\'et\'e de dimension $d$ telles que les bords $\partial V,\partial U_1,\cdots,\partial U_\ell$ soient en position g\'en\'erale au voisinage de $V$. Pour tout $k\geq k_d$ et pour toute $2$-peinture $c_0$ de $\{U_1,\cdots,U_\ell\}$ \`a valeurs dans $\{1,\cdots,k\}$, il existe un coloriage $(\cV,c)$ de $V\setminus \cup_{i=1}^\ell int(U_i)$ qui respecte $c_0$.

\end{prop}

Notons $W=V\setminus \cup_{i=1}^\ell int( U_i)$. Comme les bords $\partial V,\partial U_1,\cdots,\partial U_\ell$  sont par hypoth\`ese en position g\'en\'erale au voisinage de $V$, on v\'erifie que $W$ est une vari\'et\'e compacte \`a bord et \`a coins, c'est \`a dire que tout point $x$ de $W$ poss\`ede un voisinage $O_x$ diff\'eomorphe \`a un ouvert $\hat O_x$ du cadran positif $\RR^d_+=[0,+\infty[^d\subset \RR^d$, chaque bord $\partial V\cap O_x,\partial U_1\cap O_x,\cdots,\partial U_\ell\cap O_x$ \'etant ou bien vide, ou bien correspondant \`a l'intersection avec $\hat O_x$ de l'un des hyperplans du bord de $\RR^d_+$. 

Puisque les bords $\partial V,\partial U_1,\cdots,\partial U_\ell$ sont en position g\'en\'erale au voisinage de $V$, l'intersection de $s$ d'entre eux est au voisinage de $V$ une vari\'et\'e (\'eventuellement vide) de dimension $d-s$~; les composantes connexes de l'intersection de cette vari\'et\'e avec $W$ seront appel\'ees  {\em face de dimension $d-s$} de $W$. On v\'erifie que chaque face de dimension $d-s$ est une vari\'et\'e \`a bord et \`a coins incluse dans $\partial W$. Les composantes connexes de $W$ seront appel\'ees  les {\em faces de dimension $d$}.

Le bord de $W$ est stratifi\'e par les faces de diff\'erentes dimensions~: le bord d'une face $\Si$ de dimension $j\geq 1$ est la r\'eunion de faces de dimension strictement inf\'erieure \`a $j$. Si $\Si'$ est une face incluse dans $\partial \Si$, on dira que $\Si$ est adjacente \`a $\Si'$.

\begin{lemm}\label{l.strate} Sous les hypoth\`eses de la proposition~\ref{p.coloriage}, et avec les notations ci-dessus, il existe pour tout $j\in\{0,\dots,d\}$ une vari\'et\'e compacte \`a bord $O_j$ de dimension $d$ ayant les propri\'et\'es suivantes~:
\begin{enumerate}
\item Pour chaque composante connexe $\De$ de $O_j$, l'ouvert $int(\Delta)\cup \bigcup_{j'=0}^{j-1} int(O_{j'})$ contient  exactement une face $\Si(\Delta)$ de dimension $j$ de $W=V\setminus \cup_{i=1}^\ell int( U_i)$. De plus l'application $\De\mapsto\Si(\De)$ est une bijection de l'ensemble des composantes connexes de $O_j$ sur l'ensemble des faces de dimension $j$. On dira que $\Delta$ est la composante associ\'ee \`a $\Si(\De)$.
\item Une composante $\Delta'$ de $O_{j'}$, avec $j'<j$, rencontre $\De$ si et seulement si $\Si(\De')$ est contenue dans le bord de $\Si(\De)$. 
\item Si la composante connexe  $\De$  rencontre une  des vari\'et\'es $U_i$, $i\in\{1,\dots,\ell\}$, alors la face correspondante $\Si(\De)$ est incluse dans le bord $\partial U_i$ de cette vari\'et\'e.
\end{enumerate}
\end{lemm}
\begin{demo}
On construit les vari\'et\'es $O_j$ par r\'ecurrence sur $j$.

Pour $j=0$, $O_0$ est un voisinage des faces de dimension $0$ (i.e. des points), constitu\'e de disques centr\'e en chacun de ces points. On le choisit suffisamment petit pour qu'une composante $\Delta$ de $O_0$ ne rencontre qu'une seule face $\Sigma(\De)$ de dimension $0$ et ne rencontre que les faces de dimension plus grande et les vari\'et\'es $U_i$ qui contiennent $\Sigma(\De)$ dans leur bord.

En supposant construits les vari\'et\'es $O_0,\dots,O_{j-1}$, on construit $O_j$ de la fa\c{c}on suivante~:

Pour toute face $\Si$ de $N$ de dimension $j$, on  peut choisir un compact connexe $\Si'$ de l'int\'erieur de $\Si$ qui contient le compact $\Si\setminus \cup_{j'=0}^{j-1} int( O_{j'})$.  Les compacts  $\Si'$ que l'on obtient ainsi sont deux \`a deux disjoints. La vari\'et\'e $O_j$ est l'union de voisinages connexes deux \`a deux disjoints des compacts $\Si'$, o\`u $\Si$ parcourt l'ensemble des faces de dimension $j$,  et choisis suffisamment petits pour que~:
\begin{itemize}
\item Une composante $\De$ de $O_j$ ne rencontre qu'une seule face de dimension $j$, not\'ee $\Si(\De)$.
\item $\De$ ne rencontre une composante $\De'$ de $O_{j'}$, pour $j'<j$, seulement si $\De'\cap \Sigma(\De)$ est non-vide et par hypoth\`ese de r\'ecurrence cela signifie que $\Si(\De')$ est contenue dans $\Si(\De)$.

\item Finalement $\De$ ne rencontre que les faces de dimension sup\'erieure et les vari\'et\'es $U_i$ qui contiennent $\Si(\De)$. 
\end{itemize}
\end{demo}

\begin{demo}[D\'emonstration de la proposition~\ref{p.coloriage}] 

Comme on ne cherche pas \`a trouver $k_d$ minimal, nous allons juste montrer que $k_d=(2d+1)^2$ convient. Nous consid\'erons la famille de vari\'et\'es $\{O_i\}$ donn\'ee par le lemme~\ref{l.strate}.

Fixons donc une $2$-peinture $c_0$ de $\{U_1,\dots,U_\ell\}$ \`a valeurs dans $\{1,\dots,k\}$ avec $k\geq (2d+1)^2$. 

Montrons \`a pr\'esent, par r\'ecurrence sur $i\in \{0,\dots,d\}$, que l'on peut associer \`a toute composante $\De$ de  $O_i$ une couleur $c(\De)\in \{i\cdot (2d+1)+1,\dots, (i+1)\cdot (2d+1)\}$ de fa\c con que le coloriage de $O_i$ ainsi obtenu respecte $c_0$. Comme les composantes connexes de $O_i$ sont par d\'efinitions disjointes, il n'y a pas de compatibilit\'es \`a v\'erifier  entre les couleurs des diff\'erentes composantes.

Consid\'erons une composante $\Delta$ de $O_i$, associ\'ee \`a une face $\Si(\De)$ de dimension $i$.  

Il y a, par hypoth\`ese, au plus $d-i$ faces de dimension $d-1$ qui sont adjacentes \`a $\Si(\De)$ et, d'apr\`es l'item (3) du lemme~\ref{l.strate}, la composante $\De$ rencontre au plus $d-i$ vari\'et\'es $U_j$. L'union $p(\Delta)$ des peintures associ\'ees aux $U_j$ rencontrant $\De$ est un ensemble de couleurs de cardinal inf\'erieur \`a $2d$. On peut donc associer \`a  $\Delta$ une couleur $c(\De)$ dans $\{i\cdot (2d+1)+1,\dots, (i+1)\cdot (2d+1)\}$ (de cardinal $2d+1$)  qui ne soit pas contenue dans $p(\De)$.

Montrons que la collection des composantes $\De$ des vari\'et\'es $O_j$ et la fonction $c$ ainsi construite d\'efinissent  bien un coloriage respectant $c_0$. Si deux composante $\De$ et $\De'$ se rencontrent, les faces associ\'ees sont par construction de dimension diff\'erentes $i\neq i'$. De ce fait les ensembles  $\{i\cdot (2d+1)+1,\dots, (i+1)\cdot (2d+1)\}$ et $\{i'\cdot (2d+1)+1,\dots, (i'+1)\cdot (2d+1)\}$ sont disjoints et les couleurs $c(\De)$ et $c(\De')$ sont donc diff\'erentes.  Comme l'union des int\'erieurs des $O_i$ contient $W$ (item (1) du lemme~\ref{l.strate}), la fonction $c$ induit donc bien un coloriage de $W$. Par construction, pour tout $\De$, la couleur $c(\De)$ n'appartient pas \`a l'union des peintures des $U_j$ qui rencontrent $\De$: le coloriage $c$ respecte la peinture $c_0$.

\end{demo}

\subsection{Mise en position g\'en\'erale}
\begin{lemm} \label{l.posigene}

Soit $n$ un entier, $M$  d'une vari\'et\'e de dimension $d$ et $f$ un diff\'eomorphisme de $M$ dont l'ensemble des points p\'eriodiques de p\'eriode inf\'erieure \`a $n$ est fini. 

Soit $S$ une vari\'et\'e compacte de dimension $d-1$ et soit $\psi\colon S\to  M$ un plongement de $S$. Il existe $\psi'$ arbitrairement proche de $\psi$ en topologie $C^r$, ($r\geq 1$ arbitraire si $\psi$ est $C^r$), tel que la famille des sous vari\'et\'es $f^i\circ \psi'(S)$, $i\in\{0,\dots,n\}$ soit en position g\'en\'erale.

\end{lemm}
\begin{demo}

Quitte \`a perturber $\psi$, on peut supposer que $\psi(S)$ ne contient aucun point de $Per_n(f)$, qui est fini par hypoth\`ese. Tout point $x\in S$ poss\`ede un voisinage $U_x$ tel que $\psi(U_x)$ est disjoint de ses $n$ premiers it\'eres.
Il existe donc des collections $\{U_i\}_{i\in I}$ et $\{V_i\}_{i\in I}$ d'ouverts de $S$ telles que $\bar V_i\subset U_i$ pour $i\in I$, que les $V_i$ recouvrent $S$ et que les $U_i$ aient la propri\'et\'e suivante~:

$$
{(*)}\quad 
\mbox{Pour tous $i,j\in I$ et $t\in\{1,\dots,n\}$ si $\psi(\bar U_i)$ intersecte $f^t\circ\psi(\bar U_j)$ alors $\bar U_i\cap \bar U_j= \emptyset$.}
$$

On consid\`ere successivement les $n+1$-uplets $(i_0,\dots,i_n)$ de $I^{n+1}$. Pour chacun de ces $n+1$-uplets, nous allons mettre en position g\'en\'erale les vari\'et\'es $(\psi(U_{i_0}),f\circ\psi(U_{i_1}),\dots,f^n\circ\psi(U_{i_n}))$ au voisinage des compacts $(\psi(\bar V_{i_0}),f\circ\psi(\bar V_{i_1}),\dots,f^n\circ\psi(\bar V_{i_n}))$. Cette propri\'et\'e persiste par perturbations de $\psi$. La mise en position g\'en\'erale pour les $n+1$-uplets suivants ne la d\'etruira pas. Apr\`es un nombre fini de perturbations, on obtiendra le r\'esultat voulu.

La mise en position g\'en\'erale de $(\psi(U_{i_0}),f\circ\psi(U_{i_1}),\dots,f^n\circ\psi(U_{i_n}))$ est possible car pour toute sous-famille de ce $n+1$-uplet de vari\'et\'es, d'apr\`es la propri\'et\'e ($*$) on a :
\begin{itemize}
\item  soit l'intersection de leur adh\'erence est vide, ce qui est stable par perturbation, 
\item  soit les adh\'erences des $U_i$ correspondants sont deux \`a deux disjointes, ce qui permet de faire des perturbations ind\'ependantes au voisinage de chaque $\bar V_i$. Dans ce dernier cas la mise en position g\'en\'erale est obtenue par le lemme de transversalit\'e de Thom. 
\end{itemize}
On applique donc inductivement ce proc\'ed\'e \`a l'ensemble des sous-familles du $n+1$-uplet $(i_0,\dots,i_n)$.
\end{demo}

\begin{rema}\label{r.posigene} Dans le lemme~\ref{l.posigene}, si l'on veut seulement que la famille des sous vari\'et\'es $f^i\circ \psi'(S)$, $i\in\{0,\dots,n\}$ soit en position g\'en\'erale au voisinage d'un compact $V$, il suffit de supposer que   $f$ ne poss\`ede qu'un nombre fini de points p\'eriodiques de p\'eriode inf\'erieure \`a $n$ au voisinage de $V$ (mais peut-\^etre un nombre infini, loin de $V$).
\end{rema}

\subsection{D\'emonstration du th\'eor\`eme~\ref{t.coloriage}}\label{s.coloriage3}
\begin{lemm}\label{l.union} Soient $M$ une vari\'et\'e compacte de dimension $d$ et $k_d$ la constante de la proposition~\ref{p.coloriage}. Soient $N>0$ un entier et $f$ un diff\'eomorphisme de $M$.

Soient $U$ et $V$ deux sous-vari\'et\'es compactes \`a bord, de dimension $d$. On suppose que $U$ est disjoint de ses $N$ premiers it\'er\'es positifs, et $V$ est disjoint de ses $2(k_d +1)N$ premiers it\'er\'es positifs. Alors, il existe $W$, sous-vari\'et\'e compacte \`a bord de dimension $d$, disjointe de ses $N$ premiers it\'er\'es, telle que $U\subset W$ et $V\subset \bigcup_{i=0}^{2(k_d+1)\cdot N}f^{-i}(W)$.

\end{lemm}
\begin{demo} 
Remarquons que, comme $V$ est disjoint de ses $2(k_d +1)N$ premiers it\'er\'es positifs, il n'existe aucune orbite de p\'eriode inf\'erieure ou \'egale  \`a $2(k_d +1)N$ au voisinage de $V$. 
D'apr\`es le lemme~\ref{l.posigene} et la remarque~\ref{r.posigene}, on peut augmenter l\'eg\`erement la vari\'et\'e \`a bord $U$ pour obtenir une vari\'et\'e \`a bord $\tilde U$ telle que la famille $\partial f^{-i}(\tilde U)$, $i\in\{0,\dots,2(k_d+1) N\}$ soit en position g\'en\'erale au voisinage de $V$. De plus, comme $U$ est un compact disjoint de ses $N$ premiers it\'er\'es, quitte \`a choisir $\tilde U\supset U$ suffisamment petit, on peut supposer que $\tilde U$ est aussi disjoint de ses $N$ premiers it\'er\'es. Quitte \`a augmenter $V$, on supposera que son bord est en position g\'en\'erale avec la famille $\partial f^{-i}(\tilde U)$, $i\in\{0,\dots,2(k_d+1)N\}$.

On consid\`ere l'intersection de $V$ avec les $f^{-i}(\tilde U)$ pour $i\in\{0,\dots,2(k_d+1)N-1\}$. On associe \`a $U_i=f^{-i}(\tilde U)$ les couleurs $c_0(i)=\{E(\frac{i}{2N}),E(\frac{i}{2N})+1\}$ qui appartiennent \`a $\{0,\cdots, k_d+1\}$, o\`u $E(\frac{i}{2N})$ repr\'esente la partie enti\`ere de $\frac{i}{2N}$.

D'apr\`es le proposition~\ref{p.coloriage}, il existe un coloriage $(\{V_j\},c)$ de $V\setminus \cup_i U_i$ qui respecte $c_0$ et \`a valeurs dans $\{1,\cdots,k_d\}$. Puisque $V$ est disjoint de ses $2(k_d+1)N$ premiers it\'er\'es, on peut supposer que l'union $\cup_j V_j$ est contenue dans un voisinage de $V$, lui aussi disjoint de ses $2(k_d+1)N$ premiers it\'er\'es. On pose alors

$$W=\tilde U \cup \bigcup_{j}f^{2c(j)N}(V_j).$$

Voyons maintenant que $W$ est disjoint de ses $N$ premiers it\'er\'es. Supposons par l'absurde que $x$ et $y=f^t(x)$, avec $t\in\{1,\dots,N\}$, appartiennent \`a $W$. Nous allons raisonner en envisageant toutes les positions possibles des points $x$ et $y$~:
\begin{itemize}
\item Les points $x$ et $y$ ne peuvent appartenir simultan\'ement \`a $\tilde U$ qui est disjoint de ses $N$ premiers it\'er\'es. 
\item De m\^eme, ils ne peuvent pas appartenir \`a un m\^eme ensemble $f^{2c(j)N}(V_j)$. 
\item Si $x$ et $y$ appartiennent \`a des ensembles $f^{2c(j_1)N}(V_{j_1})$ et $f^{2c(j_2)N}(V_{j_2})$ diff\'erents, alors, puisque $y=f^t(x)$, les it\'er\'es $f^{-2c(j_1)N}(x)$ et $f^{t-2c(j_2)N}(x)$ de $x$ appartiennent tous deux \`a $\cup_j V_j$. Par ailleurs, $2|c(j_1)-c(j_2)|N+t$ est inf\'erieur \`a $2(k_d-1)N+N$. Ceci contredit donc le fait que $\cup_j V_j$ est disjoint de ses $2k_dN$ premiers it\'er\'es.

\item Traitons enfin le cas o\`u $x$ appartient \`a $f^{2c(j)N}(V_j)$ et $y$ \`a $\tilde U$ (le cas o\`u $y$ appartient \`a $f^{2c(j)N}(V_j)$ et $x$ \`a $\tilde U$ lui est similaire)~: $V_j$ et $f^{-(2c(j)N+t)}(\tilde U)$ se rencontrent. La peinture associ\'ee \`a $U_{2c(j)N+t}= f^{-(2c(j)N+t)}(\tilde U)$ est $c_0(2c(j)N+t)=\{c(j),c(j)+1\}$. Cependant, puisque $c$ est un coloriage qui respecte la peinture $c_0$ et que $V_j$ rencontre $U_ {2c(j)N+t}$, la couleur  $c(j)$ ne doit pas appartenir \`a $c_0(2c(j)N+t)=\{c(j),c(j)+1\}$ ce qui est contradictoire.
\end{itemize}

Par le m\^eme raisonnement, les vari\'et\'es $\tilde U$ et $f^{2c(j)N}(V_j)$ sont deux \`a deux disjointes. Ceci montre que $W$ est une sous-vari\'et\'e compacte \`a bord de dimension $d$. D'autre part, par construction, $U\subset W$ et $V\subset \bigcup_0^{k_d\cdot N}f^{-i}(W)$, ce qui conclut la d\'emonstration.
\end{demo}

\begin{demo}[D\'emonstration du th\'eor\`eme~\ref{t.coloriage}]
Posons $\kappa_d=2k_d+1$. Puisque les orbites p\'eriodiques de p\'eriode plus petite que $\kappa_d N_0$ sont hyperboliques, il existe une vari\'et\'e compacte \`a bord $U_0$, disjointe de ses $N_0$ premiers it\'er\'es et disjointe de $Per_{\kappa_d N_0}(f)$ telle que

\begin{itemize}
\item pour tout $x\in W^s(Per_{\kappa_d N_0}(f))\setminus W^s_\delta(Per_{\kappa_d N_0}(f))$, il existe $n>0$ tel que $f^n(x)\in int(U_0)$~;

\item pour tout $x\in W^u(Per_{\kappa_d N_0}(f))\setminus W^u_\delta(Per_{\kappa_d N_0}(f))$, il existe $n>0$ tel que $f^{-n}(x)\in int(U_0)$.

\end{itemize}

En d'autres termes, $U_0$ est un  voisinage d'une sorte de domaine fondamental (non connexe) des vari\'et\'es stables et instables des points de $Per_{\kappa_d N_0}(f)$. Comme dans la preuve du corollaire~\ref{c.coloriage}, on peut montrer que l'union des it\'er\'e positifs de $U_0$ et de la vari\'et\'e instable locale $W^u_\delta(Per_{\kappa_d N_0}(f))$ contient un voisinage de $Per_{\kappa_d N_0}(f))$ (et de m\^eme pour ses it\'er\'es n\'egatifs).

Par cons\'equent, il existe un voisinage ouvert $O_0$ de $Per_{\kappa_d N_0}(f)$ tel que tout point de $O_0\setminus W^s_\delta(Per_{\kappa_d N_0}(f))$ a un it\'er\'e positif dans $U_0$ et tout point de $O_0\setminus W^u_\delta(Per_{\kappa_d N_0}(f))$ a un it\'er\'e n\'egatif dans $U_0$. Choisissons $O$ un voisinage ouvert de $Per_{\kappa_d N_0}(f)$ tel que pour tout $j\in\{0,\dots,2\kappa_dN_0\}$, son it\'er\'e $f^j(O)$ est inclus dans $O_0$.

Soit $K$ le compact invariant de l'\'enonc\'e du th\'eor\`eme~\ref{t.coloriage}.
Tout point $x$ de $K\setminus O$ poss\`ede un voisinage $U_x$ disjoint de ses $\kappa_dN_0$ premiers it\'er\'es. Par compacit\'e de $K\setminus O$, il existe donc une collection finie de vari\'et\'es \`a bord de dimension $d$, $\{V_0,\dots,V_r\}$, chacune disjointe de ses $\kappa_d N_0$ premiers it\'er\'es telles que leurs int\'erieurs recouvrent $K\setminus O$.

Consid\'erons le lemme~\ref{l.union} appliqu\'e \`a $U=U_0$, $V=V_0$ et $N=N_0$. On obtient une vari\'et\'e compacte \`a bord $U_1$ disjointe de  ses $N_0$ premiers it\'er\'es contenant $U_0$ et telle que $V_0\subset \bigcup_{j=0}^{\kappa_dN_0}f^{-j}(U_1)$. On construit par r\'ecurrence une suite de vari\'et\'es compactes \`a bord $(U_i)$ disjointes de  leurs $N_0$ premiers it\'er\'es,  contenant $U_0$ et telles que $\bigcup_{j=0}^iV_j\subset \bigcup_{j=0}^{\kappa_dN_0}f^{-j}(U_i)$.

Pour construire $U_{i+1}$, on applique le lemme~\ref{l.union} \`a $U=U_i$, $V=V_{i+1}$ et $N=N_0$. On obtient une vari\'et\'e compacte \`a bord $U_{i+1}$ disjointe de  ses $N_0$ premiers it\'er\'es contenant $U_i$ et telle que $V_{i+1}\subset \bigcup_{j=0}^{\kappa_dN_0}f^{-j}(U_{i+1})$. L'hypoth\`ese de r\'ecurrence appliqu\'ee \`a $U_i$ permet de montrer que $U_0$ est inclus dans $U_{i+1}$ et $\bigcup_{j=0}^iV_j$ est contenue dans $\bigcup_{j=0}^{\kappa_dN_0}f^{-j}(U_{i+1})$.

Posons $U=U_r$. Tout point $x$ de $K\setminus O$ a un it\'er\'e $f^i(x)$, $i\in\{0,\dots,\kappa_dN_0\}$, appartenant \`a $U$. Tout point $x$ de $K\setminus O_0$ poss\`ede un it\'er\'e positif dans $U$ et par ailleurs, $f^{-(\kappa_dN_0+1)}(x)$ appartient \`a $K\setminus O$ par d\'efinition de $O$. Il existe donc $i\in\{0,\dots,\kappa_dN_0\}$ tel que $f^{-(\kappa_dN_0+1)+i}(x)$ appartient \`a $U$. C'est un it\'er\'e n\'egatif de $x$. Ceci montre le r\'esultat pour les points de
$K\setminus O_0$.

Par d\'efinition de $O_0$, tout point de $O_0\setminus W^s_\delta(Per_{\kappa_d N_0}(f))$ a un it\'er\'e positif dans $U_0\subset U$ et, de m\^eme, tout point de $O_0\setminus W^u_\delta(Per_{\kappa_d N_0}(f))$ a un it\'er\'e n\'egatif dans $U$. Ceci
conclut la d\'emonstration pour les points de $K\cap O_0$.

\end{demo}
\section{Lemme de connexion pour les pseudo-orbites~: preuve du th\'eor\`eme~\ref{t.connect}}\label{s.connect}

\subsection{Bo\^\i tes de perturbation recouvrant l'espace des orbites}\label{s.recouvre}

\`A partir du  th\'eor\`eme~\ref{t.coloriage} on est \`a pr\'esent capable de construire une famille finie de bo\^\i tes de perturbation  de supports disjoints et  rencontrant toute orbite. Bien s\^ur les orbites p\'eriodiques et leurs vari\'et\'es invariantes locales doivent \^etre trait\'ees s\'epar\'ement~:
\begin{coro}\label{c.coloriage}
Soit $f$ un diff\'eomorphisme d'une vari\'et\'e compacte de dimension $d$, et $\cU$ voisinage de $f$ pour la topologie $C^1$. Soient $N$  l'entier associ\'e \`a l'ouvert $\cU$ par le th\'eor\`eme~\ref{t.connecting}, $\kappa_d$ l'entier donn\'e   par le  th\'eor\`eme~\ref{t.coloriage} et soit $\delta>0$ arbitraire. On note $N_0=10dN$.

Pour tout compact invariant $K$ ne contenant pas d'orbite p\'eriodique non-hyperbolique de p\'eriode inf\'erieure \`a $\kappa_dN_0$, il existe une famille finie $\cB_0$ de  bo\^\i tes de perturbation d'ordre $N$ pour $(f,\cU)$ et de supports deux \`a deux disjoints telle que toute orbite positive de points de $K\setminus \bigcup_{p\in Per_{N_0}}W^s_{\delta}(p)$ et tout orbite n\'egative de points de $K\setminus \bigcup_{p\in Per_{N_0}}W^u_{\delta}(p)$ rencontre l'int\'erieur d'un carreau d'un quadrillage d'une des bo\^\i tes de $\cB_0$.

De plus, pour tout $x\in M$, il existe $t\in \{0,\dots, N\}$ tel que $f^t(x)$ n'appartienne pas \`a l'union des supports de bo\^\i tes de $\cB_0$.

En cons\'equence, il existe une famille finie $\cC_0$ de carreaux associ\'es aux bo\^\i tes de $\cB_0$ et une famille finie $\cD_0$ de compacts contenus dans l'int\'erieur des carreaux de $\cC_0$ telles que,
\begin{itemize}
\item chaque carreau de $\cC_0$ contient exactement un compact de $\cD_0$ dans son int\'erieur~;
\item toute orbite positive de points de $K\setminus \bigcup_{p\in Per_{N}}W^s_{\delta}(p)$, et tout orbite n\'egative de point de $K\setminus \bigcup_{p\in Per_{N}}W^u_{\delta}(p)$ rencontre l'int\'erieur d'un compact de $\cD_0$.
\end{itemize}
\end{coro}

\begin{demo}
D'apr\`es le th\'eor\`eme~\ref{t.coloriage}  et le corollaire~\ref{c.coloriage2} appliqu\'e \`a $N_0$ et \`a $\delta/2$, il existe un ouvert $U$ dont l'adh\'erence est disjointe de ses $N_0$ premiers it\'er\'es  et un compact $V\subset U$  tel que tout point $x\in K\setminus  \bigcup_{p\in Per_{N_0}}W^s_{\delta}(p)$ poss\`ede un it\'er\'e  $f^n(x)\in V$, $n>0$  et tout point $y\in K\setminus \bigcup_{p\in Per_{N_0}}W^u_{\delta}(p)$ poss\`ede un it\'er\'e $f^{-n}(y)\in V$, $n>0$.
De plus on peut choisir $U$ tel  que les composantes de son adh\'erence  soient assez petites pour que tous leurs it\'er\'es d'ordre inf\'erieur \`a  $N_0$ soient inclus dans une carte de perturbation donn\'ee par le th\'eor\`eme~\ref{t.connecting}.

Soit $U_0$ une composante de $U$ et notons $V_0=U_0\cap V$. Par hypoth\`ese, $U_0$ est contenue dans une carte de perturbation $\varphi\colon U_0\to \RR^n$. Pavons $\RR^n$ par des cubes quadrill\'es de c\^ot\'e de m\^eme longueur $\ell$. On choisit $\ell$ de fa\c con que tout cube rencontrant $\varphi(V_0)$ soit inclus dans $\varphi(U_0)$.

On fait de m\^eme pour chacune des composantes de $U$ contenant un point de $V$. On obtient ainsi une famille finie $\cP_0$ de bo\^\i tes de perturbation, deux \`a deux disjointes (une bo\^\i te de perturbation est un cube ouvert), contenues dans $U$, et dont l'union des adh\'erences contient $V$ dans son int\'erieur. Notons $\Phi_0$ la famille des cartes de perturbation servant dans la construction de $\cP_0$.

En r\'ep\'etant cette construction pour $f^{2iN}(U),f^{2iN}(V)$, $i\in\{1,\dots,10d-1\}$, on construit de m\^eme des familles $\cP_i$ de  bo\^\i tes de perturbations contenues dans $f^{2iN}(U)$, et dont l'union des adh\'erences contient $f^{2iN}(V)$ dans son int\'erieur, et l'on note $\Phi_i$ la famille des cartes de perturbation correspondantes. 

On consid\`ere les familles de cubes $f^{-2iN}(\cP_i)$ contenues dans $U$. Pour tout $i$, l'union des adh\'erences des cubes de la famille $f^{-2iN}(\cP_i)$ contient $V$ dans son int\'erieur. Quitte \`a perturber l\'eg\`erement les $\Phi_i$ (en topologie $C^1$), on peut supposer que, pour toute famille de cubes choisis dans des familles deux \`a deux diff\'erentes $f^{-2iN}(\cP_i)$, leurs bords sont en position g\'en\'erale. En particulier un point de $V$ peut appartenir au bord de cubes d'au plus $d$ familles parmi les $f^{-2iN}(\cP_i)$. Puisqu'il y a au moins $5d$ familles de cubes, tout point de $V$ appartient \`a l'int\'erieur d'au moins $4d$ de ces cubes.

On remplace \`a pr\'esent dans $\RR^n$ chaque cube quadrill\'e par un cube de m\^eme centre, homoth\'etique de rapport $\rho<1$ tr\`es proche de $1$. On obtient ainsi des familles $\cP_{i,\rho}$ de bo\^\i tes de perturbation qui cette fois sont d'adh\'erences disjointes. Pour $\rho$ assez proche de $1$, tout point de $V$ appartient encore \`a l'int\'erieur d'un cube d'au moins $4d$ familles $f^{-2iN}(\cP_{i,\rho})$.

Par compacit\'e de $V$, on peut extraire, pour tout $i$,  une famille finie $\Ga_i$ de carreaux des bo\^\i tes de $\cP_{i,\rho}$ telle que  l'union $\Si_i$ de ces carreaux poss\`ede la propri\'et\'e suivante~: 
tout point de $V$ appartient \`a l'int\'erieur d'au moins $4d$ compacts  $f^{-2iN}(\Si_i)$. Appelons carreaux de $f^{-2iN}(\Si_i)$ l'image par $f^{-2iN}$ des carreaux de la famille $\Ga_i$.

Quitte \`a perturber \`a nouveau les $\Phi_i$, on peut supposer que, pour toute famille de carreaux choisis dans des familles deux \`a deux diff\'erentes $f^{-2iN}(\Ga_i)$, leurs bords sont en position g\'en\'erale. En particulier un point de $V$ peut appartenir au bord de carreaux d'au plus $d$ familles parmi les $f^{-2iN}(\Ga_i)$. Puisqu'il y a au moins $4d$ familles de carreaux, tout point de $V$ appartient \`a l'int\'erieur d'au moins un de ces carreaux.

Ainsi tout point de $V$ appartient \`a l'int\'erieur d'un carreau d'une famille $f^{-2iN}(\Ga_i)$.  Par cons\'equent, d'apr\`es le corollaire~\ref{c.coloriage2}, tout point de $K\setminus W^s_{\delta}(Per_{N_0})$ a un it\'er\'e positif dans l'int\'erieur d'un carreau d'une des familles $\Ga_i$. De m\^eme, tout point de $K\setminus f^{N_0}(W^u_{\delta}(Per_{N_0}))$ a un it\'er\'e n\'egatif dans $f^{N_0}(V)$ et donc dans un des carreaux d'une famille $\Ga_i$. En appliquant ce r\'esultat \`a un $\delta'$ assez petit pour que $f^{N_0}(W^u_{\delta'}(Per_{N_0}))\subset W^u_{\delta}(Per_{N_0})$, on obtient la propri\'et\'e voulue.

Remarquons que par construction, les supports des bo\^\i tes sont deux \`a deux disjoints. Plus pr\'ecis\'ement, pour toutes bo\^\i tes distinctes $B$ et $B'$ choisies dans l'union  des familles $\cP_{i,\rho}$ et tous $n,n'\in\{0,\dots, 2N-1\}$, les it\'er\'es $f^{n}(B)$ et $f^{n'}(B')$ sont disjoints. En particulier, pour tout point $x$ de $M$, il existe un entier $n\in\{0,\dots,N\}$ tel que l'it\'er\'e $f^n(x)$ n'appartienne \`a aucune bo\^\i te des familles $\cP_{i,\rho}$.
L'union des familles $\cP_{i,\rho}$ et l'union des familles $\Ga_i$ d\'efinissent les familles $\cB_0$ et $\cC_0$ respectivement.

Finalement, le compact $V$ a \'et\'e recouvert par l'int\'erieur des carreaux des familles $f^{-2iN}(\Ga_i)$. D'un recouvrement d'un compact par des ouverts, on peut choisir un compact inclus dans chacun de ces ouverts de fa\c con \`a obtenir un recouvrement du compact initial par l'int\'erieur de ces compacts. Ceci permet de choisir la famille $\cD_0$ annonc\'ees. 

\end{demo}

\subsection{Au voisinage des points p\'eriodiques}

\begin{prop}\label{p.local}
Soit $f$ un diff\'eomorphisme d'une vari\'et\'e compacte de dimension $d$, et $\cU$ voisinage de $f$ pour la topologie $C^1$. Soit $N$  l'entier associ\'e \`a l'ouvert $\cU$ par le th\'eor\`eme~\ref{t.connecting}.
Pour toute orbite p\'eriodique hyperbolique $\gamma$ de $f$, et tout voisinage $V$ de $\gamma$, il existe un voisinage $W$ de $\gamma$, deux familles finies de bo\^\i tes de perturbations $\cE$ et $\cS$ d'ordre $N$ pour $(f,\cU)$, deux familles finies de carreaux $\cC(\cE)$ et $\cC(\cS)$ des bo\^\i tes de $\cE$ et $\cS$ et deux familles finies $\cD(\cE)$ et $\cD(\cS)$ de compacts,  et un entier $n_0$ tels que~:
\begin{enumerate}
\item L'ouvert $V$ contient $\bar W$ et le support des bo\^\i tes de perturbation de $\cE$ et $\cS$.
\item Les supports des bo\^\i tes de perturbations de $\cE\cup\cS$ sont deux \`a deux disjoints.
\item Chaque compact de $\cD(\cE)$ est inclus dans l'int\'erieur d'un carreau $C\in\cC(\cE)$, et chaque compact de $\cD(\cS)$ dans l'int\'erieur d'un carreau $C\in\cC(\cS)$. Chaque carreau de $\cC(\cE)$ et $\cC(\cS)$ contient exactement un carreau de $\cC(\cE)\cup\cC(\cS)$.
\item Pour toute paire $D_e\in \cD(\cE)$ et $D_s\in \cD(\cS)$ de carreaux, il existe $n\in\{0,\dots,n_0\}$ tel que $f^n(D_e)\cap D_s\neq \emptyset$.
\item Pour tout point $z$ de $W\setminus W^s_{loc}(\gamma)$, il existe $n>0$ et un compact $D\in \cD(\cS)$ tel que $f^n(z)\in int(D)$ et $f^j(z)\in V$ pour tout $0\leq j \leq n$. De plus, si $f(z)\not\in W$, l'entier $n$ est major\'e par $n_0$.
\item Pour tout point $z$ de $W\setminus W^u_{loc}(\gamma)$, il existe $n>0$ et un compact $D\in \cD(\cE)$ tel que $f^{-n}(z)\in int(D)$ et $f^{-j}(z)\in V$ pour tout $0\leq j \leq n$. De plus, si $f^{-1}(z)\not \in W$, l'entier $n$ est major\'e par $n_0$.
\end{enumerate}
\end{prop}

\begin{demo}
De m\^eme que dans l'\'enonc\'e du corollaire~\ref{c.coloriage}, on pose $N_0=10 dN$.

Soit $\gamma$ une orbite p\'eriodique hyperbolique et $V$ un voisinage de $\gamma$.  Nous choisissons une vari\'et\'e stable locale $W^s_{loc}(\gamma)$ dont tous les it\'er\'es positifs sont contenus dans $V$. Il existe un ouvert $U^s\subset V$ tel que $U^s$ et ses $N_0$ premiers it\'er\'es  sont deux \`a deux disjoints et contenus dans $V$ et tel que $U^s\cap W^s_{loc}(\gamma)$ rencontre toute orbite de $W^s(\gamma)\setminus\{\gamma\}$. De m\^eme, on choisit une vari\'et\'e instable locale $W_{loc}^u(\gamma)$ dont tous les it\'er\'es n\'egatifs sont contenus dans $V$ et un ouvert $U^u\subset V$ tel que $U^u$ et ses $N_0$ premiers it\'er\'es sont deux \`a deux disjoints, disjoints de l'union $\bigcup_{j=0}^{N_0}f^j(U^s)$, contenus dans $V$ et tels que toute orbite de $W^u(\gamma)\setminus\{\gamma\}$ rencontre $f^{N_0}(U^u)\cap W^u_{loc}(\gamma)$.

Comme dans la preuve du corollaire~\ref{c.coloriage}, on peut choisir  une famille finie $\cE$ de  bo\^\i tes de perturbation dont les supports sont deux \`a deux disjoints et contenus dans l'int\'erieur de $\bigcup_{j=0}^{N_0}f^j(U^s)$, et ayant la propri\'et\'e suivante~: il existe une famille finie $\cC(\cE)$ de carreaux  de l'ensemble de ces bo\^\i tes, et pour chaque carreau $C$ un compact $D$ inclus dans l'int\'erieur de $C$ tels que toute orbite de $W^s(\gamma)\setminus \gamma$ rencontre l'intersection de $W_{loc}^s(\gamma)$ avec l'int\'erieur d'un des compacts $D$ (le compact $D$ d\'epend de l'orbite). Nous noterons $\cD(\cE)$ l'ensemble de ces compacts. 

De m\^eme, on peut choisir une famille finie $\cS$ de bo\^\i tes de perturbation  de supports  deux \`a deux disjoints et contenus dans l'int\'erieur de $\bigcup_{j=0}^{N_0}f^j(U^u)$,   une famille finie $\cC(\cS)$ de carreaux de ces bo\^\i tes  et une famille $\cD(\cS)$ de compacts inclus dans l'int\'erieur des carreaux de $\cC(\cS)$ (un compact par carreau) tel que toute orbite de $W^u(\gamma)\setminus \gamma$ rencontre l'intersection de $W_{loc}^u(\gamma)$ avec l'int\'erieur d'un des compacts $D\in\cD(\cS)$. 

Soient $D_e\in \cD(\cE)$ et $D_s\in\cD(\cS)$ deux compacts comme ci-dessus. Remarquons que $D_e$ et $D_s$ contiennent dans leur int\'erieur un point de $W^s(\gamma)$ et de $W^u(\gamma)$, respectivement. Le $\lambda$-lemma montre qu'il existe un entier $n=n(D_e,D_s)>0$ tel que $f^n(D_e)\cap D_s\neq \emptyset$.

Notons $D(\cE)$ et $D(\cS)$ l'union des \'el\'ements de $\cD(\cE)$ et $\cD(\cS)$. Pour tout $n\geq 0$, l'ouvert $D_n(\cE)=f^n(int(D(\cE)))\cap \bigcap_{j=0}^{n-1}f^j(V)$, est l'image par $f^n$ des points de $int(D(\cE))$ dont les $n$ premiers it\'er\'es restent dans $V$. L'ensemble $W_e=W^u_{loc}(\gamma)\cup\bigcup_{n\geq 0} D_n(\cE)$ contient un voisinage de l'orbite $\gamma$. Par construction, tout point $x$ de $W_e\setminus W^u_{loc}(\gamma)$ poss\`ede un it\'er\'e n\'egatif $f^{-n}(x)$ dans l'int\'erieur d'un compact $D\in\cD(\cE)$ tel que $f^{-i}(x)$ reste dans $V$ pour tout $i\in\{0,\dots,n\}$.

De m\^eme, l'union $W_s$ de $W^s_{loc}(\gamma)$ et, pour tout $n$, de l'image par $f^{-n}$ des points de $int(D(\cS))$ dont les $n$ premiers it\'er\'es n\'egatifs restent dans $V$ contient un voisinage de l'orbite $\gamma$. On choisit un voisinage ouvert $W_0$ de $\gamma$ contenu dans l'intersection $W_e\cap W_s$.

Remarquons que l'orbite positive d'un  point $x$ de $W_0$   passe par l'int\'erieur d'un compact $D\in \cD(\cS)$ avant de sortir de $V$ si $x$ n'est pas sur la vari\'et\'e stable locale de $\gamma$, et son orbite n\'egative passe par l'int\'erieur d'un compact $D\in \cD(\cE)$ avant de sortir de $V$ si $x$ n'est pas sur la vari\'et\'e instable locale de $\gamma$. On pourrait choisir $W=W_0$ et conclure la preuve de la proposition si le temps pour atteindre $int(D(\cS))$ ou $int(D(\cE))$ \'etait uniform\'ement born\'e pour tout $x\in W_0$ tel que $f(x)\notin W_0$ ou $f^{-1}(x)\notin W_0$, respectivement. Nous allons choisir $W\subset W_0$ ayant cette derni\`ere propri\'et\'e.

Il existe des domaines fondamentaux $\Ga_e$ et $\Ga_s$ inclus dans $W_0$ des vari\'et\'es stables et instables locales de $\gamma$, tels que le temps pour atteindre $int(D(\cE))$ ou $int(D(\cS))$ est uniform\'ement born\'e par un entier $n_1$ pour tout $x\in \overline{\Ga_e}$ ou $x\in\overline{\Ga_s}$, respectivement. De plus cette borne uniforme reste valide sur des voisinages tubulaires $\De_e$ et $\De_s$ des adh\'erences $\overline{\Ga_e}$ et $\overline{\Ga_s}$.  On peut choisir un voisinage ouvert $W$ de $\gamma$ de fa\c con que tout point $x\in W$ tel que $f(x)\notin W$ ou $f^{-1}(x)\notin W$ appartient \`a de tels voisinages tubulaires $\De_s$ ou $\De_e$, respectivement.

Pour conclure la preuve de la proposition, il suffit de choisir $n_0$  majorant $n_1$ et tous les entiers $n(D_e,D_s)$ d\'efinis ci-dessus.

\end{demo}

\subsection{Preuve du th\'eor\`eme~\ref{t.connect}~: pr\'eparation}\label{ss.preparation}

Soit $f$ un diff\'eomorphisme d'une vari\'et\'e compacte $M$, dont toutes les orbites p\'eriodiques sont hyperboliques. Soit $\cU$ un voisinage de $f$ pour la topologie $C^1$. Soient $x$ et $y$ deux points de $M$ tels que $x\dashv y$.

\begin{rema}\label{r.u}
{\em \'Etant donn\'e un voisinage $\cU_0$ de $f$  il existe  un ouvert $\cU\subset \cU_0$ ayant la propri\'et\'e suivante~:

{\em Soient $V_1,\dots,V_r$ des ouvert disjoints de $M$ et $g_1,\dots, g_r\in \cU$ des diff\'eomorphismes tels que pour tout $i\in\{1,\dots,r\}$, $g_i$ co\"\i ncide avec   $f$ hors de $V_i$. Alors le diff\'eomorphisme $g$, \'egal \`a $g_i$ sur $V_i$ et \`a $f$ hors de $\bigcup_i V_i$, appartient \'egalement \`a $\cU$.}}
\end{rema}
Quitte \`a restreindre le  voisinage $\cU$ de $f$, on supposera d\'esormais que $\cU$ v\'erifie la propri\'et\'e de la remarque ci-dessus. On consid\`ere l'entier  $N$  associ\'e \`a $\cU$ par le th\'eor\`eme~\ref{t.connecting} et on pose $N_0=10dN$.   On choisit une constante $\delta>0$ telle que les vari\'et\'es invariantes locales  $W^s_\delta(z)$ et $W^u_\delta(z')$ de points $z,z'$ de $Per_{N_0}(f)$ ne se coupent qu'aux points de $Per_{N_0}(f)$.

Afin de commencer la preuve du th\'eor\`eme, nous allons maintenant fixer une famille de bo\^\i tes de perturbations associ\'ees \`a l'ouvert $\cU$ et aux points $x$ et $y$. Ceci se fait en deux \'etapes~: une premi\`ere famille de  bo\^\i tes est choisie \`a l'aide du corollaire~\ref{c.coloriage}. Cette famille sera enrichie de bo\^\i tes, au voisinage des orbites p\'eriodiques de basse p\'eriode, \`a l'aide de la proposition~\ref{p.local}. 

\vskip 2mm

 Nous allons consid\'erer quelques configurations particuli\`eres des points $x$ et $y$  et montrer comment elles se ram\`enent au cas g\'en\'eral~: 
\begin{rema} \label{r.casparticulier}
{\em
\begin{enumerate}

\item {\em Si $y$ est un it\'er\'e strictement positif de $x$}~: la conclusion du th\'eor\`eme est  alors d\'ej\`a v\'erifi\'ee. 

\item {\em Si $x$ ou $y$ appartiennent pas \`a $Per_{N_0}(f)$}~: si $x\in Per_{N_0}(f)$ et si $y$ n'est pas sur l'orbite de $x$  (voir l'item (1) ci-dessus) alors il existe un point $x'$ (arbitrairement proche de $x$) sur la vari\'et\'e instable locale de $x$ tel que $x'\dashv y$~: le point $x'$ est un point d'accumulation de points de $\varepsilon$-pseudo-orbites joignant $x$ \`a $y$ quand $\varepsilon$ tend vers $0$. Pour montrer la conclusion du th\'eor\`eme il suffit dans ce cas de la montrer pour les points $x'$ et $y$. En effet, s'il existe une $C^1$-arbitrairement petite perturbation $g$ de $f$ telle que $y$ soit sur l'orbite positive pour $g$ de $x'$, par une conjugaison proche de l'identit\'e \`a support ne contenant pas le point $y$, on envoie $x'$ sur $x$, et maintenant l'orbite positive de $x$ passe par $y$. On traite de la m\^eme fa\c con le cas o\`u $y\in Per_{N_0}(f)$.
\end{enumerate}
}

\vskip 2mm
Par la suite nous supposerons d\'esormais que $x$ et $y$ n'appartiennent pas \`a $Per_{N_0}(f)$ et que $y$ n'est pas un point de l'orbite strictement positive de $x$. 

\end{rema}

On fixe une famille $\cB_0$ de bo\^\i tes de perturbations donn\'ees par le corollaire~\ref{c.coloriage}, appliqu\'e \`a toute la vari\'et\'e  au voisinage $\cU$ et \`a la constante $\delta$. On note $\cC_0$ et $\cD_0$ l'ensemble de carreaux $C$ et  l'ensemble de compacts $D$ contenus dans l'int\'erieur des carreaux $C$, annonc\'es par ce corollaire. 

Pour toute orbite p\'eriodique $\gamma$ de p\'eriode inf\'erieure \`a $N_0$ on choisit un voisinage $V(\gamma)$ de fa\c con que les $V(\gamma)$ soient deux \`a deux disjoints, soient chacun disjoint des supports de bo\^\i tes de perturbations de $\cB_0$  et ne contiennent pas les points $x$ et $y$ (c'est possible puisque ni $x$ ni $y$ n'appartient \`a $Per_{N_0}(f)$). Notons $\cE(\gamma)$ et $\cS(\gamma)$ les familles de bo\^\i tes de perturbations contenues dans $V(\gamma)$  obtenues \`a l'aide de la proposition~\ref{p.local}. En suivant encore les notation de cette proposition, on notera $\cC(\cE,\gamma)$ et $\cC(\cS,\gamma)$, (resp. $\cD(\cE,\gamma)$ et $\cD(\cS,\gamma)$) les familles de carreaux (resp. de compacts contenus dans l'int\'erieur des carreaux) des bo\^\i tes de perturbation de $\cE(\gamma)$ et $\cS(\gamma)$, respectivement~; finalement on note $n_0(\gamma)$ et $W(\gamma)$ l'entier et le voisinage de $\gamma$ annonc\'es par cette proposition. On notera $n_0$ un majorant des entiers $n_0(\gamma)$, pour $\gamma\in Per_{N_0}(f)$.

On note $\cB$ l'ensemble $\cB_0\cup\bigcup_{\gamma\in Per_{N_0}(f)} \left(\cE(\gamma)\cup\cS(\gamma)\right)$ de toutes les bo\^\i tes de perturbation ainsi obtenues. On note de m\^eme $\cC$ et $\cD$ l'ensemble de tous les carreaux et l'ensemble de tous les compacts contenus dans l'int\'erieur des carreaux, respectivement, ainsi obtenus.

\begin{rema} \label{r.hors}

 Si $x$ ou $y$ appartient au support d'une bo\^\i te de perturbation $B\in\cB${\em   (rappelons  que  l'on a suppos\'e que $y$ n'est pas un it\'er\'e strictement positif de $x$)~: remarquons que la bo\^\i te $B$ doit appartenir \`a la famille $\cB_0$ car $x$ et $y$ n'appartiennent pas aux ouverts $V(\gamma)$.

Puisque $y$ n'est pas un it\'er\'e strictement positif de $x$, pour tout it\'er\'e $f^i(x)$ et $f^j(y)$ on a encore $f^i(x)\dashv f^j(y)$. D'autre part, d'apr\`es le corollaire~\ref{c.coloriage},   on peut choisir $i$ et $j$ dans $\{0,\dots,N\}$ tels que $f^i(x)$ et $f^j(y)$ n'appartiennent pas \`a l'union des supports des $B_i$ (si $x=y$ on choisit bien s\^ur $i=j$). Pour montrer la conclusion du th\'eor\`eme il suffit dans ce cas de la montrer pour les points $f^i(x)$ et $f^j(y)$. En effet, s'il existe une $C^1$-arbitrairement petite perturbation $g$ de $f$ telle que $f^j(y)$ soit sur l'orbite positive pour $g$ de $f^i(x)$, par une conjugaison proche de l'identit\'e, on envoie $g^{-i}(f^i(x))$ sur $x$ et $g^{-j}(f^j(y))$ sur $y$, et maintenant l'orbite positive de $x$ passe par $y$.
}
\vskip 2mm

On supposera d\'esormais de plus que les points $x$ et $y$ n'appartiennent pas aux supports des bo\^\i tes de perturbation de $\cB$. Cette hypoth\`ese est sans perte de g\'en\'eralit\'e d'apr\`es la remarque ci-dessus. 
\end{rema}

\subsection{Retour des pseudo-orbites dans les carreaux}\label{s.regroupe}

Si une pseudo-orbite retourne en temps born\'e dans les compacts $D\in \cD$ des carreaux et si ses sauts sont suffisamment petits, on peut regrouper les sauts au moment des passages dans les carreaux, obtenant ainsi une pseudo-orbite qui pr\'eserve les quadrillages des bo\^\i tes de perturbation (voir la section~\ref{s.boites} et la remarque~\ref{r.tuboite}). Le lemme~\ref{l.singularites2} montre que si le retour dans les carreaux se fait en temps tr\`es long, l'orbite passe pr\`es des points p\'eriodiques de $Per_{N_0}$ et donc dans les carreaux des familles $\cE(\gamma)$ et $\cS(\gamma)$. On remplacera alors le segment de pseudo-orbite pr\`es d'une orbite p\'eriodique $\gamma$ par un morceau de vraie orbite entrant et sortant du voisinage de $\gamma$ par les m\^emes carreaux de $\cE(\gamma)$ et $\cS(\gamma)$ que la pseudo-orbite initiale.

\begin{lemm}\label{l.singularites2}

Avec les notation de la section~\ref{ss.preparation}, il existe $n_1>0$ et $\varepsilon_0>0$ avec la propri\'et\'e suivante~:

Pour toute $\varepsilon_0$-pseudo-orbite $x_0=x,\dots,x_k=y$, il existe une suite strictement croissante $t_0=0,t_1,\dots,t_r=k$ telle que, pour tout $s\notin\{ 0,r\}$, le point $x_{t_s}$ appartient \`a un compact $D_s$ de la famille $\cD$. De plus, pour tout $s\in\{0,\dots,r-1\}$, ou bien $t_{s+1}-t_s\leq n_1$, ou bien il existe une orbite p\'eriodique  $\gamma\in Per_{N_0}(f)$ telle que $D_s\in \cD(\cE,\gamma)$ et $D_{s+1}\in \cD(\cS,\gamma)$ (dans ce cas, $s\notin\{ 0,r-1\}$). 
\end{lemm}
\begin{demo}

Consid\'erons tout d'abord le cas des vraies orbites~: tout point $z$ a des it\'er\'es positifs et n\'egatifs dans l'union des ouverts $W(\gamma)$ et des int\'erieurs des compacts $D$ de $\cD$. Les temps de retour sont uniform\'ement born\'es par un entier $n_1$, par compacit\'e de $M$ . On choisit $n_1$ sup\'erieur \`a l'entier $n_0$ introduit \`a la fin de la section~\ref{ss.preparation} (de ce fait, $n_1$ majore tous les entiers $n_0(\gamma)$ pour chaque orbite p\'eriodique $\gamma$ de $Per_{N_0}(f)$).

Afin d'assurer que les pseudo-orbites visitent l'int\'erieur des compacts $D\in\cD$, nous introduisons des compacts plus petits $\tilde D$ inclus dans l'int\'erieur des compacts $D$, et visit\'es par les vraies orbites.

Par compacit\'e de $M$, tout point de $M$ retourne, en temps born\'e par $n_1$, dans des compacts $\tilde W(\gamma)$ et $\tilde D$ inclus dans  les ouverts $W(\gamma)$ et les int\'erieurs des compacts $D$.  En utilisant la d\'efinition de $n_0(\gamma)$ (voir la proposition~\ref{p.local}), on montre que l'on peut choisir la famille $\tilde D$ de fa\c con que, si $z\in W(\gamma)$ mais que $f(z)\notin W(\gamma)$, il existe un it\'er\'e positif inf\'erieur \`a $n_0(\gamma) < n_1$ de $z$ qui appartient \`a l'int\'erieur d'un  de ces compacts $\tilde D$ contenu dans l'int\'erieur d'un compact $D\in \cD(\cS,\gamma)$. De plus, le segment d'orbite joignant $z$ \`a $\tilde D$ est inclus dans $V(\gamma)$. De m\^eme, si $z\in W(\gamma)\setminus f(W(\gamma))$, son orbite n\'egative rencontre, en temps born\'e par $n_1$, l'int\'erieur d'un compact $\tilde D$ contenu dans un compact $D\in\cD(\cE,\gamma)$ avant de sortir de $V(\gamma)$.

On note $\eta$ l'infimum des distances entre un point de $\bigcup_{\gamma\in Per_{N_0}(f)} \tilde W(\gamma)\cup \bigcup_{D\in \cD} \tilde D$ et un point du compl\'ementaire de $\bigcup_{\gamma\in Per_{N_0}(f)} W(\gamma)\cup \bigcup_{D\in\cD} D$.

Il existe alors $\varepsilon_0>0$, tel que pour toute $\varepsilon_0$-pseudo-orbite $z_0,\dots,z_{n_1}$,  et tout $i\in\{0,\dots,n_1\}$, les distances $d(f^{i}(z_0),z_i)$ et $d(f^{-i}(z_i),z_0)$ sont inf\'erieures \`a $\eta$. Par ce qui pr\'ec\`ede, il existe $i\in\{0,\dots,n_1\}$, tel que $f^i(z_0)$ appartienne \`a un compact $\tilde W(\gamma)$ ou $\tilde D$. On en d\'eduit que le point $z_i$ appartient \`a $W(\gamma)$ ou \`a l'int\'erieur du compact $D$ correspondant.

Soit $x_0,\dots,x_k$ une $\varepsilon_0$-pseudo-orbite dont les extr\'emit\'es $x_0,x_k$ ne sont dans aucun $V(\gamma)$. Si un point $x_j$ appartient \`a un $W(\gamma)$ et si $f(x_{j})$ n'appartient pas \`a $W(\gamma)$, il existe $x_{j'}$ avec $0\leq j'-j \leq n_1$ tel que $x_{j'}$ appartienne \`a l'int\'erieur d'un compact $D\in \cD(\cS,\gamma)$. En effet,  par d\'efinition de $n_0(\gamma)<n_1$, il existe un it\'er\'e $f^{j'-j}(x_j)$ de $x_j$ qui appartient \`a l'int\'erieur de $\tilde D$. De m\^eme, si un point $x_j$ appartient \`a un $W(\gamma)$ et si $f^{-1}(x_{j})$ n'appartient pas \`a $W(\gamma)$, il existe $x_{j'}$ avec $-n_1\leq j'-j \leq 0$ tel que $x_{j'}$ appartienne \`a l'int\'erieur d'un compact $D\in \cD(\cE,\gamma)$.

\`A tout point $x_j$ appartenant \`a l'un des $W(\gamma)$, on associe un segment de la pseudo-orbite $x_{e(j)},\dots,x_{s(j)}$ contenant $x_j$ de la fa\c con suivante~: consid\'erons le segment maximal de la pseudo-orbite contenant $x_j$ et inclus dans $V(\gamma)$. Alors, $e(j)$ est le plus petit indice dans ce segment tel que $x_{e(j)}$ appartient \`a l'int\'erieur d'un compact $D\in \cD(\cE,\gamma)$ et $s(j)$ est le plus grand indice dans ce segment tel que $x_{s(j)}$ appartient \`a l'int\'erieur d'un compact $D\in \cD(\cS,\gamma)$.

\begin{affi}
Pour tout $j$ tel que $x_j$ appartient \`a l'union des $W(\gamma)$, le segment $\{e(j),\dots,s(j)\}$ est bien d\'efini et
contient $j$. Deux segments $\{e(j),\dots,s(j)\}$ et $\{e(j'),\dots,s(j')\}$ sont disjoints ou confondus.
\end{affi}
\begin{demo}
D'apr\`es ce qui pr\'ec\`ede, il existe un segment de pseudo-orbite $x_{j_0},\dots,x_j$, avec $j_0<j$, inclus dans $V(\gamma)$ et tel que $x_{j_0}$ appartienne \`a l'int\'erieur d'un compact $D\in \cD(\cE,\gamma)$. Ceci montre que $e(j)$ est bien d\'efini et inf\'erieur \`a $j$. De m\^eme $s(j)$ est bien d\'efini et sup\'erieur \`a $j$.

Soient $I=\{e(j),\dots,s(j)\}$ et $I'=\{e(j'),\dots,s(j')\}$ deux segments se rencontrant en un point. Les $V(\gamma)$ \'etant deux \`a deux disjoints, les points $x_j$, $j\in I$, et $x_{j'}$, $j'\in I'$ appartiennent au m\^eme $V(\gamma)$. On obtient $I=I'$ par maximalit\'e.
\end{demo}

Soient $t_1<t_2<\cdots<t_{r-1}$ l'ensemble des indices $j$ pour lesquels~:
\begin{itemize}
\item soit il existe $j'$ tel que $j=e(j')$ ou $j=s(j')$,
\item soit $x_j$ appartient \`a l'int\'erieur d'un compact $D\in \cD$ et $j$ n'appartient \`a aucun segment $\{e(j'),\dots,s(j')\}$. 
\end{itemize}
En cons\'equence, pour tout $s\in\{1,\dots,r-2\}$, ou bien le segment de pseudo-orbite $\{x_{t_s},\dots,x_{t_{s+1}}\}$ rencontre un ouvert $W(\gamma)$ et dans ce cas il existe $j'$ tel que $t_s=e(j')$ et $t_{s+1}=s(j')$, ou bien ce segment d'orbite ne rencontre aucun $W(\gamma)$ et dans ce cas le segment de pseudo-orbite $\{x_{t_s+1},\dots,x_{t_{s+1}-1}\}$ ne rencontre aucun int\'erieur de compact $D\in\cD$.

Nous devons montrer que tout intervalle de temps $\{t_s,\dots,t_{s+1}\}$ qui n'est pas de la forme $\{e(j),\dots,s(j)\}$ est de longueur $t_{s+1}-t_s\leq n_1$. Le segment de pseudo-orbite correspondant n'a, par construction, de point dans aucun ouvert $W(\gamma)$. On consid\`ere le point $x_{t_s}$. Par d\'efinition de $n_1$, ce point poss\`ede un premier it\'er\'e positif $f^{\ell}(x_{t_s})$, avec $0\leq \ell\leq n_1$, appartenant \`a l'un des $\tilde W(\gamma)$ ou \`a l'un des $\tilde D$. Par le choix de $\varepsilon_0$, le point $x_{t_s+\ell}$ appartient \`a $W(\gamma)$ ou \`a l'int\'erieur d'un compact $D\in \cD$.  

Montrons que $t_{s+1}\leq t_s+\ell$ ce qui conclura puisque $\ell\leq n_1$~: si $x_{t_s+\ell}$ appartient \`a un $W(\gamma)$, on a remarqu\'e que la suite $\{x_{t_s},\dots,x_{t_{s+1}} \}$ ne contient pas de point dans $W(\gamma)$ et donc $t_{s+1}<t_s+\ell$. Dans l'autre cas, $x_{t_s+\ell}$ appartient \`a l'int\'erieur d'un compact $D$~; nous avons vu que le segment de pseudo-orbite $\{x_{t_s+1},\dots,x_{t_{s+1}-1}\}$ ne rencontre pas l'int\'erieur de $D$, ce qui implique que $t_{s+1}\leq t_s+\ell$.

Posons \`a pr\'esent $t_0=0$ et $t_r=k$.

Le point $x_0$ n'appartient \`a aucun $V(\gamma)$ et il existe $\ell\leq n_1$ tel que $x_\ell$ appartient \`a l'int\'erieur d'un compact $D\in \cD$ ou \`a un $W(\gamma)$. Dans le premier cas, $t_1\leq \ell$. Dans le second cas, $0<e(\ell)<\ell$ et $t_1\leq e(\ell)$.

Dans tous les cas, $t_1\leq n_1$. On montre de m\^eme que $t_{r}-t_{r-1}\leq n_1$.

\end{demo}

\subsection{Regroupement des sauts dans les carreaux } \label{s.regroupeencore}

Notons $\eta_1$ l'infimum des distances entre un point d'un compact $D\in\cD$ et le bord du carreau $C\in\cC$ qui contient $D$. Il  existe $\varepsilon_1\in]0,\varepsilon_0[$ tel que pour toute $\varepsilon_1$-pseudo-orbite $x_0,\dots,x_m$ de longueur $m$ inf\'erieure \`a $n_1$, les distances $d\left(f^{-m}(x_m),x_0\right)$ et $d\left(f^m(x_0),x_m\right)$ sont plus  petites que $\eta_1$.

Consid\'erons \`a pr\'esent  une $\varepsilon_1$-pseudo-orbite $x_0=x,\dots,x_k=y$ joignant $x$ \`a $y$. Comme $\varepsilon_1$ \`a \'et\'e choisi inf\'erieur \`a $\varepsilon_0$, le lemme~\ref{l.singularites2} lui associe une suite    $t_0=0<t_1<\cdots<t_r=k$. Par le choix de $\varepsilon_1$ et $\eta_1$, et comme les points $x_{t_i}$, $i\in\{1,\dots,r-1\}$, appartiennent tous \`a un compact $D\in\cD$, on obtient~:

\begin{lemm}\label{l.carreau} Sous ces hypoth\`eses, pour tout $s\in\{1,\dots,r-1\}$ tel que $t_{s+1}-t_s\leq n_1$, les points $f^{t_s-t_{s+1}}(x_{t_{s+1}})$ et $x_{t_s}$ sont dans l'int\'erieur d'un m\^eme carreau $C\in\cC$~; de m\^eme, pour tout $s\in\{0,\dots,r-2\}$ tel que $t_{s+1}-t_s\leq n_1$, les points $f^{t_{s+1}-t_{s}}(x_{t_{s}})$ et $x_{t_{s+1}}$ sont dans l'int\'erieur d'un m\^eme carreau $C\in\cC$.
\end{lemm}

Rappelons que les temps $t_s$ correspondent \`a des retours de la pseudo-orbite dans les compacts $D$, et que si un segment d'orbite entre deux temps $t_s,t_{s+1}$ cons\'ecutifs est de longueur plus grand que $n_1$, il est contenu dans un voisinage $V(\gamma)$ d'une orbite p\'eriodique $\gamma$ et \`a extr\'emit\'es dans les carreaux $\cC(\cE,\gamma)$ et $\cC(\cS,\gamma)$. 
Nous allons traiter diff\'eremment les segments d'orbite de longueur inf\'erieure \`a $n_1$ et ceux contenus dans un voisinage $V(\gamma)$, pour obtenir finalement une pseudo-orbite $y_0=x,\dots,y_\ell=y$ qui pr\'eserve le quadrillage des bo\^\i tes de perturbation.

Pour les segments $\{x_{t_s},\dots,x_{t_{s+1}}\}$ de longueur inf\'erieur \`a $n_1$, nous modifions la pseudo-orbite $x_0,\dots,x_k$ et d\'efinissons une nouvelle pseudo-orbite $\tilde x_0=x,\dots,\tilde x_k=y$  de la fa\c con suivante~:

\begin{itemize}
\item Pour tout $i\in\{0,\dots,t_1\}$ on pose $\tilde x_i= f^i(x_0)$. Remarquons que, comme $t_1\leq n_1$, le point $\tilde x_{t_1}$ appartient \`a l'int\'erieur du m\^eme carreau que $x_{t_1}$. 
\item Pour tout $s\in\{1,\dots,r-1\}$ tel que $t_{s+1}-t_s > n_1$ et tout $i\in\{t_s+1,\dots,t_{s+1}\}$, on pose $\tilde x_i=x_i$.
\item Pour tout $s\in\{1,\dots,r-1\}$ tel que $t_{s+1}-t_s \leq n_1$ et tout $i\in\{t_s+1,\dots,t_{s+1}\}$, on pose $\tilde x_i=f^{i-t_{s+1}}(x_{t_{s+1}})$.
\end{itemize}

Remarquons que, pour tout $s\in\{0,\dots,r-1\}$ tel que $t_{s+1}-t_s\leq n_1$, alors le segment $\tilde x_{t_s},\dots,\tilde x_{t_{s+1}}$ est une pseudo-orbite n'ayant qu'un seul saut~: entre $\tilde x_{t_s}$ et $\tilde x_{t_s+1}$. De plus ce saut pr\'eserve le quadrillage~: en effet, d'apr\`es le lemme~\ref{l.carreau}, $\tilde x_{t_s}$ et $f^{-1}(\tilde x_{t_s+1})$ appartiennent tous deux \`a l'int\'erieur du carreau qui contient $x_{t_s}$.

 Nous pouvons enfin construire la suite $y_0=x,y_1,\dots, y_\ell=y$ annonc\'ee~:

Pour tout $s$ tel que $t_{s+1}-t_s >n_1$, alors il existe une orbite p\'eriodique $\gamma\subset Per_{N_0}(f)$ telle que les points $x_{t_s}$ et $x_{t_{s+1}}$ appartiennent \`a des carreaux $C_e\in\cC(\cE,\gamma)$ et $C_s\in\cC(\cS,\gamma)$, respectivement (d'apr\`es le lemme~\ref{l.singularites2}). Les points  $x_{t_s}$ et $x_{t_{s+1}}$ appartiennent \`a ces m\^emes carreaux, par construction. D'apr\`es la proposition~\ref{p.local}, il existe un point $z$ de l'int\'erieur de $C_e$
et un entier $n\leq n_1$ tel que $f^n(z)$ appartienne \`a l'int\'erieur de $C_s$. 

On remplace alors le segment de pseudo-orbite $\{\tilde x_{t_s+1},\dots,\tilde x_{t_{s+1}}\}$ par le segment de vraie orbite 
$\{f(z),\dots,f^n(z)\}$. 

La suite obtenue en appliquant ce proc\'ed\'e \`a chaque segment $\{\tilde x_{t_s+1},\dots,\tilde x_{t_{s+1}}\}$ de longueur sup\'erieure \`a $n_1$ est not\'ee  $y_0=x,y_1,\dots, y_\ell=y$. On notera $\tau_0=0,\tau_1,\dots,\tau_r$ les temps de cette suite correspondants aux temps $t_0,\dots,t_r$ de la suite $\tilde x_i$~: formellement, la suite $\{y_0,\dots,y_\ell\}$ a \'et\'e obtenue en concat\'enant des segments de la forme $\{y_{\tau_0},\dots,y_{\tau_1}\}$ ou $\{y_{\tau_s+1},\dots,y_{\tau_{s+1}}\}$ qui remplacent les segments de la forme $\{\tilde x_{t_0},\dots,\tilde x_{t_1}\}$ ou $\{\tilde x_{t_s+1},\dots,\tilde x_{t_{s+1}}\}$.

\begin{lemme}
La suite $y_0=x,\dots,y_\ell=y$  d\'efinie ci-dessus est une pseudo-orbite pr\'eservant le quadrillage de toutes les bo\^\i tes $B$ de $\cB$, et n'ayant pas de saut hors des bo\^\i tes de $\cB$.
\end{lemme}
\begin{demo} Nous devons v\'erifier que pour tout $i\in\{0,\dots,\ell-1\}$, si $f(y_i)\neq y_{i+1}$ alors $y_i$ et $f^{-1}(y_{i+1})$ appartiennent \`a l'int\'erieur d'un m\^eme carreau $C\in \cC$. Par construction, ceci n'arrive que lorsque $i$ est de la forme $\tau_s$, $s\in\{1,\dots,r-1\}$. De nouveau par construction, les points $y_{\tau_s}$ et $f^{-1}(y_{\tau_{s}+1})$ appartiennent \`a l'int\'erieur du m\^eme carreau que $x_{t_s}$ (et aussi que $\tilde x_{t_s}$), ce qui conclut. 
\end{demo}

En r\'esum\'e, nous avons montr\'e~:
\begin{prop} \label{p.regroupe}
Il existe $\varepsilon_1>0$ tel que, pour tous points $x,y$ hors des voisinages $V(\gamma)$ et des supports des bo\^\i tes de perturbations $B\in\cB$ on ait la propri\'et\'e suivante~:

 Si $x\dashv_{\varepsilon_1} y$,  alors il existe une pseudo-orbite pr\'eservant le quadrillage de toute les bo\^\i tes de $\cB$ et n'ayant pas de sauts hors des bo\^\i tes de $\cB$ joignant $x$ \`a $y$.
\end{prop}

\subsection{Fin de la preuve du th\'eor\`eme~\ref{t.connect}}\label{s.fin}
Nous pouvons maintenant terminer la preuve du th\'eor\`eme~\ref{t.connect}~: 
soient $x$ et $y$ deux points tels que $x\dashv y$. D'apr\`es la section~\ref{ss.preparation}, on peut supposer que ces points n'appartiennent pas \`a $Per_{N_0}(f)$ et que $y$ n'est pas sur l'orbite positive de $x$. Soit $\cU$ un $C^1$ voisinage de $f$, que l'on peut choisir comme dans la remarque~\ref{r.u}. On note $\cB=\{B_1,\dots,B_m\}$ la famille de bo\^\i tes de perturbation que l'on a associ\'ee, dans la section~\ref{ss.preparation}, \`a  $\cU$ et \`a ces deux points. D'apr\`es la remarque~\ref{r.hors}, on peut supposer que les points $x$ et $y$ sont situ\'es hors des bo\^\i tes $B\in \cB$ (ils sont par construction hors des voisinages $V(\gamma)$).

D'apr\`es la proposition~\ref{p.regroupe}, il existe alors une pseudo-orbite $y_0=x,\dots,y_\ell=y$ pr\'eservant le quadrillage de toutes les bo\^\i tes de $\cB$ et n'ayant pas de sauts hors des bo\^\i tes de $\cB$.

Par le choix de l'ouvert $\cU$ et puisque les bo\^\i tes de $\cB$ sont de supports disjoints, un dif\-f\'e\-o\-mor\-phis\-me obtenu en composant des perturbations dans $\cU$ \`a support dans chacune des bo\^\i tes $B\in\cB$ est encore un diff\'eomorphisme de $\cU$.

D'apr\`es le lemme~\ref{l.tuboite} appliqu\'e  \`a la suite $y_0,\dots,y_\ell$ et \`a la bo\^\i te $B_1$, on obtient un dif\-f\'e\-o\-mor\-phis\-me $f_1\in \cU$ et une pseudo-orbite de $f_1$ joignant les points $x$ et $y$, pr\'eservant le quadrillage de toutes les bo\^\i tes de $\{B_2,\dots,B_m\}$ et n'ayant pas de sauts hors des bo\^\i tes de $\{B_2,\dots,B_m\}$.

On applique \`a nouveau ce proc\'ed\'e   successivement \`a  chaque bo\^\i te de $\{B_2,\dots,B_m\}$. \`A la $i^{\mbox{\`eme}}$-\'etape, on associe une perturbation de $f$ dans $\cU$, respectant les quadrillages de toutes les bo\^\i tes et n'ayant aucun saut hors des bo\^\i tes $B_{i+1},\dots,B_m$. La $m^{\mbox{\`eme}}$ \'etape cr\'ee un diff\'eomorphisme $g\in \cU$ tel que l'orbite positive de $x$ passe par $y$, ce qui conclut la d\'emonstration.
\section{Cons\'equences du lemme de connexion}\label{s.consequences}
\subsection{Preuve du th\'eor\`eme~\ref{t.equivalence}~: argument de g\'en\'ericit\'e}\label{s.equivalence}

Soit $\cB$ une base d\'enombrable d'ouverts de $M$.

Pour tous $U,V\in\cB$ on note $\cO(U,V)$ l'ensemble des diff\'eomorphismes tels que, ou bien, de fa\c con robuste, il existe un it\'er\'e positif de $U$ qui rencontre $V$, ou bien, de fa\c con robuste, tout it\'er\'e positif de $U$ est disjoint de $V$. Plus pr\'ecis\'ement, $\cO(U,V)$ est l'union $\cO_\infty(U,V)\cup \bigcup_0^\infty \cO_n(U,V)$, o\`u $\cO_n(U,V)$ est l'ensemble des diff\'eomorphismes pour lequel $f^n(U)$ rencontre $V$, et o\`u $\cO_\infty(U,V)$ est l'int\'erieur du compl\'ementaire de $\bigcup_0^\infty \cO_n(U,V)$.

Chacun des $\cO_n(U,V)$ est ouvert car $U$ et $V$ sont ouverts et $\cO_\infty(U,V)$ est ouvert par d\'efinition. On en d\'eduit que $\cO(U,V)$ est un ouvert de $\diff^1(M)$. Il est dense par construction.

Notons $\cG_0$ l'ensemble des diff\'eomorphismes Kupka-Smale~: toutes les orbites p\'eriodiques sont hyperboliques et leurs vari\'et\'es invariantes sont transverses. Par le th\'eor\`eme de Kupka-Smale, $\cG_0$ est un $G_\delta$ dense de $\diff^1(M)$. Soit $\cG_1$ l'intersection $\cG_0\cap\bigcap_{U,V\in \cB}\cO(U,V)$. C'est un $G_\delta$ dense de $\diff^1(M)$ d'apr\`es le th\'eor\`eme de Baire.

Soit $f$ un diff\'eomorphisme dans $\cG_1$, et soit $(x,y)$ une paire de points tels que $x\dashv_f y$. Soient $U\in \cB$ et $V\in \cB$ des voisinages de $x$ et $y$, respectivement. D'apr\`es le th\'eor\`eme~\ref{t.connect}  tout voisinage de $f$ contient un diff\'eomorphisme $g$ pour lequel un it\'er\'e positif de $U$ rencontre $V$.  Par d\'efinition de $\cG_1$, $f$ appartient \`a $\cO(U,V)$. L'existence de $g$ arbitrairement proche de $f$ interdit que $f$ appartienne \`a l'ouvert $\cO_\infty(U,V)$. Il existe donc un it\'er\'e positif de $U$ par $f$ qui rencontre $V$. On a montr\'e $x\prec_f y$, ce qui conclut la preuve du th\'eor\`eme~\ref{t.equivalence}.

\subsection{Cons\'equences  du th\'eor\`eme~\ref{t.equivalence}}\label{s.consequences2}

Le corollaire~\ref{c.generic} est un cons\'equence directe du th\'eor\`eme~\ref{t.equivalence}~: soit $\cG_1$ le $G_\delta$ dense de $\diff^1(M)$ construit \`a la section~\ref{s.equivalence}, sur lequelle les relations $\dashv$ et $\prec$ co\"\i ncident. $\Om(f)$ est l'ensemble des points $x$ tels que $x\prec x$ donc co\"\i ncide avec l'ensemble des points $x$ tels que $x\dashv x$ qui est $\cR(f)$.

La premi\`ere partie du corollaire~\ref{c.homocline} se d\'eduit du th\'eor\`eme~\ref{t.equivalence} et du lemme (sans doute classique) suivant~:
\begin{lemm}\label{l.dashv1} Si $f$ est un diff\'eomorphisme d'une vari\'et\'e connexe $M$ tel que $\Om(f)=M$ alors, pour tous $x,y\in M$ on a $x\dashv y$. 
\end{lemm}
\begin{demo} Pour tout $\varepsilon>0$ et tous $x,y\in M$ choisissons une suite de points $x_0=x, x_1,\dots,x_k=y$ tels que $d(x_i,x_{i+1})\leq \frac\varepsilon2$, pour tout $i$. Pour tout $i$ il existe une $\varepsilon$ pseudo-orbite joignant $x_i$ \`a $x_{i+1}$~: en effet, comme $x_i$ est non-errant il existe un point $\tilde x_i$ arbitrairement proche de $x_i$ ayant un retour $f^{n_i}(\tilde x_i)$ lui aussi proche de $x_i$~; la suite $x_i, f(\tilde x_i),f^2(\tilde x_i),\dots,f^{n_i-1}(\tilde x_i), x_{i+1}$ est alors une $\varepsilon$-pseudo-orbite. En mettant bout \`a bout ces pseudo-orbites on obtient une $\varepsilon$-pseudo-orbite joignant $x$ \`a $y$. Donc $x\dashv_\varepsilon y$ pour tout $\varepsilon>0$ et finalement $x\dashv y$.
\end{demo}

Si $f$ appartient \`a la partie r\'esiduelle $\cG_1$ de $\diff^1(M)$ annonc\'ee par le th\'eor\`eme~\ref{t.equivalence} et si $\Om(f)=M$, alors pour tous points $x,y\in M$ on a $x\dashv y$ et donc $x\prec y$. On en d\'eduit que pour tout couple d'ouvert $U,V$ de $M$ il existe un it\'er\'e positif de $U$ qui rencontre $V$, c'est \`a dire que $f$ est transitif~: ceci montre la premi\`ere partie du corollaire~\ref{c.homocline}.

Il existe d'autre part, comme cons\'equence du lemme de fermeture de Pugh, une partie r\'esiduelle $\cG_2$ de $\diff^1(M)$ sur laquelle l'ensemble non-errant $\Om(f)$, $f\in\cG_2$ est l'adh\'erence de l'ensemble des points p\'eriodiques hyperboliques de $f$. D'autre part, en cons\'equence du lemme de connexion d'Hayashi, \cite{BD2} montre qu'il existe une partie r\'esiduelle $\cG_3$ de $\diff^1(M)$ sur laquelle deux points p\'eriodiques hyperboliques appartiennent \`a un m\^eme ensemble transitif si et seulement s'ils ont m\^eme classe homocline. En cons\'equence, $\cG_1\cap\cG_2\cap\cG_3$ est une partie r\'esiduelle de $\diff^1(M)$ et, si $f\in\cG_1\cap\cG_2\cap\cG_3$ v\'erifie que $\Om(f)=M$, alors il est transitif et poss\`ede au moins un point p\'eriodique $p$ dont la classe homocline doit co\"\i ncider avec celle de tous les autres points p\'eriodiques, mais ceux-ci \'etant denses dans $M$, on obtient que $H(p,f)=M$. Ceci conclut la preuve du corollaire~\ref{c.homocline}.

\subsection{D\'ecomposition en pi\`eces \'el\'ementaires}

Le corollaire~\ref{c.coincide} d\'ecoule directement du th\'eor\`eme~\ref{t.equivalence} et du fait que, quand la relation $\prec$ est transitive (par exemple sur la partie r\'esiduelle $\cG_1$ introduite ci-dessus), les ensembles faiblement transitifs maximaux sont les classes d'\'equivalence de la relation  d'\'equivalence obtenue sur $\Om(f)$ en  sym\'etrisant la relation $\prec$. Sur $\cG_1$ cette relation co\"\i ncide avec la relation d'\'equivalence $\vdashv$ d\'efinie sur $R(f)=\Om(f)$~: les classes de r\'ecurrences par cha\^\i nes d'un diff\'eomorphisme  $f\in\cG_1$ sont donc exactement les ensembles faiblement transitifs maximaux, ce qui prouve le corollaire~\ref{c.coincide}.

\vskip 2mm

Nous utiliserons plusieurs fois par la suite le lemme suivant~:
\begin{lemm}\label{l.dashv2} Soit $\cG_1$ la partie r\'esiduelle de $\diff^1(M)$ sur laquelle $\prec_f$ et $\dashv_f$ co\"\i ncident. Pour tout $f\in \cG_1$ et tout compact invariant $K$, Lyapunov stable pour $f$, on a~:
\begin{equation}\label{e.dashv}
x\in K \mbox{ et } x\dashv_f y \Longrightarrow y\in K.
\end{equation}
Si de plus $K$ est r\'ecurrent par cha\^\i nes, $K$ est une classe de r\'ecurrence par cha\^\i nes.
\end{lemm}
\begin{demo}
Comme $K$ est stable au sens de Lyapunov, il poss\`ede des voisinages $U$ arbitrairement petits qui sont positivement invariant par $f$ (i.e. $f(U)\subset U$). On en d\'eduit que, si $x\in K$ et $y\in M$ v\'erifient $x\prec y$ alors $y\in K$. Comme les relations $\prec$ et $\dashv$ co\"\i ncident sur $M$, on a montr\'e l'implication~(\ref{e.dashv}).
L'ensemble $K$ contient donc toute classe de r\'ecurrence par cha\^\i nes qui le rencontre. Si $K$ est lui-m\^eme r\'ecurrent par cha\^\i nes, c'est exactement une classe de r\'ecurrence par cha\^\i nes.
\end{demo}

Afin de montrer le corollaire~\ref{c.hurley} (conjecture de Hurley pour la topologie $C^1$) nous devons pr\'ealablement montrer la proposition~\ref{p.quasiattracteur}.  Le premier item de cette proposition est presque une tautologie~: un quasi-attracteur $K$ poss\`ede des voisinages ouverts $U$ arbitrairement petits qui sont strictement invariants (i.e. $f(\overline{U})\subset U$). Il est donc stable au sens de Lyapunov et de plus les $\varepsilon$-pseudo-orbites, pour $\varepsilon$ assez petit, ne peuvent pas sortir de $U$. Un quasi-attracteur contient donc toute classe de r\'ecurrence par cha\^\i nes qu'il rencontre. S'il est lui m\^eme r\'ecurrent par cha\^\i nes, il est de ce fait une classe de r\'ecurrence par cha\^\i nes, ce qui montre le premier item.   

Le second item de la proposition~\ref{p.quasiattracteur} est la r\'eciproque g\'en\'erique du premier item~: soit $f$ un diff\'eomorphisme appartenant \`a la partie g\'en\'erique $\cG_1$ donn\'ee par le th\'eor\`eme~\ref{t.equivalence}, et soit $K$ une classe de r\'ecurrence par cha\^\i nes de $f$ qui est stable au sens de Lyapunov. On d\'efinit pour $\varepsilon>0$ l'ensemble $V(K,\varepsilon)=\{y\in M \; |\; \exists x\in K , x\dashv_\varepsilon y\}$. On remarque que pour tout $\varepsilon_1<\varepsilon_2$, l'adh\'erence de $V(K,\varepsilon_1)$ est incluse dans $V(K,\varepsilon_2)$. On en d\'eduit que $\bigcap_\varepsilon \overline{V(K,\varepsilon)}$ est une intersection d\'ecroissante de compacts qui est \'egale \`a $\bigcap_\varepsilon V(K,\varepsilon)$ et donc \'egale \`a $\{y\in M \; | \; \exists x\in K , x\dashv y\}$. Finalement, puisque $K$ est suppos\'e \^etre stable au sens de Lyapunov, $\bigcap_\varepsilon \overline{V(K,\varepsilon)}$  co\"\i ncide avec $K$ d'apr\`es le lemme~\ref{l.dashv2}.

On en d\'eduit que, pour tout voisinage $U$ de $K$, il existe $\varepsilon>0$ tel que l'ensemble 
$V(K,\varepsilon)=\{y\in M \; | \; \exists x\in K , x\dashv_\varepsilon y\}$ est inclus dans $U$. Remarquons que $V(K,\varepsilon)$ est un voisinage de $K$~; de plus, si $y$ est un point de $V(K,\varepsilon)$ et si $x_0\in K, x_1,\dots, x_k=y$ est une $\varepsilon$-pseudo-orbite joignant un point $x_0\in K$ \`a $y$, alors pour tout point $z$ tel que $d(f(y),z)<\varepsilon$, la suite $x_0\in K, x_1,\dots, x_k=y, x_{k+1}=z$ est une $\varepsilon$-pseudo-orbite joignant $x_0$ \`a $z$. Ceci montre $f(\overline{V(K,\varepsilon)})\subset int(V(K,\varepsilon))$. La classe de r\'ecurrence par cha\^\i nes $K$ poss\`ede donc une base de voisinages strictement invariants par $f$, et est donc un quasi-attracteur r\'ecurrent par cha\^\i nes, ce qui conclut la preuve de la proposition~\ref{p.quasiattracteur}. 

D'apr\`es \cite{MP}, il existe une partie r\'esiduelle $\cG_4$ de $\diff^1(M)$ telle que, pour tout $f\in\cG_4$, il existe une partie r\'esiduelle $R$ de $M$ telle que, pour tout point $x\in R$, l'ensemble $\omega(x,f)$ est un ensemble r\'ecurrent par cha\^\i nes, stable au sens de Lyapunov. Consid\'erons $f\in\cG_4\cap\cG_1$. Pour tout $x\in R$, l'ensemble $\omega(x,f)$ \'etant stable au sens de Lyapunov et r\'ecurrent par cha\^\i nes, d'apr\`es le lemme~\ref{l.dashv2}, il est  donc une classe de r\'ecurrence par cha\^\i nes stable au sens de Lyapunov, et donc un quasi-attracteur r\'ecurrent par cha\^\i nes, d'apr\`es la proposition~\ref{p.quasiattracteur}. L'union des bassins d'attraction des quasi-attracteurs r\'ecurrents par cha\^\i nes de $f$ contiennent donc la partie r\'esiduelle $R$ de $M$, ce qui termine la preuve de la conjecture de Hurley  en topologie $C^1$.

\subsection{Classes de r\'ecurrence par cha\^\i nes et orbites p\'eriodiques}

Montrons le corollaire~\ref{c.interieur}~: 
pour tout $f$, d'apr\`es le th\'eor\`eme fondamental des syst\`emes dynamiques de Conley, toute composante connexe de $\cR(f)$ est incluse dans une classe de r\'ecurrence par cha\^\i nes. Si celle-ci est d'int\'erieur non vide, elle contient, pour $f$ g\'en\'erique, une orbite p\'eriodique hyperbolique. Pour $f$ g\'en\'erique, la classe homocline de cette orbite p\'eriodique est une classe de r\'ecurrence par cha\^\i nes (voir la remarque~\ref{r.homoclasse}). Nous avons montr\'e que pour $f$ g\'en\'erique, toute composante connexe de $\cR(f)$ d'int\'erieur non-vide est incluse dans une classe homocline, ce qui conclut la preuve du corollaire~\ref{c.interieur}.

\vskip 2mm
Finalement, nous montrons le corollaire~\ref{c.isolee}.
Nous aurons besoin du lemme suivant, dont la preuve est calqu\'ee sur celle de \cite[Theorem A]{Ab}~:

\begin{lemm}\label{l.abdenur} Il existe un ensemble $\cG_5$ de $\diff^1(M)$ tel que, pour tout $f\in \cG_5$ et toute classe homocline isol\'ee $H(p,f)$ d'une orbite p\'eriodique hyperbolique $p$ de $f$, il existe un voisinage $\cU$ de $f$ dans $\diff^1(M)$ et $U$ de $H(p,f)$ dans $M$ tels que pour tout $g\in \cG_5\cap \cU$, $U$ ne contient qu'une seule classe homocline de $g$.
\end{lemm}
\begin{demo} Fixons une base d\'enombrable d'ouverts $\cO=\{O_n\}$ de $M$. Notons $\cG_6$ la partie r\'esiduelle de $\diff^1(M)$, donn\'ee par \cite{CMP} sur laquelle deux classes homoclines sont disjointes ou confondues et toute classe homocline varie contin\^ument avec le diff\'eomorphisme (celui-ci variant dans $\cG_6$). 

\`A tout diff\'eomorphisme $f\in \cG_6$ et tout ouvert $O_n$, nous associons le nombre $N(n,f)\in\NN\cup\{+\infty\}$ de classes homoclines de $f$ rencontrant $O_n$. Par le choix de $\cG_6$ les applications $f\mapsto N(n,f)$ sont semi-continues en restriction \`a $\cG_6$. La base d'ouverts $O_n$ \'etant d\'enombrable, il existe une partie r\'esiduelle $\cG_5\subset \cG_6$ en restriction \`a laquelle les fonctions $f\mapsto N(n,f)$ sont continues donc localement constantes. 

Soit $H(p,f)$ une classe homocline isol\'ee d'un diff\'eomorphisme $f\in\cG_5$. Soit $U_0$ un voisinage de $H(p,f)$ qui ne rencontre aucune autre classe homocline. Soient $O_{i_1},\dots,O_{i_k}$ un recouvrement fini de $H(p,f)$ par des ouverts de la base $\cO$, contenus dans $U_0$,  et soit $U$ un voisinage de $H(p,f)$ contenu dans l'union $\bigcup_1^k O_{i_j}$. Il existe un voisinage  ouvert $\cU$ de $f$ tel que, pour tout $g\in\cU\cap\cG_5$ et tout $j\in\{1,\dots, k\}$, $N(i_j,g)=N(i_j,f)$, c'est-\`a-dire que $O_{i_j}$ ne rencontre pas d'autre classe homocline que $H(p,g)$. Ainsi $U$ ne contient qu'une seule classe homocline de $g$.
  
\end{demo}

\begin{demo}[D\'emonstration du corollaire~\ref{c.isolee}]
Soient $f\in\cG_1\cap\cG_5$ et $H(p,f)$ une classe homocline isol\'ee de $f$.
Supposons par l'absurde que $H(p,f)$ n'est pas robustement r\'ecurrente par cha\^\i nes. Soit $V$ un voisinage filtrant de $H(p,f)$ (construit \`a l'aide d'une fonction de Lyapunov) et soit $\cV$ un voisinage de $f$, choisit assez petit pour que $V$ soit toujours un voisinage filtrant (c'est-\`a-dire qu'il est intersection d'un ouvert strictement invariant pour tout $g\in \cV$ et d'un ouvert strictement invariant pour tout $g^{-1}$, $g\in\cV$).

Pour tout $g\in\cV$, l'ouvert $V$ \'etant filtrant, l'ensemble $\cR(g)\cap V$ est inclus dans l'ensemble maximal invariant $\Lambda_g$ de $g$ dans $\bar V$ et dans $V$. De plus, les classes de r\'ecurrence par cha\^\i nes de $g$ incluses dans $\Lambda_g$ co\"\i ncident avec les classes de r\'ecurence par cha\^\i nes de la restriction de $g$ \`a $\La_g$~: en effet, $\La_g$ poss\`ede  une base de voisinages, qui sont des voisinages filtrants. Les pseudo-orbites joignant deux points de $\La_g$ ne peuvent donc pas s'\'eloigner de $\La_g$ et peuvent donc \^etre remplac\'ees par des pseudo-orbites de points de $\La_g$.

Par hypoth\`ese, il existe $g\in\cV$ tel que $\Lambda_g$ n'est pas r\'ecurrent par cha\^\i nes. La th\'eorie de Conley permet de montrer que la restriction de $g$ \`a $\La_g$ poss\`ede au moins deux classes distinctes de r\'ecurrence par cha\^\i nes, qui sont, comme nous l'avons vu, deux classes distinctes de r\'ecurrence par cha\^\i nes de la dynamique de $g$ sur $M$. Il existe donc une classe de r\'ecurrence par cha\^\i nes $K$ de $g$, incluse dans $\La_g$ et disjointe de $H(p,g)$. La th\'eorie de Conley montre l'existence d'un ouvert de $M$ strictement invariant $V'$ tel que $V'$ et le compl\'ementaire de $\bar V'$ contiennent chacun un et un seul des compacts $K$ et $H(p,g)$~: nous supposerons, pour fixer les id\'ees que $K\subset V'$ et $H(p,g)\subset M\setminus \bar V'$. L'ouvert $V'$ est encore strictement invariant pour tout diff\'eomorphisme $g'$ proche de $g$. Le closing lemma appliqu\'e \`a une orbite r\'ecurrente de $K$ permet de construire un diff\'eomorphisme $g'$ arbitrairement proche de $g$, poss\'edant une orbite p\'eriodique hyperbolique $p'$
contenue dans $V'\cap V$. Pour tout diff\'eomorphisme $\tilde g$ proche de $g$, les deux classes homoclines $H(p,\tilde g)$
et $H(p',\tilde g)$ sont distinctes et contenues dans $V\cap (M\setminus \bar V')$ et $V\cap V'$ respectivement.
Nous avons montr\'e que pour tout diff\'eomorphisme $\tilde g$ contenu dans un ouvert de diff\'eomorphismes arbitrairement proche de $f$, le voisinage $V$ de $H(p,f)$, choisi arbitrairement petit, contient deux classes homoclines distinctes. Ceci contredit le choix initial du diff\'eomorphisme $f\in \cG_5$.
\end{demo}
\section{Le cas conservatif}\label{s.conservatif}
\subsection{Les th\'eor\`emes ~\ref{t.connect},~\ref{t.equivalence} et \ref{t.conservatif} lorsque $dim(M)\geq 3$}\label{s.conservatif3}
 Le lemme de connection (connecting lemma) de Hayashi, tel qu'il est \'enonc\'e au th\'eor\`eme~\ref{t.connecting}, reste valide dans le monde conservatif, sa preuve restant inchang\'ee  (voir les remarques~\ref{r.connecting} et \ref{r.perturb}).

Soit $f$ un diff\'eomorphisme d'une vari\'et\'e compacte $M$ de dimension $d$ pr\'eservant une forme volume $\omega$.
On d\'efinit une {\em bo\^\i te de perturbation d'ordre $N$ pour $(f,\cU,\omega)$}, o\`u $\cU$ est un voisinage ouvert de $f$ dans $\diff^1(M)$, en demandant que le diff\'eomorphisme $g$, dans la d\'efinition~\ref{d.boites}, pr\'eserve $\omega$ (c'est-\`a-dire que $g$ appartient \`a $\cU\cap \diff^1_\omega(M)$). 


\begin{theo}\label{t.connectingconservatif} Soit $f$ un diff\'eomorphisme d'une vari\'et\'e compacte $M$ de dimension $d$ pr\'e\-ser\-vant une forme volume $\omega$. Pour tout voisinage $\cU$ de $f$ dans $\diff^1(M)$ il existe $N>0$ v\'erifiant~: pour tout point  $x\in M$, il existe  une carte locale $\varphi:U_x\to\RR^d$ en $x$  telle que tout  cube quadrill\'e de $(U_x,\varphi)$, disjoint de ses $N$ premiers it\'er\'es est une bo\^\i te de perturbation d'ordre $N$ pour $(f, \cU,\omega)$.

De plus, cette propri\'et\'e est encore v\'erifi\'ee par les cartes proches de $\varphi$ en topologie $C^1$.
\end{theo}
La seule pr\'ecaution additionnelle (par rapport au cas non-conservatif) consiste \`a v\'erifier que le lemme~\ref{l.perturb} reste valide pour les diff\'eomorphismes pr\'eservant $\omega$~; nous renvoyons \`a \cite[proposition 5.1.1]{Ar2} o\`u ce lemme est montr\'e dans le cadre conservatif.

On transpose alors sans modification la preuve du th\'eor\`eme~\ref{t.connect} (en utilisant le th\'eor\`eme~\ref{t.connectingconservatif} au lieu du th\'eor\`eme~\ref{t.connecting}) et on obtient:
\begin{theo}\label{t.connectconservatif}
Soit $f$ un diff\'eomorphisme d'une vari\'et\'e compacte  $M$  pr\'eservant une forme volume $\omega$. On suppose de plus que toutes les orbites p\'eriodiques de $f$ sont hyperboliques.
Soit $\cU$ un $C^1$-voisinage de $f$ dans $\diff^1_\omega(M)$. Alors, pour toute paire $(x,y)$ de points de $M$ telle que $x\dashv y$, il existe un diff\'eomorphisme $g$ dans $\cU$ et un entier $n>0$ tel que $g^n(x)=y$.
\end{theo}

Soit $M$ une vari\'et\'e compacte connexe de dimension $d\geq 3$, munie d'une forme volume $\omega$. Il existe alors (voir \cite{Ro1}) une partie r\'esiduelle $\cG_\omega$ de $\diff^1_\omega(M)$ telle que tout point p\'eriodique d'un diff\'eomorphisme $f\in\cG_\omega$ est hyperbolique. Ceci permet de r\'ep\'eter sans modificaton la preuve du th\'eor\`eme~\ref{t.equivalence} donn\'ee \`a la section~\ref{s.equivalence}, et on obtient~:
\begin{theo}\label{t.equivalenceconservatif} Il existe une partie r\'esiduelle $\cG$ de $\diff^1_\omega(M)$  telle que pour tout diff\'eomorphisme $f$ de $\cG$ et tout couple $(x,y)$ de points de $M$ on a~:
$$x\dashv_f y\Longleftrightarrow x\prec_fy .$$
\end{theo}
 
\begin{demo}[Fin de la preuve du th\'eor\`eme~\ref{t.conservatif} dans le cas o\`u $\dim M \geq 3$]
Finalement, pour tout diff\'eomorphisme $f\in \diff^1_\omega(M)$, l'ensemble non-errant $\Om(f)$ co\"\i ncide avec la vari\'et\'e~; on a vu que ceci implique que $M$ est alors une classe de r\'ecurrence par cha\^\i nes pour $f$. Si $f$ est g\'en\'erique, on a$\prec_f=\dashv_f$ d'apr\`es le th\'eor\`eme~\ref{t.equivalenceconservatif}~; ceci implique que $f$ est un diff\'eomorphisme transitif~; de plus, d'apr\`es~\cite{PuRo,Ro1}, ses points p\'eriodiques hyperboliques sont denses. D'apr\`es \cite[corollaire 19]{Ar} dans le cadre conservatif (quand tous les points p\'eriodiques sont hyperboliques, voir \cite[paragraphe 1.5]{Ar}), il existe une partie r\'esiduelle de $\diff^1_\omega(M)$ sur laquelle tout compact invariant transitif de $M$ contenant un point p\'eriodique hyperbolique $p$ est contenu dans la classe homocline de $p$.  Ainsi, pour $f\in\diff^1_\omega(M)$ g\'en\'erique (avec $\dim M\geq 3$), la vari\'et\'e $M$ toute enti\`ere est incluse dans la classe homocline de chacune de ses orbites p\'eriodiques (qui sont toutes hyperboliques). 

Nous avons donc montr\'e le th\'eor\`eme~\ref{t.conservatif}, dans le cas o\`u $\dim M \geq 3$.
\end{demo}

Pour les sections suivantes, nous utiliserons le r\'esultat classique suivant (adapt\'e de \cite{F})~:
\begin{prop}[Lemme de Franks conservatif, voir par exemple~\cite{BDP}, Lemma 7.6]\label{p.franks} Soit $f$ un diff\'eomorphisme de classe $C^r$ d'une vari\'et\'e compacte $M$ pr\'eservant une forme volume $\omega$. Pour tout $C^1$-voisinage $\cU$ de $f$ dans $\diff^r_\omega(M)$, il existe un $C^1$-voisinage $\cO$ de $f$ dans $\diff^r_\omega(M)$ et $\varepsilon>0$ avec la propri\'et\'e suivante~:

\'Etant donn\'es un diff\'eomorphisme $\tilde f\in\cO$,  une partie finie $E$ de $M$, un voisinage $V$ de $E$ et en tout point $x\in E$ une $\varepsilon$-perturbation $B_x\colon T_xM\to T_{\tilde f(x)}M$ de $D_x \tilde f$, il existe $g\in \cU$ co\"\i ncidant avec $\tilde f$ hors de $V$ et sur $E$, tel que pour tout point $x\in E$, $D_xg=B_x$.
\end{prop}
Bien que l'\'enonc\'e soit l\'eg\`erement plus fort que celui de~\cite{BDP} ($g$ est de classe $C^r$, la partie $E$ n'est pas suppos\'e invariante et le r\'esultat est uniforme pour $\tilde f\in \cO$), la preuve est rigoureusement identique.

\subsection{Diff\'eomorphismes conservatifs des surfaces}\label{s.conservatif2}
\subsubsection{Au voisinage des points elliptiques}

Soit $S$ une surface compacte munie d'une forme volume (forme d'aire) $\omega$.  Il existe une partie r\'esiduelle $\cG_\omega$ de $\diff^1_\omega$ telle que, pour tout $f\in \cG_\omega$ toute orbite p\'eriodique est ou bien hyperbolique, ou bien la diff\'erentielle (\`a la p\'eriode) de cette orbite est conjug\'ee \`a une rotation irrationnelle (voir \cite{Ro1}). Par cons\'equent l'ensemble $Per_{N_0}(f)$ des orbites p\'eriodiques de p\'eriode inf\'erieure \`a $N_0$ est fini.  

\begin{figure}
\begin{center}
\input{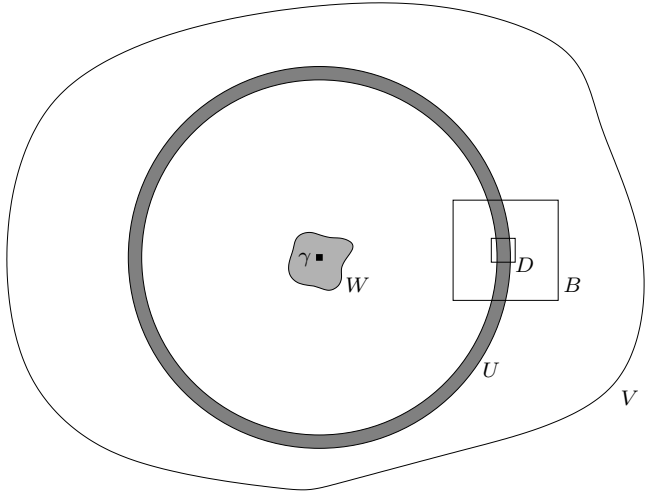}
\end{center}
\caption{Bo\^\i te de perturbation au voisinage d'un point fixe elliptique. \label{f.elliptique}}
\end{figure}

Au voisinage des points p\'eriodiques elliptiques, nous avons besoin d'un r\'esultat \'equivalent \`a la proposition~\ref{p.local}.
\begin{prop}\label{p.localelliptique} Soit $f$ un diff\'eomorphisme d'une surface compacte pr\'eservant une forme volume $\omega$, et $\cU$ voisinage de $f$ dans $\diff^1_\omega(M)$ pour la topologie $C^1$. Soit $N$  l'entier associ\'e \`a l'ouvert $\cU$ par le th\'eor\`eme~\ref{t.connectingconservatif}.
Pour toute orbite p\'eriodique $\gamma$ de $f$ dont la diff\'erentielle (\`a la p\'eriode) est conjugu\'ee \`a une rotation irrationnelle, et tout voisinage $V$ de $\gamma$, il existe un ouvert $U$, un voisinage $W$ de $\gamma$, une bo\^\i te de perturbation $B$, un compact D contenu dans l'int\'erieur du carreau central de $B$ et un entier $n_0$ tels que (voir la figure~\ref{f.elliptique})~:
\begin{enumerate}
\item L'ouvert $V$ contient les compacts $\bar W$, $\bar U$ et le support de la bo\^\i te de perturbation $B$. De plus $\bar U$ et $\bar W$ sont disjoints.
\item Pour tout point $z$ de $U$, il existe un it\'er\'e positif $f^n(z)$ avec $n\in\{0,\dots, n_0\}$
qui appartient \`a l'int\'erieur de $D$. De m\^eme, il existe un it\'er\'e n\'egatif $f^{-n}(z)$ avec $n\in\{0,\dots, n_0\}$ qui appartient \`a l'int\'erieur de $D$. (On pense \`a $U$ comme \'etant une r\'eunion de couronnes autour de chaque point de $\gamma$.)
\item Pour toute orbite finie $z,\cdots,f^n(z)$ rencontrant $M\setminus V$ et $W$, il existe un it\'er\'e
$f^{m}(z)$ avec $m\in\{0,\dots,n\}$ qui appartient \`a $U$.
\end{enumerate}
\end{prop}
\begin{demo}
Nous noterons $q$ la p\'eriode de l'orbite $\gamma$.
Consid\'erons une carte de perturbation $\varphi$ donn\'ee par le th\'eor\`eme~\ref{t.connectingconservatif} appliqu\'e en un point $z_0$ de l'orbite $\gamma$. De ce fait, nous pouvons travailler dans $\mathbb{R}^2$ au voisinage du point $0=\varphi(z_0)$ avec l'application $F=\varphi\circ f^q\circ \varphi^{-1}$.

Nous consid\'erons n'importe quel carr\'e quadrill\'e $B_1$, contenu dans un petit voisinage de $0$. Pour tout $t\in ]0,1[$, l'image $B_t$ de $B_1$ par une homoth\'etie de rapport $t$ centr\'ee en $0$ est encore un carr\'e quadrill\'e.
Quitte \`a choisir $B_1$ assez petit, pour tout $t\in ]0,1[$ les bo\^\i tes $\varphi^{-1}(B_t)$ sont disjointes de leur $N$ premiers it\'er\'es et d\'efinissent des  bo\^\i tes de perturbation d'ordre $N$ pour $(f,\cU,\omega)$. Consid\'erons un compact $D_1$ d'int\'erieur non vide contenu dans l'int\'erieur du carreau central de la bo\^\i te $B_1$ et $D_t$ son image par l'homoth\'etie de rapport $t$.

Par hypoth\`ese, la diff\'erentielle $A\colon \RR^2\mapsto \RR^2$ de l'application $F$ en $0$ est conjugu\'ee \`a une rotation irrationnelle $R$ de $\RR^2$ par un automorphisme $H\in GL(2,\RR)$ (ainsi $A=H^{-1}RH$).
On note $C_t$ l'image par $H^{-1}$ du cercle de $\RR^2$ de rayon $t$. Quitte \`a multiplier $H$ par une homoth\'etie, on peut supposer \'egalement que $C_1$ contient un point de l'int\'erieur de $D_1$ (et $C_t$ un point de l'int\'erieur de $D_t$ pour tout $t\in]0,1[$). On fixe une constante $\eta\in]0,1[$ telle que, pour tout $t\in]0,1[$ et tout $t'\in ]t,(1+\eta)t[$, la courbe $C_{t'}$ rencontre l'int\'erieur de $D_t$.

L'application $A$ pr\'eserve chaque courbe $C_t$.
D'autre part, il existe un entier $n_0\geq 0$ tel que pour tout $t\in ]0,1[$ et tout $t'\in ]t,(1+\eta)t[$, tout point $z\in C_{t'}$ poss\`ede un it\'er\'e positif $A^n(z)$ par $A$, avec $n\leq n_0$, qui appartient \`a l'int\'erieur de $D_t$. Pour $\varepsilon>0$ suffisamment proche de $0$, l'application $F$ est proche de $A$ en topologie $C^1$ sur la boule de centre $0$ et de rayon $\varepsilon$. On en d\'eduit que pour $t$ suffisamment petit les deux propri\'et\'es suivantes sont satisfaites~:
\begin{enumerate}
\item Pour tout $t'\in ]t,(1+\eta)t[$, tout point $z$ de $C_{t'}$ poss\`ede un it\'er\'e positif $F^n(z)$ par $F$ avec $n\in\{0,\dots, n_0\}$ appartenant \`a l'int\'erieur de $D_t$.
\item Pour tout $n\in\{0,\dots,n_0\}$ et tout point $z\in C_t$, le point $F^n(z)$ appartient \`a une courbe $C_{t'}$ avec $t'<(1+\eta)t$.
\end{enumerate}
Les m\^emes propri\'et\'es sont encore v\'erifi\'ees par $F^{-1}$.
Fixons un choix $t_0$ de $t$ de sorte que ces deux propri\'et\'es soient satisfaites simultan\'ement pour $F$ et $F^{-1}$.

Soit $W_{0}$ un petit voisinage connexe de de $0$ qui ne rencontre pas $C_{t_0}$, et $z$ un point de $W_0$. Alors, on est dans un des cas suivants~:
\begin{itemize}
\item
Ou bien l'orbite positive de $z$ par $F$ ne rencontre aucune courbe $C_t$ avec $t\in ]t_0,(1+\eta)t_0[$. Dans ce cas, la seconde propri\'et\'e montre que l'orbite de $z$ par $F$ reste contenue dans l'union des courbes $C_t$ avec $t\leq t_0$.
\item
Ou bien l'orbite positive de $z$ par $F$ rencontre une courbe $C_t$ avec $t\in]t_0,(1+\eta)t_0[$ et, d'apr\`es la premi\`ere propri\'et\'e, intersecte l'int\'erieur de $D_{t_0}$.
\end{itemize}
On obtient encore le m\^eme r\'esultat pour l'orbite n\'egative de $z$ par $F$.
Soit $U_0$ l'union des courbes $C_t$ avec $t\in ]t_0,(1+\eta)t_0[$.

On peut supposer dans cette construction que tous les cercles $\varphi^{-1} (C_t)$, avec $t\in ]0,1[$, et leurs $q$
premiers it\'er\'es par $f$ sont contenus dans $V$ ($q$ est la p\'eriode de $\gamma$). Il suffit pour conclure la d\'emonstration de choisir pour $B$ et $D$ l'image de $B_{t_0}$ et $D_{t_0}$ par $\varphi^{-1}$. L'ouvert $U$ est l'union de $\varphi^{-1}(U_{0})$ et de ses $q$ premiers it\'er\'es par $f$. On d\'efinit de la m\^eme fa\c con $W$ \`a partir de $W_0$.
\end{demo}

\begin{rema}\label{r.boiteelliptique} Dans l'\'enonc\'e de la proposition~\ref{p.localelliptique}, on peut demander de plus que l'adh\'erence de la bo\^\i te de perturbation $B$ soit disjointe de ses $m$ premiers it\'er\'es o\`u $m$ est un entier arbitraire.
\end{rema}

\subsubsection{Preuve du th\'eoreme~\ref{t.connect} pour les diff\'eomorphismes conservatifs des surfaces}

Afin d'appliquer le corollaire~\ref{c.coloriage}, il faut se ramener \`a un compact invariant ne contenant pas
d'orbite p\'eriodique elliptique de basse p\'eriode. C'est pourquoi cette fois, nous devons commencer par traiter le cas
des orbites elliptiques en construisant, au voisinage, des bo\^\i tes de perturbation gr\^ace \`a la proposition~\ref{p.localelliptique}. Si $K$ est l'ensemble maximal invariant hors d'un voisinage $W$ des orbites elliptiques de basse p\'eriode, le corollaire~\ref{c.coloriage} lui associe un nombre fini de bo\^\i tes de perturbation de supports disjoints couvrant l'espace des orbites. En proc\'edant de cette fa\c con, il est n\'ecessaire de s'assurer que les nouvelles bo\^\i tes ainsi construites ont leur support disjoint des bo\^\i tes que l'on avait pr\'elablement introduites pr\`es des orbites elliptiques~; ceci nous am\`enera \`a reprendre l\'eg\`erement la preuve du corollaire~\ref{c.coloriage}. Comme dans la preuve de la section~\ref{s.connect}, on peut alors gr\^ace \`a la proposition~\ref{p.local} ajouter de nouvelles bo\^\i tes de perturbation au voisinage des orbites p\'eriodiques hyperboliques de basse p\'eriode. Cette nouvelle collection de bo\^\i tes de perturbation permet de terminer la preuve comme dans le cas non-conservatif.

Nous allons d\'etailler maintenant les \'etapes de la d\'emonstration qui sont diff\'erentes du cas g\'en\'eral.

\noindent {\bf Adaptation de la preuve du corollaire~\ref{c.coloriage} aux diff\'eomorphismes conservatifs des surfaces~:}

Consid\'erons un diff\'eomorphisme $f\in \cG_\omega$ de $S$ et $\cU$ un voisinage de $f$ dans $\diff^1(M)$. Nous supposerons, quitte \`a le r\'eduire, que $\cU$ satisfait la propri\'et\'e de la remarque~\ref{r.u}. Soient enfin $x$ et $y$ deux points de $S$ tels que $x\dashv y$. Comme toujours, nous poserons $N_0=10dN$ o\`u $N$ est l'entier associ\'e \`a $\cU$ par le th\'eor\`eme~\ref{t.connecting}. On choisit une constante $\delta>0$ telle que les vari\'et\'es invariantes locales $W^s_\delta(z)$ et $W^u_\delta(z')$ de points $z$ et $z'$ de $Per_{3N_0}(f)$ soient deux \`a deux disjointes.

Nous supposerons (comme cela a \'et\'e expliqu\'e par la remarque~\ref{r.casparticulier}) que $x$ et $y$ n'appartiennent pas \`a $Per_{3\kappa_dN_0}(f)$ et que $y$ n'est pas un point de l'orbite positive de $x$ (la raison pour laquelle nous avons remplac\'e $Per_{N_0}(f)$ par $Per_{3\kappa_dN_0}(f)$ par rapport \`a la section~\ref{s.connect} vient de l'ordre l\'eg\`erement diff\'erent dans lequel nous faisons la preuve).

Pour toute orbite p\'eriodique elliptique $\gamma$ de p\'eriode inf\'erieure \`a $3\kappa_d N_0$, on choisit un voisinage $V(\gamma)$ de fa\c con \`a ce que les $V(\gamma)$ soient deux \`a deux disjoints et ne contiennent pas les points $x$ et $y$. Notons $B(\gamma)$ la bo\^\i te de perturbation ainsi que $D(\gamma)$ le compact contenu dans le carreau central $C(\gamma)$ de $B(\gamma)$ obtenus \`a l'aide de la proposition~\ref{p.local}. Gr\^ace \`a la remarque~\ref{r.boiteelliptique} et en choisissant les voisinages $V(\gamma)$ suffisamment petits, on peut demander que les adh\'erences des   bo\^\i tes $B(\gamma)$ et leurs $4N_0$ premiers it\'er\'es soient deux \`a deux disjointes. 

\begin{lemm}\label{l.separation} Sous ces hypoth\`eses, il existe $r>0$ tel que, pour tout point $z\in M$, il existe $n\in\{0,\dots,2N_0\}$ tel que les it\'er\'es $f^i(B(z,r))$, avec $i\in\{n,\dots,n+N_0\}$, de la boule $B(z,r)$ sont disjoints des adh\'erences des supports des bo\^\i tes $B(\gamma)$.
\end{lemm}
\begin{demo}
Soit $z$ un point de $M$. 
\begin{itemize}
\item Ou bien, pour tout $i\in\{0,\dots,N_0\}$, le point $f^i(z)$ ne rencontre pas l'union des adh\'erences des supports des bo\^\i tes $B(\gamma)$.
\item Ou bien il existe $i\in\{0,\dots, N_0\}$ tel que $f^i(z)$ appartient \`a l'adh\'erence du support d'une bo\^\i te $B(\gamma)$. Alors, les points $f^{i+N+1+k}(z)$ pour $k\in \{0,\dots, N_0\}$ sont hors des adh\'erences des supports des bo\^\i tes $B(\gamma)$, par construction de ces bo\^\i tes.
\end{itemize}

Nous avons montr\'e qu'il existe pour tout point $z$ un entier $n\in\{0,\dots,2N_0\}$ tel que les it\'er\'es $f^i(z)$, avec $i\in\{n,\dots,n+N_0\}$, de $z$ sont disjoints des adh\'erences des supports des bo\^\i tes $B(\gamma)$. Il existe donc un rayon $r(z)$ tel que les it\'er\'es $f^i(B(z,r(z)))$, avec $i\in\{n,\dots,n+N_0\}$, de la boule $B(z,r(z))$ sont disjoints des adh\'erences des supports des bo\^\i tes $B(\gamma)$. Un argument de compacit\'e de $M$ permet de choisir un rayon $r$ uniforme.

\end{demo}

Finalement, on note aussi $W(\gamma)$, $U(\gamma)$ et $n_0(\gamma)$ les ouverts et entiers annonc\'es par la proposition~\ref{p.localelliptique}.

Soit $K$ l'ensemble des points dont l'orbite n'intersecte pas l'union de ces ouverts $W(\gamma)$, avec $\gamma$ orbite elliptique de p\'eriode inf\'erieure \`a $3\kappa_d N_0$. Le compact $K$ est invariant et ne contient pas d'orbite p\'eriodique non-hyperbolique de p\'eriode inf\'erieure \`a $3\kappa_d N_0$. On applique le th\'eor\`eme~\ref{t.coloriage}
et le corollaire~\ref{c.coloriage2} \`a l'entier $3N_0$ afin d'obtenir un ouvert $U$ disjoint de ses $3N_0$ premiers it\'er\'es. Le th\'eor\`eme~\ref{t.coloriage} permet de choisir les composantes connexes de $U$ de diam\`etre arbitrairement petit. On les choisit de diam\`etre plus petit que le rayon $r$ du lemme~\ref{l.separation}. Soit $C$ l'une de ces composantes, et $n(C)\in\{0,\dots,2N_0\}$ un it\'er\'e tel que $f^{n(C)+i}(C)$ est disjoint des supports des boites $B(\gamma)$ pour tout $i\in\{0,\dots,N_0\}$. On consid\`ere $\tilde U$ l'ouvert obtenu en rempla\c cant chaque composante $C$ de $U$ par $f^{n(C)}(C)$. C'est un ouvert, disjoint de ses $N_0$ premiers it\'er\'es et disjoint du support de toutes les bo\^\i tes $B(\gamma)$, qui v\'erifie les conclusions du th\'eor\`eme~\ref{t.coloriage} (la constante $\delta$ a pu \^etre modifi\'ee).

On peut reprendre la preuve du corollaire~\ref{c.coloriage} avec le compact $K$ et l'ouvert $\tilde U$ afin d'obtenir une famille de bo\^\i tes de perturbation $\cB_0$ v\'erifiant les conclusions du corollaire~\ref{c.coloriage} et disjointes des  supports des bo\^\i tes $B(\gamma)$. On applique enfin la proposition~\ref{p.local} au voisinage des orbites p\'eriodiques hyperboliques de p\'eriode inf\'erieure ou \'egale \`a $3N_0$ et on note $\cB$ l'ensemble des bo\^\i tes de perturbations introduites.

\medskip
\noindent
{\bf Adaptation de la fin de la d\'emonstration du th\'eor\`eme~\ref{t.connect}~:}

Afin d'uniformiser les notations avec le cas des orbites p\'eriodiques hyperboliques,
pour toute orbite elliptique $\gamma$ (de p\'eriode inf\'erieure \`a $3\kappa_d N_0$), on note $\cE(\gamma)=\cS(\gamma)=\{B(\gamma)\}$, de m\^eme $\cC(\cE,\gamma)=\cC(\cS,\gamma)=\{C(\gamma)\}$ et
$\cD(\cE,\gamma)=\cD(\cS,\gamma)=\{D(\gamma)\}$.
La fin de la d\'emonstration est alors identique au cas sans orbite elliptique (sections~\ref{s.regroupe}, \ref{s.regroupeencore} et \ref{s.fin}).

\subsubsection{Preuve du th\'eor\`eme~\ref{t.conservatif} dans le cas des surfaces}

La transitivit\'e annonc\'ee par le th\'eor\`eme~\ref{t.conservatif} est une cons\'equence imm\'ediate du th\'eor\`eme~\ref{t.equivalence} comme en dimension $d\geq 3$. La nouvelle difficult\'e du th\'eor\`eme~\ref{t.conservatif} consiste \`a montrer que la surface est une unique classe homocline. Pour cela, nous avons besoin de cr\'eer des orbites p\'eriodiques hyperboliques.

Voici l'une des clefs de la d\'emonstration que nous \'enon\c cons en classe $C^r$, avec $r\geq 1$, pour pr\'eparer les sections suivantes\footnote{M.-C. Arnaud nous a signal\'e que \cite{N2} permet de remplacer avantageusement la proposition~\ref{p.elliptiquehyperbolique} ainsi que son corollaire principal, le corollaire~\ref{c.elliptiquehyperbolique}. En effet, S. Newhouse montre le corollaire~\ref{c.elliptiquehyperbolique} en topologie $C^r$, ce qui est suffisant pour les applications que nous donnons (voir la fin de la d\'emonstration du th\'eor\`eme~\ref{t.conservatif} pour les surfaces ci-dessous et la d\'emonstration du lemme~\ref{l.homoclinedense}. }.

\begin{prop}\label{p.elliptiquehyperbolique} Soient $r\geq 1$, $f\in \diff^r_\omega(S)$ et $\cU$ un $C^1$-voisinage de $f$ dans $\diff^r_\omega(S)$. Il existe $N\geq 1$ tel que pour toute orbite p\'eriodique $\gamma$ de $f$ de p\'eriode sup\'erieure \`a $N$ et tout voisinage $V$ de $\gamma$, il existe $g\in\cU$, co\"\i ncidant avec $f$ hors de $V$ et le long de $\gamma$, pour lequel $\gamma$ est une orbite p\'eriodique hyperbolique.
\end{prop}
\begin{demo}
D'apr\`es le lemme de Franks (proposition~\ref{p.franks}), il suffit de montrer que pour tout $\varepsilon>0$, il
existe $N\geq 1$ tel que pour toute orbite p\'eriodique de p\'eriode sup\'erieure \`a $N$, il existe une $\varepsilon$-perturbation de la diff\'erentielle le long de l'orbite de fa\c con que le produit (le long de l'orbite) des applications lin\'eaires ainsi obtenues soit hyperbolique. Il s'agit en fait d'un r\'esultat sur les suites finies, de longueur sup\'erieure \`a $N$, de matrices de $SL(2,\RR)$ de norme uniform\'ement born\'ee~: le lemme~\ref{l.elliptiquehyperbolique} ci-dessous permet d'obtenir par une perturbation de la diff\'erentielle, une valeur propre r\'eelle. Une perturbation aussi petite que l'on veut permet alors d'obtenir une valeur propre (r\'eelle) de module diff\'erent de $1$ ce qui, dans $SL(2,\RR)$, assure l'hyperbolicit\'e.
\end{demo}

\begin{lemm}\label{l.elliptiquehyperbolique} Pour tout $\varepsilon>0$, il existe $N\geq 1$ tel que, pour tout $n\geq N$ et toute suite finie $A_0,\dots,A_n$ d'\'el\'ements finis de $SL(2,\RR)$, il existe $\alpha_0,\dots,\alpha_n$ dans $]-\varepsilon,\varepsilon[$ tels que l'on a la propri\'et\'e suivante~:

Pour tout $i\in \{0,\dots,n\}$, on note $B_i=R_{\alpha_i}\circ A_i$ la compos\'ee de $A_i$ avec la rotation $R_{\alpha_i}$ d'angle $\alpha_i$. Alors,  la matrice
$B_n\circ B_{n-1}\circ \dots \circ B_0$ a ses valeurs propres r\'eelles.
\end{lemm}
\begin{demo}
Les matrices de $SL(2,\RR)$ agissent sur le cercle $\SS^1$ des demi-droites vectorielles de $\RR^2$. La rotation d'angle $\alpha$ s'\'ecrit $x\mapsto x+\alpha$ et le demi-tour $R_{1/2}$ s'\'ecrit $x\mapsto x+\frac{1}{2}$.

{\bf Affirmation.} {\it Pour tout $\varepsilon>0$, il existe $\eta>0$ tel que, pour toute matrice $A$ de $SL(2,\RR)$, l'une des deux propri\'et\'es suivantes est v\'erifi\'ee~:
\begin{enumerate}
\item il existe $\varphi\in \SS^1$ et $\theta\in ]0,\varepsilon[$ tels que
$$R_\theta \circ A \circ R_\theta (\varphi)=-R_{-\theta} \circ A \circ R_{-\theta} (\varphi)~;$$
\item pour tout $\varphi\in \SS^1$, la distance entre $A(\varphi+\varepsilon)$ et $A(\varphi)$ est plus grande que $\eta$.
\end{enumerate}}
\begin{demo} En effet, $A$ est la compos\'ee $O_1\circ D \circ O_2$ o\`u $O_1$ et $O_2$ sont des rotations et $D$ est une matrice diagonale. Puisque les rotations commutent entre elles et sont des isom\'etries de $\SS^1$, il suffit de montrer l'affirmation pour les matrices diagonales, de valeurs propres $\lambda, \lambda^{-1}$.

Si $\lambda$ est  sup\'erieure \`a une constante $\lambda_0$ (tr\`es grande par rapport \`a $\varepsilon^{-1}$), il existe $\theta\in]0,\varepsilon[$ tel que le vecteur $R_\theta\circ D\circ R_\theta ((0,1))$ est colin\'eaire \`a $(1,0)$. Par sym\'etrie  on obtient que le vecteur  $R_{-\theta}\circ D\circ R_{-\theta} ((0,1))$ est lui aussi colin\'eaire \`a $(1,0)$ et plus pr\'ecis\'ement~: $R_{-\theta}\circ D\circ R_{-\theta}((0,1))=-R_\theta\circ D\circ R_\theta ((0,1))$, ce qui
donne le premier cas de l'affirmation.

L'ensemble des matrices diagonales de valeurs propres born\'ees par $\lambda_0$ \'etant compact, il existe $\eta>0$ tel que pour toute matrice $A$ dans cet ensemble, le second cas de l'affirmation est v\'erifi\'e.
\end{demo}
Posons $N>\frac1{2\eta}$. Soit $A_0, \dots ,A_n$, $n\geq N$, une famille finie de matrice de $SL(2,\RR)$. 

Supposons d'abord qu'il existe $1\leq k\leq \ell\leq n$ tels que la matrice $A=A_\ell\circ A_{\ell-1}\circ \cdots \circ A_k$ v\'erifie le premier cas de l'affirmation pr\'ec\'edente. Il existe alors $\theta\in ]0,\varepsilon[$ et une demi-droite $\xi\in \SS^1$ telle que $M_1(\xi)=-M_{-1}(\xi)$ o\`u l'on a pos\'e, pour $t\in [-1,1]$,
$$M_t=(A_n\circ\cdots\circ A_{\ell+1})\circ R_{t\theta}\circ (A_{\ell}\circ \cdots \circ A_k)\circ R_{t\theta}\circ (A_{k-1}\circ \cdots\circ A_0).$$
Par continuit\'e de la famille, quand $t$ parcourt $[-1,1]$, la demi-droite $M_t(\xi)$ contient un demi-cercle de $\SS^1$ et donc contient $\xi$ ou $-\xi$. Il existe donc une valeur $t\in [-1,1]$ telle que $M_t(\xi)$ est colin\'eaire \`a $\xi$. La matrice $M_t$ a donc ses valeurs propres r\'eelles. On pose alors $\alpha_{k-1}=\alpha_{\ell}=t\theta$ et $\alpha_i=0$ pour $i\not\in\{k-1,\ell\}$.

Dans le cas contraire, on pose pour $t\in [0,1]$ et $i\in \{0,\dots,n\}$.
$$M_{i,t}=(A_n\circ\cdots\circ A_{i+1})\circ (R_{t\varepsilon}\circ A_i)\circ(R_{\varepsilon}\circ A_{i-1}\cdots \circ R_{\varepsilon}A_1\circ R_{\varepsilon}\circ A_0).$$
Puisque $A_n\circ\cdots\circ A_{i+1}$ v\'erifie le second cas de l'affirmation (par hypoth\`ese, il ne v\'erifie pas le premier cas), on peut montrer que pour tout $\varphi\in \SS^1$ la diff\'erence $M_{i,1}(\varphi)-M_{i,0}(\varphi)$ est minor\'ee par $\eta$. En remarquant que $M_{i,0}=M_{i-1,1}$, on forme une famille continue de matrice joignant $M_{0,0}$ \`a $M_{n,1}$.
La variation totale de l'image d'un point $\varphi\in \SS^1$ le long de cette famille est minor\'ee par $N\eta>1$. Il existe donc une application $M_{i,t}$ de cette famille qui poss\`ede un point fixe sur $\SS^1$. Cette matrice a donc une valeur propre r\'eelle. On pose donc $\alpha_i=t\varepsilon$, $\alpha_j=\varepsilon$ pour $j<i$ et $\alpha_j=0$ pour $j>i$.
\end{demo}

\begin{coro}\label{c.elliptiquehyperbolique} Soit $S$ une surface compacte munie d'une forme volume $\omega$. Il existe un ensemble r\'esiduel $\cG$ de $\diff^1_\omega(S)$ tel que, pour tout $f\in \cG$,
l'ensemble des orbites p\'eriodiques hyperboliques de $f$ est dense dans $S$.
\end{coro}
\begin{demo}
D'apr\`es~\cite{Ro1}, les orbites p\'eriodiques d'un diff\'eomorphisme g\'en\'erique $f\in\diff^1_\omega(S)$ sont toutes hyperboliques ou elliptiques. Par cons\'equent, pour tout $n\geq 1$, 
il existe un ouvert dense de $\diff^1_\omega(S)$ sur lequel les orbites de p\'eriode $n$ sont en nombre fini. D'autre part, le closing lemma de Pugh entra\^\i ne que les orbites p\'eriodiques d'un diff\'eomorphisme g\'en\'erique sont denses dans $S$. Donc pour un diff\'eomorphisme g\'en\'erique, les orbites p\'eriodiques de p\'eriode plus grande qu'une constante arbitraire sont denses dans $S$.

Fixons $\varepsilon>0$ et montrons que l'ensemble des diff\'eomorphismes dont l'ensemble des orbites p\'eriodiques hyperboliques est $\varepsilon$-dense dans $S$ est un ouvert dense de $\diff^1_\omega(S)$. Puisque les orbites p\'eriodiques hyperboliques peuvent \^etre suivies localement contin\^ument avec $f$, cet ensemble est clairement ouvert~; il reste \`a montrer la densit\'e.

Fixons donc un diff\'eomorphisme g\'en\'erique $f$ et un voisinage $\cU$ de $f$ qui v\'erifie la propri\'et\'e de la remarque~\ref{r.u}. Soit $N$ l'entier donn\'e par la proposition~\ref{p.elliptiquehyperbolique}. L'ensemble des orbites de $f$ de p\'eriode plus grande que $N$ est dense dans $S$. On consid\`ere donc un nombre fini d'entre elles dont l'union est $\varepsilon$-dense. La proposition~\ref{p.elliptiquehyperbolique} permet de les rendre chacune hyperbolique en perturbant $f$ sur des voisinages disjoints de ces orbites. Le choix de $\cU$ montre que le diff\'eomorphisme ainsi obtenu appartient \`a $\cU$. Ceci montre la densit\'e.

L'intersection de ces ouverts denses pour une suite $(\varepsilon_n)$ tendant vers $0$ est l'ensemble r\'esiduel annonc\'e.
\end{demo}

\begin{demo}[Conclusion de la d\'emonstration du th\'eor\`eme~\ref{t.conservatif} dans le cas des surfaces]
Nous avons vu qu'il existe une partie r\'esiduelle $\cG_0$ de $\diff^1_\omega(S)$ form\'ee de diff\'eomorphismes transitifs pour lesquels les orbites p\'eriodiques hyperboliques sont denses dans $S$. Puisque nous sommes en dimension $2$, les orbites p\'eriodiques hyperboliques sont des selles de m\^eme indice. Il suffit de voir que pour $f$ g\'en\'erique, elles sont deux \`a deux homocliniquement reli\'ees\footnote{Deux orbites p\'eriodiques selle de m\^eme indice sont {\em homocliniquement reli\'ees} si la vari\'et\'e stable de chacune d'entre elles coupe transversalement la vari\'et\'e instable de l'autre. Dans ce cas, leurs classes homoclines co\"\i ncident.}.

Pour $n\geq 1$, voyons que l'ensemble $O_n$ des diff\'eomorphismes dont les orbites p\'eriodiques hyperboliques de p\'eriode inf\'erieure ou \'egale $n$ sont deux \`a deux homocliniquement reli\'ees contient un ouvert dense de $\diff^1_\omega(S)$. Comme il a \'et\'e rappel\'e dans la preuve du corollaire~\ref{c.elliptiquehyperbolique}, il existe un ouvert dense $U_n\subset\diff^i_\omega(S)$ de diff\'eomorphismes dont les orbites p\'eriodiques de p\'eriode inf\'erieure \`a $n$ est fini~; le m\^eme argument montre que, quitte \`a restreindre $U_n$, ce nombre est localement constant, chacune des orbites peut \^etre suivie localement contin\^ument en restant, suivant le cas, elliptique ou hyperbolique. Pour deux orbites p\'eriodiques selles, la propri\'et\'e d'\^etre homocliniquement reli\'ees est ouverte~: on en d\'eduit que $O_n\cap U_n$ est un ouvert. Il nous reste \`a montrer la densit\'e de $O_n\cap U_n$. Pour cela fixons $f\in U_n$ et un voisinage ouvert $\cV$ de $f$ dans $U_n$ sur lequel toute orbite p\'eriodique de p\'eriode inf\'erieure \`a $n$ peut \^etre suivie contin\^ument. Soient $\gamma$ et $\sigma$ deux orbites selles de p\'eriode inf\'erieure \`a $n$ de $f$. Nous avons vu qu'il existe $\tilde f\in \cV$ arbitrairement proche de $f$ qui est transitif. Le connecting lemma de Hayashi (dans le cas conservatif) permet de cr\'eer une intersection transverse entre la vari\'et\'e instable de $\gamma$ et la vari\'et\'e stable de $\sigma$, par une petite perturbation $\hat f\in\cV$ de $\tilde f$. Une nouvelle perturbation $\bar f\in\cV$ de $\hat f$ est transitive tout en conservant cette intersection transverse. Le connecting lemma d'Hayashi permet enfin de cr\'eer l'autre intersection transverse afin de lier homocliniquement $\gamma$ et $\sigma$.  En proc\'edant ainsi pour toutes les paires de selles de p\'eriode inf\'erieure \`a $n$, on obtient la densit\'e de $O_n$. 

L'intersection $\cG_1$ des $O_n$ est donc un ensemble r\'esiduel. Pour tout $f\in\cG_0\cap\cG_1$, toutes les orbites p\'eriodiques hyperboliques sont homocliniquement reli\'ees et leur union est dense, chacune des classes homoclines co\"\i ncide donc avec $S$.
\end{demo}

\subsection{Exposants stables et d\'ecomposition domin\'ee~: preuve du th\'eor\`eme~\ref{t.herman}}

Soit $r\geq 1$ et $f$ un diff\'eomorphisme de classe $C^r$ d'une vari\'et\'e compacte $M$ qui pr\'eserve une forme volume
$\omega$. La d\'emonstration suit facilement du lemme~\ref{l.homoclinedense} suivant et du th\'eor\`eme~\ref{t.bdp} qui est une adaptation directe d'un r\'esultat de~\cite{BDP}~:

\begin{lemm}\label{l.homoclinedense}
Soit $r\geq 1$ et $f$ un diff\'eomorphisme de classe $C^r$ d'une vari\'et\'e compacte $M$ qui pr\'eserve une forme volume
$\omega$. Il existe alors une suite $(g_n)$ de diff\'eomorphismes de classe $C^r$ de $M$, convergeant en topologie $C^1$ vers $f$, et pour tout $n$ un point p\'eriodique $p_n\in M$ de $g_n$ tel que la classe homocline $H(p_n,g_n)$ soit non-triviale\footnote{La classe homocline d'une orbite p\'eriodique $p$ est dite non-triviale si elle n'est pas r\'eduite \`a l'orbite p\'eriodique elle-m\^eme.} et $\frac{1}{n}$ dense dans $M$.
\end{lemm}

\begin{theo}[\cite{BDP}, Theorem 5]\label{t.bdp} Soit $f\in \diff^r_\omega(M)$ un diff\'eomorphisme d'une vari\'et\'e compacte $M$.
Alors~:
\begin{enumerate}
\item ou bien il existe un $C^1$-voisinage $\cU$ de $f$ dans $\diff^r_\omega(M)$ et $\ell\geq 1$ tels que, pour tout dif\-f\'eo\-morphisme $g\in \cU$, toute classe homocline non-triviale $H(p,g)$ admet une d\'ecomposition $\ell$-domin\'ee~;
\item ou bien pour tout $C^1$-voisinage $\cV$ de $f$ dans $\diff^r_\omega(M)$, il existe $g\in \cV$ et un point p\'eriodique $p\in M$ pour $g$ de p\'eriode $n$ tel que $D g^n(p)=Id$.
\end{enumerate}
\end{theo}
Avant de montrer le lemme~\ref{l.homoclinedense} et d'expliquer comment obtenir cette adaptation de~\cite{BDP}, nous terminons la preuve du th\'eor\`eme~\ref{t.herman}.

\begin{demo}[D\'emonstration du th\'eor\`eme~\ref{t.herman}]
Consid\'erons $f\in \diff^r_\omega(M)$ un diff\'eomorphisme d'une vari\'et\'e compacte $M$.
D'apr\`es le th\'e\-o\-r\`e\-me~\ref{t.bdp},
\begin{itemize}
\item ou bien (d'apr\`es l'item 2) $f$ est approch\'e en topologie $C^1$ par une suite de diff\'eomorphismes $g$ de classe $C^r$ poss\'edant une orbite p\'eriodique $p$ de p\'eriode $n$ telle que $D g^n(p)=Id$~;
\item ou bien (item 1) il existe un $C^1$-voisinage $\cU$ de $f$ dans $\diff^r_\omega(M)$ et $\ell\geq 1$ tels que, pour tout diff\'eomorphisme $g\in \cU$, toute classe homocline non-triviale $H(p,g)$ admet une d\'ecomposition $\ell$-domin\'ee.
Consid\'erons des suites (donn\'ees par le lemme~\ref{l.homoclinedense})~: $(g_n)$ et $(p_n)$.
Pour $n$ assez grand, $g_n$ est dans le voisinage $\cU$ de $f$. La classe homocline $H(p_n,g_n)$ est donc $\ell$-domin\'ee. La suite $(H(p_n,g_n))$ converge vers $M$ pour la topologie de Hausdorff et la suite $(g_n)$ vers $f$ en topologie $C^1$. La $\ell$-domination passe \`a la limite (voir~\cite[Corollary 1.5]{BDP}) ce qui montre que $M$ poss\`ede une d\'ecomposition $\ell$-domin\'ee pour $f$.
\end{itemize}
Nous avons ainsi montr\'e la dichotomie annonc\'ee par le th\'eor\`eme~\ref{t.herman}.
\end{demo}

Pour montrer le lemme~\ref{l.homoclinedense}, nous utiliserons les lemmes suivants~:
\begin{lemm}\label{l.laceau} Soient $f\in \diff^r_\omega(M)$ et $\gamma$ une orbite p\'eriodique et $x$ un point de $M$.
Pour tout $\eta>0$, il existe un diff\'eomorphisme $g\in \diff^r_\omega(M)$ arbitrairement $C^1$-proche de $f$ et
une orbite p\'eriodique $\sigma$ de $g$ qui contient $\gamma\cup\{x\}$ dans son $\eta$-voisinage.
\end{lemm}
\begin{demo} Notons $m$ la p\'eriode de $\gamma$. Commen\c cons par perturber $f$ en diff\'eomorphisme $\tilde f\in \diff^r_\omega(M)$ arbitrairement $C^1$-proche de $f$
qui v\'erifie les hypoth\`eses du th\'eor\`eme~\ref{t.connect} (les orbites p\'eriodiques sont hyperboliques ou elliptiques).
Soit $y$ un point proche de $\gamma$ qui n'est pas sur l'orbite positive de $x$. Remarquons que l'on a~: $x\dashv_{\tilde f} y$ (voir le lemme~\ref{l.dashv1}). D'apr\`es le th\'eor\`eme~\ref{t.connect} et la remarque~\ref{r.connect}, il existe un diff\'eomorphisme $h\in \diff^r_\omega(M)$ arbitrairement $C^1$-proche de $\tilde f$ (et donc de $f$) tel que $y$ appartient \`a l'orbite positive de $x$. Le segment d'orbite joignant $x$ \`a $y$ \'etant de longueur plus grande que $m$ (sa longueur tend vers $\infty$ lorsque $h$ tend vers $f$), il contient les points $h^{-m}(y),\dots,y$ qui sont arbitrairement proches des points de $\gamma$. Puisque l'on a $x\dashv_h x$, il existe un diff\'eomorphisme $g\in \diff^r_\omega(M)$ arbitrairement $C^1$-proche de $h$ tel que le point $x$ est p\'eriodique. L'orbite $\sigma$ de $x$ par $g$ passe arbitrairement pr\`es de tout point du segment $x,h(x)\dots,y$ et contient donc dans son $\eta$-voisinage \`a la fois le point $x$ et l'orbite $\gamma$.
On conclut en remarquant que, par construction, le diff\'eomorphisme $g$ est arbitrairement $C^1$-proche de $f$.
\end{demo}

\begin{lemm}\label{l.intersectionhomocline} Soient $f\in \diff^r_\omega(M)$ et $\sigma$ une orbite p\'eriodique hyperbolique. Alors, il existe $g\in\diff^r_\omega(M)$ arbitrairement $C^1$-proche de $f$ tel que la classe homocline de $\sigma_g$ est non-triviale, o\`u $\sigma_g$ est la continuation hyperbolique de $\sigma$ pour $g$. 
\end{lemm}
\begin{demo} Quitte \`a perturber $f$, on peut supposer qu'il v\'erifie les hypoth\`eses du th\'eor\`eme~\ref{t.connect}.
Soit $\cU$ un petit $C^1$-voisinage de $f$ dans $\diff^r_\omega(M)$ sur lequel la continuation de $\sigma$ est bien d\'efinie. Soit $N$ l'entier associ\'e \`a ce voisinage par le th\'eor\`eme~\ref{t.connecting}. Quitte \`a r\'eduire $\cU$, on peut supposer que $N$ est sup\'erieur \`a la p\'eriode de $\sigma$. On reprend la preuve du th\'eor\`eme~\ref{t.connect}~: consid\'erons une famille de bo\^\i tes de perturbation $\cB=\cB_0\cup\bigcup_{\gamma\in Per_{N_0}(f)}(\cE(\gamma)\cup\cS(\gamma))$ pour $(f,\cU)$ comme aux sections~\ref{ss.preparation} et~\ref{s.conservatif2} (les bo\^\i tes de $\cE(\gamma)$ et de $\cS(\gamma)$ correspondent aux voisinages des orbites p\'eriodiques hyperboliques ou elliptiques de basse p\'eriode). Puisque $\sigma$ est de p\'eriode inf\'erieure ou \'egale \`a $N$, elle est disjointe des bo\^\i tes de perturbation.

On choisit deux points $x$ et $y$ de $W^u_\delta (\sigma)$ et $W^s_\delta (\sigma)$ contenus dans le voisinage $W(\sigma)$ de $\sigma$, tels que les it\'er\'es n\'egatifs de $x$ et positifs de $y$ par $f$ restent dans $W(\gamma)$. Il existe deux entiers $n_x,n_y\geq 1$ tels que les points $x_1=f^{n_x}(x)$ et $y_1=f^{-n_y}(y)$ ne soient contenus dans aucune bo\^\i te de perturbation de $\cB$ et n'appartiennent \`a aucun ouvert $V(\gamma)$. Puisque $x_1\dashv y_1$ (lemme~\ref{l.dashv1}), il existe une pseudo-orbite, joignant $x_1$ \`a $y_1$, qui pr\'eserve les quadrillages de bo\^\i tes de $\cB$ et n'ayant aucun saut en dehors de ces bo\^\i tes  (proposition~\ref{p.regroupe}). On la compl\`ete par les segments d'orbite joignant $x$ \`a $x_1$ et $y_1$ \`a $y$ pour obtenir une pseudo-orbite qui pr\'eserve les quadrillages de bo\^\i tes de $\cB$ et n'ayant aucun saut en dehors de ces bo\^\i tes joignant $x$ \`a $y$. Comme pour la section~\ref{s.fin}, on supprime bo\^\i te par bo\^\i te tous les sauts de cette pseudo-orbite et on obtient un diff\'eomorphisme $g\in \cU$ tel que $y$ est sur l'orbite positive de $x$ et co\"\i ncidant avec $f$ sur $W(\sigma)$.

Comme les it\'er\'es n\'egatifs de $x$ et positifs de $y$ par $f$  restent dans $W(\gamma)$, ils n'ont pas \'et\'e modifi\'es par la perturbation. Les points $x$ et $y$ appartiennent donc aux vari\'et\'es instable  et  stable de $\sigma_g$ respectivement. Ces vari\'et\'es s'intersectent donc le long de l'orbite de $x$. On peut rendre ces intersections transverses par une nouvelle petite perturbation. 
\end{demo}

\begin{demo}[D\'emonstration du lemme~\ref{l.homoclinedense}]
En appliquant un nombre fini de fois le lemme~\ref{l.laceau} aux points $\{x_i\}_{i\in I}$ d'une famille $\frac1{4n}$ dense dans $M$ on  montre qu'il existe un diff\'eomorphisme $h_n\in\diff^r_\omega(M)$ arbitrairement $C^1$-proche de $f$ et poss\`edant une orbite p\'eriodique $\alpha_n$ passant arbitrairement proche des points $x_i$, et donc $\frac1{3n}$-dense dans $M$. Une nouvelle $C^1$-petite perturbation  $\phi_n\in\diff^r_\omega(M)$ de $h_n$ (donn\'ee par le lemme de Franks si $dim M\geq 3$ et par la proposition~\ref{p.elliptiquehyperbolique} dans le cas des surfaces)  permet alors de cr\'eer une orbite p\'eriodique hyperbolique $\beta_n$ de $\phi_n$, contenant $\alpha_n$ dans son $\frac1{3n}$-voisinage~: en particulier, $\beta_n$ est $\frac2{3n}$-dense dans $M$.   

Pour terminer la preuve du lemme~\ref{l.homoclinedense}, il suffit \`a pr\'esent de cr\'eer, \`a l'aide du lemme~\ref{l.intersectionhomocline} une intersection homocline transverse entre les vari\'et\'es invariantes de $\beta_n$ de fa\c con que la classe homocline de cette orbite soit non-triviale. 
\end{demo}

Il reste \`a expliquer comment obtenir le th\'eor\`eme~\ref{t.bdp} \`a partir de~\cite{BDP}.
La d\'emonstration utilise deux r\'esultats. Le premier est le lemme de Franks (proposition~\ref{p.franks}). Le second est la proposition suivante de \cite{BDP} que nous traduisons dans le langage des classes homoclines.

\begin{prop}[\cite{BDP}, Proposition 7.3]  \label{p.bdp} Soient $K>0$ et  $\cU\subset \diff^1_\omega(M)$ un ouvert tels que, pour tout $f\in\cU$, les normes  de $Df$ et $Df^{-1}$ sont uniform\'ement born\'ees par $K$. Pour tout $\varepsilon>0$ il existe un entier $\ell>0$ tel que, pour tout $f\in \cU$ et toute classe homocline $H(p,f)$ non-triviale, l'une au moins des deux propri\'et\'es suivantes est v\'erifi\'ee~:
\begin{itemize}
\item  $H(p,f)$ admet une d\'ecomposition $\ell$-domin\'ee~;
\item  il existe une orbite p\'eriodique $\gamma$ homocliniquement reli\'ee \`a $p$ et une $\varepsilon$-perturbation $(B_x)_{x\in\gamma}$ de la diff\'erentielle de $f$ le long de $\gamma$ telles que le produit des applications lin\'eaires $B_x$ le long de $\gamma$ est l'identit\'e. 
\end{itemize}
\end{prop}

\begin{demo}[D\'emonstration du th\'eor\`eme~\ref{t.bdp}] Soit $f\in\diff^r_\omega(M)$ un diff\'eomorphisme ne v\'erifiant pas l'item (1) du th\'eor\`eme~\ref{t.bdp}. Fixons un $C^1$-voisinage $\cV\subset\diff^r_\omega(M)$ de $f$. Quitte \`a r\'eduire $\cV$, on peut supposer que les normes de $Dg$ et $Dg^{-1}$ sont uniform\'ement born\'ees par une constante $K$ pour $g\in\cV$. Soient $\cO$ et $\varepsilon>0$ le voisinage de $f$ et la constante donn\'es par le lemme de Franks appliqu\'e \`a l'ouvert $\cV$, et $\ell$ l'entier associ\'e \`a la constante $K$, l'ouvert $\cV$ et $\varepsilon$ par la proposition~\ref{p.bdp}.

Comme $f$ ne v\'erifie pas l'item (1) du th\'eor\`eme~\ref{t.bdp}, il existe $\tilde f\in \cO$ et une classe homocline $H(p,\tilde f)$ qui ne poss\`ede aucune d\'ecomposition $\ell$-domin\'ee. D'apr\`es la proposition~\ref{p.bdp}, il existe donc une orbite p\'eriodique $\gamma$ de $\tilde f$ et une $\varepsilon$-perturbation $(B_x)_{x\in \gamma}$ de la diff\'erentielle de $\tilde f$ le long de $\gamma$ telle que le produit des applications lin\'eaires $B_x$ le long de $\gamma$ soit l'identit\'e. Finalement, comme $\tilde f$ appartient \`a $\cO$, le lemme de Franks permet de r\'ealiser la perturbation $(B_x)$ des diff\'erentielles par une perturbation $g\in\cV$ de $\tilde f$~: l'ensemble $\gamma$ est une orbite p\'eriodique de $g$, dont la diff\'erentielle  \`a la p\'eriode est l'identit\'e.  
\end{demo}

\appendix
\section{Le Connecting Lemma}\label{a.A}
Dans cet appendice, nous expliquons comment montrer le th\'eor\`eme~\ref{t.connecting}. Cet \'enonc\'e n'existant pas sous cette forme dans la litt\'erature (et puisqu'il est difficile d'expliquer au lecteur comment modifier certaines hypoth\`eses dans des preuves longues et techniques), nous avons choisi d'en redonner une preuve compl\`ete qui est cependant tr\`es proche de~\cite[pages 359--368]{Ar}.

\subsection{Cubes perturbatifs uniformes}\label{a.cube}

L'une des \'etapes clef de la d\'emonstration est un r\'esultat perturbatif (que \cite{Ar} attribue pour l'essentiel \`a \cite{Pu,PuRo}).

\'Etant donn\'ee une carte $\varphi\colon V\to \RR^d$, nous appellerons {\em cube} de la carte tout cube $C$ plong\'e dans $V$ tel que $\varphi(C)$ soit obtenu \`a partir du cube standard $[-1,1]^d$ par une homoth\'etie-translation. Nous noterons alors $(1+\varepsilon) C$ le cube de m\^eme centre qui lui est homoth\'etique de rapport $(1+\varepsilon)$ (dans les coordonn\'ees donn\'ees par $\varphi$).
Pour un cube $C$ de la carte $\varphi$, le cube $(1+\varepsilon)C$ est bien d\'efini si $\varepsilon$ est suffisamment petit. Par la suite (en particulier lorsque nous consid\'ererons les cubes perturbatifs d\'efinis ci-apr\`es) nous supposerons implicitement que le cube $3C$ est encore un cube de la carte $\varphi$.

Soit $\cU$ un $C^1$-voisinage d'un diff\'eomorphisme $f$. Fixons $N\in\NN$ et des constantes $\varepsilon,\eta\in]0,1[$. On dira que $C$ est un {\em cube $(f,\cU,\varphi,N,\varepsilon,\eta)$-perturbatif\/} si pour toute paire $(p,q)$ de points du cube $C$ il existe $g\in\cU$ qui  a les propri\'et\'es suivantes:

\begin{itemize}
\item $g$ co\"\i ncide avec $f$ hors de $\bigcup_{t=0}^{N-1}(f^t((1+2\varepsilon)C))$~;
\item $g^N(p)=f^N(q)$~;
\item pour tout $t\in\{0,dots,N-1\}$, $g^t(p)$ appartient \`a $f^t((1+\varepsilon) C)$~;
\item pour tout $t\in\{0,dots,N-1\}$, $g$ ne diff\`ere de $f$ sur $f^{t}((1+2\varepsilon)C)$ que sur une boule centr\'ee en  $g^t(p\in f^t((1+\varepsilon) C)$ et  de rayon inf\'erieur \`a $\eta$ fois la distance entre  $f^t((1+\varepsilon) C)$ et le compl\'ementaire de $f^{t}((1+2\varepsilon)C)$.  
\end{itemize}

\begin{theo}[\cite{Ar}, th\'eor\`eme 22 et son addendum]
Soit $f$ un diff\'eomorphisme d'une vari\'et\'e compacte $M$ et $\cU$ un $C^1$-voisinage de $f$. Fixons deux nombres $\varepsilon,\eta \in ]0,1[$. 

Pour tout point $p_0\in M$ non-p\'eriodique, il existe un entier $N$ et,  au voisinage de $p_0$, une carte $\varphi\colon V\to \RR^d$ telle que  tout cube $C$ de $V$ est un  cube $(f,\cU,\varphi,N,\varepsilon,\eta)$-perturbatif.
\label{t.rp}
\end{theo}

Pour la version annonc\'ee du connecting lemma que nous utilisons ici, nous avons besoin d'une version tr\`es l\'eg\`erement plus forte de cet \'enonc\'e, mais dont la preuve est, en fait, contenue dans celle donn\'ee dans \cite{Ar}. Plus pr\'ecis\'ement nous avons besoin de l'uniformit\'e de l'entier $N$ et d'une formulation qui soit aussi valable au voisinage des points p\'eriodiques.

\begin{theo} Soit $f$ un diff\'eomorphisme d'une vari\'et\'e compacte $M$ et $\cU$ un $C^1$-voisinage de $f$. Fixons deux nombres $\varepsilon,\eta \in ]0,1[$. 

Il existe un entier $N$ et,  au voisinage de tout point  $p_0\in M$, une carte $\varphi\colon V\to \RR^d$ v\'erifiant la propri\'et\'e suivante~:

Tout cube $C$ de $V$ tel que le cube $(1+2\varepsilon)C$ est disjoint de ses $N$ premiers it\'er\'es est un  cube $(f,\cU,\varphi,N,\varepsilon,\eta)$-perturbatif.
\label{t.rpu}
\end{theo}

Appelons {\em forme de cube de $\RR^d$} la donn\'ee d'une base orthogonale (pas n\'ecessairement orthonorm\'ee) $b$ de $\RR^d$. Le cube standard $C(b)$ de forme $b$ sera l'image du cube standard de $\RR^d$ par l'application lin\'eaire envoyant la base standard de $\RR^d$ sur $b$.

Le th\'eor\`eme~\ref{t.rp}  s'appuie essentiellement sur la proposition suivante (due initialement \`a \cite{PuRo})~:

\begin{prop}[voir~\cite{Ar}, proposition 26] Soit $(T_t)_{t\in\NN}$ une suite dans $GL(\RR,d)$. Fixons $\varepsilon,\eta\in]0,1[$. Il existe $N$ et il existe une forme de cube $b$ de $\RR^d$ tels que pour toute paire $(p,q)$ de points du cube $C(b)$, il existe une suite $(x_t)_{t\in\{0,\dots,N\}}$ de points de $(1+\varepsilon)C(b)$ telle que $x_0=p$ et $x_N=q$ et, pour tout $t\in\{1,\dots,N\}$, la distance $d(T_t(x_{t-1}),T_t(x_t))$ est inf\'erieure \`a $\eta$ fois la distance entre $T_t((\varepsilon+1) C(b))$ et le compl\'ementaire de $T_t((1+2\varepsilon)C(b))$.  
\label{p.rp}
\end{prop} 

Dans la d\'emonstration du th\'eor\`eme~\ref{t.rp}, la suite $(T_t)$ de la proposition~\ref{p.rp}, correspond \`a la suite $(Df^t(p_0))$ exprim\'ee dans des bases orthonorm\'ees de $T_{p_0}M$ et $T_{f^t(P_0)}M$. 

Remarquons que quitte \`a modifier l\'eg\'erement les constantes $\varepsilon$ et $\eta$, la conclusion de la proposition reste v\'erifi\'ee par des suites $(\tilde T_t)$ proches de $(T_t)$ pour la topologie produit (il suffit bien s\^ur que les $N$ premiers termes $\tilde T_t$ soient proches des $T_t$ correspondants).  

Le th\'eor\`eme~\ref{t.rpu} a la m\^eme d\'emonstration que le th\'eor\`eme~\ref{t.rp},  \`a partir d'une version uniforme de la proposition~\ref{p.rp}~:

\begin{prop} Soient $K>0$ et $\varepsilon,\eta\in]0,1[$ fix\'es. Il existe $N$ tel que, pour toute suite $(T_t)_{t\in\{0,\dots N\}}$ avec $T_t\in GL(\RR,d)$ v\'erifiant $\|T_t\circ T_{t-1}^{-1}\|\leq K$ et $\|T_{t-1}\circ T_t{-1}\|\leq K$ on a la propri\'et\'e suivante~:

il existe une forme de cube $b$ de $\RR^d$ telle que, pour toute paire $(p,q)$ de points du cube $C(b)$, il existe une suite $(x_t)_{t\in\{0,\dots,N\}}$ de points de $(1+\varepsilon)C(b)$ telle que $x_0=p$ et $x_N=q$ et, pour tout $t\in\{0,\dots,N-1\}$, la distance $d(f^t(x_{t}),f^t(x_{t+1}))$ est inf\'erieure \`a $\eta$ fois la distance entre  $f^t((1+\varepsilon) C(b))$ et le compl\'ementaire de $f^t((1+2\varepsilon)C(b))$.
\label{p.rpu}
\end{prop} 

(Cette version uniforme se d\'eduit de la proposition~\ref{p.rp} par compacit\'e de l'ensemble des suites de matrices born\'ees et d'inverses born\'es par $K$). 

On en d\'eduit la proposition suivante qui avec le lemme~\ref{l.perturb} \'enonc\'e ci-dessous implique le th\'eor\`eme~\ref{t.rpu}. C'est cet \'enonc\'e que nous utiliserons finalement par la suite.

\begin{prop}\label{p.suite} Soit $f$ un diff\'eomorphisme d'une vari\'et\'e compacte $M$ et $\cU$ un $C^1$-voisinage de $f$. Fixons deux nombres $\varepsilon,\eta \in ]0,1[$. Il existe un entier $N$ et, au voisinage de tout point $p_0\in M$, une carte $\varphi\colon V\to \RR^d$ telle que pour tout cube $C$ de $V$ disjoint de ses $N$ premiers it\'er\'es et toute paire $(a,b)$ de points de $C$, il existe une suite $(a_t)_{t\in\{0,\dots,N\}}$ de points de $(1+\varepsilon)C$ telle que $a_0=a$ et $a_N=b$ et, pour tout $t\in\{0,\dots,N-1\}$, la distance $d(f^t(a_{t}),f^t(a_{t+1}))$ est inf\'erieure \`a $\eta$ fois la distance entre  entre  $f^t((\varepsilon+1) C)$ et le compl\'ementaire de $f^t((1+2\varepsilon)C)$.
\end{prop}

Nous utiliserons \'egalement un lemme de perturbation classique.

\begin{lemm}[Voir~\cite{Ar2}, proposition 5.1.1]\label{l.perturb}
Soit $M$ une vari\'et\'e compacte $M$ et $\cV$ un voisinage de l'identit\'e dans $\diff^1(M)$. Il existe des constantes $\lambda>1$ et $\delta>0$ telles que pour tous points $p,q$ de $M$ satisfaisant $d(p,q)<\delta$, il existe une perturbation $h$ de l'identit\'e, \`a support dans la boule de centre $p$ et de rayon $\lambda d(p,q)$, telle que $h\circ f$ appartienne \`a $\cU$ et $h(p)=q$.
\end{lemm}
\begin{rema}\label{r.perturb} La $C^1$-perturbation $h$ de l'identit\'e  peut \^etre choisie de classe $C^{\infty}$. De plus si $M$ est munie d'une forme volume $\omega$, alors  on peut imposer que $h$ pr\'eserve $\omega$ .
\end{rema}

\subsection{Fin de la preuve du connecting lemma}\label{s.connecting}

On commence par choisir quelques constantes. Nous allons devoir appliquer la proposition~\ref{p.suite} aux carreaux d'un quadrillage. Chaque perturbation peut d\'eborder d'un carreau. Nous aurons donc besoin de majorer le nombre de carreaux adjacents \`a un carreau donn\'e. Il existe une borne uniforme $N^*=4^d$ o\`u $d$ est la dimension de $M$.

Les supports des perturbations associ\'ees \`a des carreaux ne pourront se chevaucher que si les carreaux sont adjacents. Pour cela nous choisissons la constante $\varepsilon$ de la proposition~\ref{p.suite} inf\'erieure \`a $\frac{1}{10}$~: pour tous carreaux $C_1,C_2$ du quadrillage standard, si $(1+2\varepsilon)C_1\cap(1+2\varepsilon)C_2\neq\emptyset$, alors $C_1\cap C_2\neq \emptyset$.

Au voisinage $\cU$ de $f$, on peut associer un voisinage $\cV$ de l'identit\'e dans $\diff^1(M)$ avec la propri\'et\'e suivante~: pour toute famille finie $\{\psi_1,\dots,\psi_k\}$ de perturbations de l'identit\'e, contenues dans $\cV$ et de supports deux \`a deux disjoints, le diff\'eomorphisme $g$, qui co\"\i ncide avec $f$ hors de l'union des supports des $\psi_i$ et avec $f\circ \psi_i$ sur le support de $\psi_i$, pour tout $i\in\{1,\cdots,k\}$, appartient \`a $\cU$.
Le lemme~\ref{l.perturb} associe au voisinage $\cV$ des constantes $\lambda>1$ et $\delta>0$. Nous fixons ensuite la constante $\eta>0$ de la proposition~\ref{p.suite} pour que

$$\frac{1}{\eta}>(6\lambda N^*)^{N^*+1}.$$
Ce choix sera justifi\'e par la proposition~\ref{p.boules}.

On applique alors la proposition~\ref{p.suite}. Ceci d\'efinit l'entier $N$. Pour tout point $x=p_0$ de $M$, on obtient aussi une carte locale $\varphi : V\to \RR^d$ en $x$. La carte locale $\varphi : U_x\to \RR^d$ est la restriction de $\varphi$ \`a un voisinage $U_x$ de $x$ sur lequel les applications $f,\dots,f^N$ sont proches des diff\'erentielles $D f(x),\dots,D f^N(x)$.

\begin{rema}\label{r.alpha}
%
%

On choisit $U_x$ assez petit pour que l'on ait de plus que le diam\`etre de $U_x$ et de ses $N$ premiers it\'er\'es soit major\'e par  la constante  $\delta>0$ du lemme~\ref{l.perturb}. 
\end{rema}

Soit $\cC$ un cube quadrill\'e de $(U_x,\varphi)$ disjoint de ses $N$ premiers it\'er\'es. Pour voir que c'est une bo\^\i te de perturbation pour $(f,\cU)$, nous consid\'erons jusqu'\`a la fin de cette section, une suite de paires de points $\{(x_i,y_i)\}_{i\in \{1,\dots,\ell\}}$ de $\cC$ telle que pour tout $i\in \{1,\dots,\ell\}$, les points $x_i$ et $y_i$ appartiennent au m\^eme carreau de $\cC$.

\begin{lemm}\label{l.boules}
Il existe une sous-suite $i_0=0,\dots,i_m$, telle que pour tout $n\in\{0,\dots,m-1\}$ les points $a^n=x_{i_n}$ et $b^n=y_{i_{n+1}-1}$ appartiennent \`a un m\^eme carreau $C_n$. De m\^eme, les points $a^m=x_{i_m}$ et $b^m=y_\ell$ appartiennent \`a un m\^eme carreau $C_m$. De plus, si $n\neq n'$ alors $C_n\neq C_{n'}$.
\end{lemm}
\begin{demo} On construit les suites $(i_n)$ et $(C_n)$ par r\'ecurrence sur $n$. On pose $i_0=0$ et on appelle $C_0$ le carreau contenant $x_{i_0}$. Soit $k_0\geq i_0$ le plus grand indice tel que les points $x_{k_0}$ et $y_{k_0}$ appartiennent \`a $C_0$. On pose $i_1=k_0+1$ et on appelle $C_1$ le carreau contenant $x_{i_1}$ et $y_{i_1}$. N\'ecessairement, $C_1$ est diff\'erent de $C_0$.

Supposons maintenant que $i_n$ et $C_n$ ont \'et\'e d\'efinis. Soit $k_n$ le plus grand indice tel que les points $x_{k_n}$ et $y_{k_n}$ appartiennent \`a $C_n$. On pose $i_{n+1}=k_n+1$ et on appelle $C_{n+1}$ les carreaux contenant $x_{i_n}$ et $y_{i_n}$. N\'ecessairement, $C_{n+1}$ est diff\'erent des carreaux $C_0,\dots,C_n$.
\end{demo}

D'apr\`es la proposition~\ref{p.suite}, pour tout $n\in \{0,\dots,m\}$, il existe une suite $a^n_0=a^n,\dots,a^n_N=b^n$ de points dans $(1+\varepsilon) C_n$ telle que pour tout $t\in\{0,\dots, N-1\}$, la distance $d(f^t(a^n_{t}),f^t(a^n_{t+1}))$ est inf\'erieure \`a $\eta$ fois la distance entre  $f^t((1+\varepsilon) C_n)$ et le compl\'ementaire de $f^t((1+2\varepsilon)C_n)$.

L'id\'ee na\"\i ve (qui ne marche pas) est de consid\'erer des perturbations $\varphi^n_t$ \`a support dans $f^t(\cC)$ telles que $\varphi^n_t(f^t(a^n_t)))=f^t(a^n_{t+1})$ et de remplacer, sur ces supports, $f$ par $f\circ \varphi^n_t$. Mais en g\'en\'eral, les supports des perturbation $\varphi^n_t$ s'intersectent. Remarquons que si les supports de $\varphi^n_t$ et $\varphi^{n'}_t$ avec $n<n'$ s'intersectent, les points $f^t(a^n_t)$ et $f^t(a^{n'}_{t+1})$ sont proches. Dans ce cas, il semble avantageux de consid\'erer une nouvelle perturbation \`a support dans $f^t(\cC)$ envoyant $f^t(a^n_t)$ sur $f^t(a^{n'}_{t+1})$ et d'oublier les points interm\'ediaires entre $a^n_t$ et $a^{n'}_{t+1}$ dans la suite~:

$$a^0_0=a^0,\dots,a^0_t,\dots,a^0_N,a^1_0,\dots,a^1_N,\dots,a^m_0,\dots,a^m_N=b^m.$$

Le ph\'enom\`ene pouvant se r\'epetter plusieurs fois, le lemme suivant produit une nouvelle suite plus courte pour laquelle les supports des perturbations sont deux \`a deux disjoints.

\begin{lemm}\label{l.suite2}
Soit $a^0_0=a^0,\dots,a^0_t,\dots,a^0_N,a^1_0,\dots,a^1_N,\dots,a^m_0,\dots,a^m_N=b^m$ une suite telle que pour tout $n\in \{0,\dots,m\}$, les points $a^n_0,\dots,a^n_N$ appartiennent \`a un m\^eme cube $(1+\varepsilon)C_n$, les cubes $C_n$ \'etant deux \`a deux distincts, et telle que pour tout $n\in \{0,\dots,m\}$ et $t\in\{0,\dots, N-1\}$, la distance $d(f^t(a^n_{t}),f^t(a^n_{t+1}))$ est inf\'erieure \`a $\eta$ fois la distance entre  $f^t((1+\varepsilon) C_n)$ et le compl\'ementaire de $f^t((1+2\varepsilon)C_n)$.

Alors, il existe $s\in \{0,\dots,m\}$, et une suite croissante (non-n\'ecessairement strictement croissante) $1=\nu(0),\dots,\nu(s (N+1)+N)=m$ telle que~:

\begin{enumerate}
\item Pour tout $k\in\{1,\dots,s\}$, on a $\nu(k (N+1))=\nu(k (N+1)-1)+1$.
\item Pour tout $t\in\{0,\dots,N-1\}$, les boules $B_t^{k}$ centr\'ees aux points $f^t(a^{\nu(k (N+1)+t)}_t)$, de rayon $\lambda d(f^t(a^{\nu(k (N+1)+t)}_t),f^t(a^{\nu(k (N+1)+t+1)}_{t+1}))$, sont deux \`a deux disjointes, pour $k\in\{0,\dots,s\}$.
\item Pour tout $t\in\{0,\dots,N-1\}$, les boules $B_t^{k}$  sont contenues dans $f^t(\cC)$.
\end{enumerate}
\end{lemm}

Nous pouvons maintenant terminer la preuve du th\'eor\`eme~\ref{t.connecting}. La d\'emonstration du lemme \ref{l.suite2} est report\'ee \`a la section~\ref{s.suite2}. Elle utilise la proposition~\ref{p.boules} de la section~\ref{s.boules}

\begin{demo}[D\'emonstration du th\'eor\`eme~\ref{t.connecting}]
De la suite $(x_i,y_i)_{i\in\{1,\dots,\ell\}}$, le lemme~\ref{l.boules} a extrait deux sous-suites $(a^n)_{n\in\{0,\dots,m\}}$ et $(b^n)_{n\in\{0,\dots,m\}}$ qui nous ont permis d'obtenir la suite $a^0_0=a^0,\dots,a^0_t,\dots,a^0_N,a^1_0,\dots,a^1_N,\dots,a^m_0,\dots,a^m_N=b^\ell$. Nous appliquons ensuite le lemme~\ref{l.suite2} qui en extrait la suite $$a^{\nu(0)}_0=a^0=x_1,a^{\nu(1)}_1,\dots,a^{\nu(N)}_N,a^{\nu(N+1)}_0,\dots,a^{\nu(2N+1)}_N,\dots,a^{\nu((s)(N+1))}_0,\dots,a^{\nu(s(N+1)+N)}_N=y_\ell.$$

L'item 3 du lemme~\ref{l.suite2} et la remarque~\ref{r.alpha} entra\^\i nent que les boules $B^k_t$ sont de diam\`etre inf\'erieur \`a $\delta$.  Gr\^ace au lemme~\ref{l.perturb}, il existe alors des diff\'eomorphismes $\psi_t^{k}$ avec $t\in\{0,\dots,N-1\}$ et $k\in\{0,\dots,s\}$, appartenant \`a $\cV$ et \`a support dans la boule $B^k_t$ d\'efinie au lemme~\ref{l.suite2} tels que $\psi^k_t(f^t(a^{\nu(k (N+1)+t)}_t))=f^t(a^{\nu(k (N+1)+t+1)}_{t+1})$. D'apr\`es le lemme~\ref{l.suite2}, les supports des $\psi_t^{k}$ sont donc deux \`a deux disjoints. Notons $g$ le diff\'eomorphisme co\"\i ncidant avec $f$ hors de l'union des supports des $\psi_t^{k}$ et avec $f\circ \psi_t^k$ sur chaque $B^k_t$. Par le choix de $\cV$ et de $\cU$, le diff\'eomorphisme $g$ appartient \`a $\cU$. Par construction, pour tout $k\in\{0,\dots,s\}$ et $t\in\{0,\dots,N-1\}$, on a $g^t(a_0^{\nu(k(N+1))})=f^t(a_t^{\nu(k(N+1)+t)})$.

Voyons que $g$ est la perturbation annonc\'ee~: on rappelle que $a^{\nu(k(N+1))}_0=a^{\nu(k (N+1))}=x_{i_{n_k}}$ et $a^{\nu(k(N+1)+N)}_N=b^{\nu(k (N+1)+N)}=y_{i_{m_k}}$ avec $n_k=i_{\nu(k (N+1))}$ et $m_k=i_{\nu(k (N+1)+N)+1}-1$.

Montrons $n_{k+1}=m_k+1$~: d'apr\`es l'item 1 du lemme~\ref{l.suite2},
$$n_{k+1}= i_{\nu((k+1) (N+1))}=i_{\nu((k+1) (N+1)-1)+1}=i_{\nu(k (N+1)+N)+1}=m_k+1.$$
Les couples $(x_{n_k},y_{m_k})_{k\in\{0,\dots,s\}}$ et la suite $n_0=1<n_1<\cdots<n_s\leq \ell$ v\'erifient donc~:

\begin{enumerate}
\item $x_{n_0}=x_1$ et $y_{m_{s}}=y_\ell$. 
\item $g^N(x_{n_r})=f^N(y_{m_r})=f^N(y_{n_{r+1}-1})$ pour $r\neq s$ et $g^N(x_{n_{s}})=f^N(y_\ell)$.
\end{enumerate}
Ceci montre que $\cC$ est une bo\^\i te de perturbation.
\end{demo}

\subsection{Regroupement de points proximaux}\label{s.boules}
Nous \'enon\c cons un r\'esultat g\'en\'eral de regroupement de points qui sera utilis\'e pour la d\'emons\-tra\-tion du lemme~\ref{l.suite2}.

\begin{prop} Soit $(X,d)$ un espace m\'etrique, $\{x_i\}_{i\in I}$ un ensemble de points et pour chaque $i\in I$ un r\'eel $r_i>0$. On suppose qu'il existe des constantes $K>0$, $\lambda>1$ et un entier $N^*$ tels que 
\begin{itemize}
\item pour tout $i\in I$, le cardinal de l'ensemble $\{i'\in I, B(x_i,Kr_i)\cap B(x_{i'},Kr_{i'})\neq \emptyset\}$ est major\'e par $N^*$.
\item $(6\lambda N^*)^{N^*+1}<K$
\end{itemize}
 Alors, il existe une partition de $I$ en classes $\cE=\{E_j, j\in J\}$, telle que l'on ait la propri\'et\'e suivante~:

Pour tout $j\in J$, notons $\delta_{1,j}$ le diam\`etre de l'ensemble $\{x_i, i\in E_j\}$, et $\delta_{2,j}=\sup \{r_i, i\in E_j\}$.  Finalement notons $\delta_j=\delta_{1,j}+\delta_{2,j}$. Soit $B_j$ l'union des boules $B(x_i, \lambda \delta_j)$, pour $i\in E_j$. Alors les $B_j$ sont deux \`a deux disjoints.

De plus pour tout $j\in J$, il existe $i_0\in E_j$ tel que l'ensemble $B_j$ est contenu dans la boule $B(x_{i_0},Kr_{i_0})$. 
\label{p.boules}
\end{prop}

\begin{demo} Nous allons d\'efinir par r\'ecurrence sur $k$ des partitions $\cE^k=\{E_j^k, j\in J^k\}$  de $I$. Ceci nous permet d'introduire, pour tout $j\in J^k$,   le diam\`etre $\delta^k_{1,j}$ de l'ensemble $\{x_i, i\in E^k_j\}$. On notera aussi  $\delta^k_{2,j}=\sup \{r_i, i\in E_j\}$ et $\delta^k_j=\delta^k_{1,j}+\delta^k_{2,j}$. Finalement on notera  $B^k_j$ l'union des boules $B(x_i, \lambda\delta^k_j)$, pour $i\in E^k_j$.

La partition $\cE^0$ est la partition en singleton $\left\{\{j\},j\in I\right\}$, c'est-\`a-dire que $J^0=I$. On a donc $\delta^0_j=r_j$, pour tout $j\in J^0=I$.

Supposons la partition $\cE^k$ construite. Nous dirons que $x_{i_1}$ et $x_{i_2}$ sont $k+1$-adjacents s'il existe une suite finie $j_1,\dots, j_\ell\in J_k$ telle que $x_{i_1}\in E^k_{j_1}$, $x_{i_2}\in E^k_{j_\ell}$ et, pour tout $s\in\{1,\dots,\ell-1\}$,  on ait  $B^k_{j_s}\cap B^k_{j_{s+1}}\neq \emptyset$. La relation de $k+1$-adjacence est une relation d'\'equivalence, qui est  clairement plus grossi\`ere que la $k$-adjacence~: tout classe de $k$-adjacence est incluse dans une classe de $k+1$-adjacence. 

On  note $\cE^{k+1}$ la partition de $I$ induite par la relation de $k+1$-adjacence. 
La proposition~\ref{p.boules} est une cons\'equence des lemmes suivants~: le lemme~\ref{l.boules5} montre que la partition $\cE^{N^*}$ produit des ensembles $(B^{N^*}_j)_{j\in J^{N^*}}$ disjoints. Les lemmes~\ref{l.boules1} et~\ref{l.boules3} entra\^\i nent, pour tout $j\in J^{N^*}$, l'existence d'un indice $i_0\in E^{N^*}_j$ tel que $B_j^{N^*}$ est contenu dans $B(x_{i_0},Kr_{i_0})$.
\end{demo}

\begin{lemm}\label{l.boules1}
Soient $k\in\{0,\cdots,N^*\}$ et $j\in J^k$. Supposons que l'on a
$$\delta^k_{1,j}\leq \left( 6\lambda N^*\right)^k\delta^k_{2,j},$$
et que le cardinal de $E^k_j$ est fini. Soit $i_0\in E_j^k$ tel que
$r_{i_0}=\delta^k_{2,j}$ (ceci est possible puisque $\delta^k_{2,j}$ est la borne sup\'erieur d'un ensemble fini). 

Alors, le diam\`etre  de $B^k_j$ est major\'e par $6\lambda \left(6\lambda N^*\right)^k r_{i_0}$. En cons\'equence $B^k_j$ est inclus dans la boule $B(x_{i_0}, Kr_{i_0})$.
\end{lemm}
\begin{demo}
En effet si $i\in E_j^k$ on a $d(x_i,x_{i_0})\leq \delta^k_{1,j}$.  

La boule $B(x_i, \lambda\delta^k_j)$ est donc incluse dans la boule $B(x_{i_0}, \delta^k_{1,j}+ \lambda\delta^k_j)$. Or,
$$\delta^k_{1,j}+ \lambda\delta^k_j = (1+\lambda)\delta_{1,j}^k +\lambda \delta_{2,j}^k\leq \left[(1+\lambda)\left(6\lambda N^*\right)^k+\lambda \right]\delta_{2,j}^k\leq \left[(1+\lambda)\left(6\lambda N^*\right)^k +\lambda\right]r_{i_0}.$$
En particulier $B^k_j$ est inclus dans la boule $B\left(x_{i_0},\left[(1+\lambda)\left(6\lambda N^*\right)^k +\lambda\right]r_{i_0}\right)$ et donc~:
$$diam(B^k_j)\leq 2\left[(1+\lambda)\left(6\lambda N^*\right)^k +\lambda\right]r_{i_0}.$$

En remarquant que $1<\lambda$ et que $1<\left(6\lambda N^*\right)^k$ on peut majorer $2\left[(1+\lambda)\left(6\lambda N^*\right)^k +\lambda\right]r_{i_0}$ par $6\lambda\left(6\lambda N^*\right)^k r_{i_0}$.

D'autre part notre choix de $K$ implique $K> 6\lambda\left(6\lambda N^*\right)^k$. Par cons\'equent, $\delta_{1,j}^k+\lambda\delta_j^k\leq K r_{i_0}$ ce qui conclut la preuve du lemme.
\end{demo}

\begin{lemm}\label{l.boules2}
Soit $k\in\{0,\cdots,N^*\}$. Supposons que pour tout $j\in J^k$, on a
$$\delta^k_{1,j}\leq \left( 6\lambda N^*\right)^k\delta^k_{2,j},$$
et que le cardinal de $E^k_j$ est fini.

Alors, pour tout $j\in J^{k+1}$, le cardinal de $E^{k+1}_j$ est major\'e par $N^*$.

De plus, consid\'erons $i_0\in E^{k+1}_j$ tel que $\delta_{j,2}^{k+1}=r_{i_0}$, alors pour tout $i\in E^{k+1}_j$ on a~:
$$ x_{i}\in B(x_{i_0},Kr_{i_0}).$$
\end{lemm}
\begin{demo}
Soit $I'$ une partie finie de $I$ contenue dans un $E\in\cE^{k+1}$. 

\begin{affi}Il existe $i_0\in E$ tel que pour tout $i\in I'$, le point $x_{i}$ est contenu dans $B(x_{i_0},Kr_{i_0})$
\end{affi}
\begin{demo} Remarquons d'abord qu'il suffit de montrer cette affirmation pour une partie finie  de $E$  contenant $I'$. Puisque les $E_j^k$ sont des ensembles finis, quitte \`a rajouter \`a $I'$ un nombre fini d'\'el\'ements de $E$, on peut supposer que~:
\begin{itemize}
\item $I'$ est une r\'eunion finie de classes $(E^k_j)_{j\in J'}$ de $\cE^k$~;
\item pour tous $i_1$ et $i_2$ de $I'$, il existe une suite $j_1,\dots, j_s$ de $J'$ telle que, pour tout $r\in\{1,\dots, s-1\}$, on ait $B^{k}_{j_r}\cap B^k_{j_{r+1}}\neq\emptyset$  et telle que ${i_1}\in E^k_{j_1}$ et ${i_2}\in E^k_{j_s}$.
\end{itemize}

On choisit maintenant l'indice $i_0\in I'$ de fa\c con \`a ce que $r_{i_0}$ soit maximal parmi les $r_i$, avec $i\in I'$. 
Nous allons montrer par l'absurde que l'affirmation est vraie avec ce choix de $i_0$. Supposons donc par l'absurde qu'il existe $i\in I'$ tel que $x_i$ n'appartienne pas \`a $B(x_{i_0},Kr_{i_0})$. On fixe une suite $j_1,\dots, j_s$ de $J'$ telle que, pour tout $r\in\{1,\dots, s-1\}$, on ait $B^{k}_{j_r}\cap B^k_{j_{r+1}}\neq\emptyset$  et telle que ${i_0}\in E^k_{j_1}$ et ${i}\in E^k_{j_s}$.

Alors il existe un plus petit entier $r$ tel que $B^k_{j_{r}}$ n'est pas inclus dans $B(x_{i_0},Kr_{i_0})$. Par cons\'equent  
$B^k_{j_{r}}\cap B(x_{i_0},Kr_{i_0})\neq \emptyset$~: en effet, si $r=1$ le point $x_{i_0}$ appartient \`a $B^k_{j_1}$ par d\'efinition, et si $r>1$ alors $B^k_{j_{r-1}}\subset B(x_{i_0},Kr_{i_0})$, par d\'efinition de $r$, et $B^{k}_{j_r}\cap B^k_{j_{r-1}}\neq\emptyset$, par choix de la suite $j_1,\dots,j_s$. 

Remarquons que $r\leq N^*$: en effet l'union des $E^k_{j_t}$ pour $t\leq r$ est de cardinal inf\'erieur \`a $N^*$, car pour tout $t$, l'ensemble $B^k_{j_t}$ intersecte la boule $B(x_{i_0},Kr_{i_0})$~; d'apr\`es le lemme~\ref{l.boules1}, il existe donc un point $x_{i'}$ avec $i'\in E^k_{j_t}$ tel que les boules $B(x_{i'},Kr_{i'})$ et $B(x_{i_0},Kr_{i_0})$ se rencontrent mais ceci ne peut arriver plus de $N^*$ fois d'apr\`es les hypoth\`eses de la proposition~\ref{p.boules}.

De plus le diam\`etre de chaque $B^k_{j_t}$, pour $t\leq r$ est major\'e par $6\lambda\left(6\lambda N^*\right)^k r_{i_0}$ d'apr\`es le lemme~\ref{l.boules1} et par d\'efinition de $i_0$.

De l'hypoth\`ese $B^{k}_{j_r}\cap B^k_{j_{r+1}}\neq\emptyset$, on d\'eduit alors que le diam\`etre de l'union des $B^k_{j_t}$, pour $t\leq r$, est major\'e par $6\lambda N^*\left(6\lambda N^*\right)^k r_{i_0}$ et donc inf\'erieur \`a $Kr_{i_0}$. Ceci montre $B^k_{j_r}\subset B(x_{i_0},Kr_{i_0})$ contredisant la d\'efinition de $r$, et donc  notre hypoth\`ese par l'absurde.
\end{demo}

On en d\'eduit finalement que le cardinal de $I'$ est major\'e par $N^*$. Ceci implique que $E$ lui-m\^eme est fini et de cardinal major\'e par $N^*$. On peut donc appliquer le raisonnement pr\'ec\'edent \`a l'ensemble $I'=E$ et \`a tout indice $i_0$ tel que $r_{i_0}=\sup \{r_i,i\in E\}$.  On montre ainsi que $\{x_i,i\in E\}$ est contenu dans $B(x_{i_0},Kr_{i_0})$, ce qui termine la d\'emonstration du lemme.
\end{demo}

\begin{lemm} Pour tout $k\in\{0,\dots,N^*+1\}$ et tout $j\in J^k$, on a
$$\delta^k_{1,j}\leq \left(6\lambda N^*\right)^k\delta^k_{2,j},$$
et le cardinal de $E^k_j$ est born\'e par $N^*$.
\label{l.boules3}
\end{lemm}
\begin{demo} La preuve se fait par r\'ecurrence sur $k$, sachant qu'elle est v\'erifi\'ee pour $k=0$. Supposons donc la propri\'et\'e montr\'ee pour $k\leq N^*$. 

Soit $j\in J^{k+1}$. D'apr\`es le lemme~\ref{l.boules2}, le cardinal de $E_j^{k+1}$ est major\'e par $N^*$. Soit $i_0\in E_j^{k+1}$ tel que $\delta^{k+1}_{2,j}=r_{i_0}$. Le lemme~\ref{l.boules1} montre que pour tout ensemble $E_{j'}^k\in \cE^k$ contenu dans $E_j^{k+1}$, le diam\`etre de $B^k_{j'}$ est major\'e par $$diam(B^k_{j'})\leq 6\lambda\left(6\lambda N^*\right)^k\delta^k_{2,j'}\leq 6\lambda\left(6\lambda N^*\right)^kr_{i_0}.$$
L'ensemble $E_j^{k+1}$ est une classe de $k+1$-adjacence qui regroupe au plus $N^*$ classes $E_{j'}^k$ telles que les $B_{j'}^k$ sont de diam\`etre major\'e par $6\lambda\left(6\lambda N^*\right)^kr_{i_0}$. Par cons\'equent, la r\'eunion des $B^k_{j'}$, et donc a fortiori $\{x_i, i\in E^{k+1}_j\}$, est de diam\`etre $\delta_{1,j}^{k+1}$ major\'e par
$$6\lambda N^*\left(6\lambda N^*\right)^kr_{i_0}= \left(6\lambda N^*\right)^{k+1}r_{i_0}.$$
On conclut en remarquant que $r_{i_0}= \delta^{k+1}_{2,j}$, par d\'efinition.
\end{demo}

\begin{lemm}\label{l.boules4}
Supposons que pour un entier $k$ il existe une classe $E^k_j$ qui n'appartienne pas \`a la partition $\cE^{k-1}$. Alors le cardinal de $E^k_j$ est strictement sup\'erieur \`a $k$. 
\end{lemm}
\begin{demo} La classe $E^k_j$ est partitionn\'ee par les classes de $\cE^{k-1}$ qui sont adjacentes entre elles, par d\'efinition. Si toutes ces classes \'etaient d\'ej\`a des classes de $\cE^{k-2}$, la d\'efinition de l'adjacence montrerait que  tous les points de $E^k_j$ seraient $k-1$-adjacents, et donc $E^k_j$ appartiendrait \`a $\cE^{k-1}$. Ceci montre que l'une des classes de $\cE^{k-1}$ incluse dans $E^k_j$ n'appartient pas \`a $\cE^{k-2}$.

L'affirmations se montre \`a pr\'esent par une r\'ecurrence facile sur $k$. 
\end{demo}

\begin{lemm} Les partitions $\cE^{N^*}$ et $\cE^{N^*+1}$ sont \'egales.
\label{l.boules5}
\end{lemm}
\begin{demo} 
Supposons, par l'absurde, que ces partitions sont diff\'erentes. D'apr\`es le lemme~\ref{l.boules4}, il existe alors dans $\cE^{N^*+1}$ une classe de cardinal strictement sup\'erieur \`a $N^*$, ce qui contredit le lemme~\ref{l.boules3}.
\end{demo}

\subsection{Preuve du lemme~\ref{l.suite2}}\label{s.suite2}

 Afin de montrer le lemme~\ref{l.suite2} on consid\`ere une suite  $$a^0_0=a^0,\dots,a^0_t,\dots,a^0_N,a^1_0,\dots,a^1_N,\dots,a^m_0,\dots,a^m_N=b^m$$  
v\'erifiant~:
\begin{itemize}
\item  pour tout $n\in \{0,\dots,m\}$, les points $a^n_0,\dots,a^n_N$ appartiennent \`a un m\^eme cube $(1+\varepsilon)C_n$, les cubes $C_n$ \'etant deux \`a deux distincts, 
\item pour tout $n\in \{0,\dots,m\}$ et $t\in\{0,\dots, N-1\}$, la distance $d(f^t(a^n_{t}),f^t(a^n_{t+1}))$ est inf\'erieure \`a $\eta$ fois la distance entre  $f^t((1+\varepsilon) C_n)$ et le compl\'ementaire de $f^t((1+2\varepsilon)C_n)$.
\end{itemize}

Les points de la forme $f^N(a^n_N)$ ont un r\^ole particulier. De ce fait nous consid\'erons la suite des $f^t(a^n_t)$ pour $0\leq t\leq N-1$ ordonn\'ee de la fa\c con suivante~:
  $$a^0_0,\dots,f^t(a^0_t),\dots,f^{N-1}(a^0_{N-1}),a^1_0,\dots,f^{N-1}(a^1_{N-1}),\dots,a^m_0,\dots,f^{N-1}(a^m_{N-1}).$$

Afin de simplifier les notations, nous notons cette suite ordonn\'ee $\{x_i\}_{i\in I}$, index\'ee de fa\c con croissante, o\`u $I=\{1,\dots,(m+1)N\}$. 

Nous allons appliquer la proposition~\ref{p.boules} \`a l'ensemble des points $x_i$.  Pour cela, l'espace m\'etrique $(X,d)$ de cette proposition sera la vari\'et\'e $M$ munie de sa distance. Les constantes $\lambda$ et $N^*$ de cette proposition sont celles fix\'ees au d\'ebut de la section~\ref{s.connecting}. La constante $K$ est choisie \'egale \`a $\frac1\eta$, o\`u $\eta$ \`a \'et\'e fix\'ee au d\'ebut de la section~\ref{s.connecting}.

Pour tout point $x_i=f^t(a^n_t)$ on d\'efinit $r_i$ comme $d(f^t(a^n_{t}),f^t(a^n_{t+1}))$. 

Ces choix nous donnent la majoration $K>(6\lambda N^*)^{N^*+1}$. L'item (1) du lemme ci-dessous montre alors que les hypoth\`eses de la proposition~\ref{p.boules} sont v\'erifi\'ees.
\begin{lemm}\label{l.hypprop}Avec les notations ci-dessus on a~:
\begin{enumerate}
\item pour tout $i\in I$ , le cardinal de l'ensemble $\{i'\in I, B(x_i,Kr_i)\cap B(x_{i'},Kr_{i'})\neq \emptyset\}$ est major\'e par $N^*$~;
\item si $x_{i}=f^{t}(a^{n}_{t})$, la boule $B(x_i,Kr_i)$ est contenue dans l'image $f^t(\cC)$ du cube quadrill\'e.
\end{enumerate}
\end{lemm}
\begin{demo} Par d\'efinition de $r_i$ et de $K$,  si $x_{i}=f^{t}(a^{n}_{t})$, la boule $B(x_i,Kr_i)$ est contenue dans l'image $f^t((1+2\varepsilon)C_n)$ o\`u $C_n$ est le carreau associ\'e \`a $a^n_t$  (en particulier cette boule est incluse dans $f^t(\cC)$ ce qui montre l'item (2). 

Donc si $B(x_{i_1},Kr_{i_1})\cap B(x_{i_2},Kr_{i_2})\neq \emptyset$ alors les cubes $f^{t_1}((1+2\varepsilon)C_{n_1})$ et $f^{t_2}((1+2\varepsilon)C_{n_2})$ s'intersectent, ce qui implique que $t_1=t_2$ (les $N$ premiers it\'er\'es de $\cC$ sont disjoints) et, par le choix de $\varepsilon$, que les carreaux $C_{n_1}$ et $C_{n_2}$ sont adjacents.  Par d\'efinition de $N^*$, ceci montre l'item (1).
\end{demo}

La proposition~\ref{p.boules} nous donne donc une partition $\cE=\{E_j, j\in J\}$ de $I$ v\'erifiant~: 
\begin{itemize}
\item Pour tout $j\in J$, notons $\delta_{1,j}$ le diam\`etre de l'ensemble $\{x_i, i\in E_j\}$, et $\delta_{2,j}=\sup \{r_i, i\in E_j\}$.  Finalement notons $\delta_j=\delta_{1,j}+\delta_{2,j}$. Soit $B_j$ l'union des boules $B(x_i, \lambda \delta_j)$, pour $i\in E_j$. Alors les $B_j$ sont deux \`a deux disjoints.

\item Pour tout $j\in J$, il existe $i_0\in E_j$ tel que l'ensemble $B_j$ est contenu dans la boule $B(x_{i_0},Kr_{i_0})$.
\end{itemize}

L'item (2) du lemme~\ref{l.hypprop} implique donc que, pour tout $j\in J$, il existe $t$ tel que  $B_j\subset f^t(\cC)$. 
Notons $P_j$ l'ensemble des points $x_i$ avec $i\in E_j$. En particulier $P_j\subset f^t(\cC)$ et tous les points $x_i\in P_j$ sont de la forme $f^t(a^n_t)$. 

Notons $\cP$ l'ensemble $\cP=\{P_j, j\in J\}$~; c'est une partition de $\{x_i,i\in I\}$. 
Pour tout $t\in\{0,\dots,N-1\}$, la famille des $P_j\in\cP $ inclus dans $f^t(\cC)$ induit donc une partition $\cP_t$ de l'ensemble 
$\{a_t^n, n\in\{0,\dots,m\}\}$ (en effet cet ensemble est $f^{-t}\left(f^t(\cC)\cap\{x_i,i\in I\}\right)$). Ainsi, si $P$ est un \'el\'ement de la partition $\cP_t$, alors $f^t(P)$ est un \'el\'ement de $\cP$.

On d\'efinit par r\'ecurrence une suite $\nu(0),\nu(1),\dots,\nu(s(N+1)+N)$ d'entier de $\{0,\dots,m\}$
de la fa\c con suivante~:  on pose $\nu(0)=0$. On suppose \`a pr\'esent que la suite est construite jusqu'\`a un nombre de la forme $k(N+1)+t$, avec $t\in\{0,\dots, N\}$.

\begin{enumerate}
\item\label{i.t<N} Si $t<N$, on consid\`ere le point $a_t^{\nu(k(N+1)+t)}$. Il appartient \`a un unique \'el\'ement $P$ de la partition $\cP_t$. Alors $\nu(k(N+1)+t+1)$ est le plus grand entier $r\in\{0,\dots,m\}$ tel que $a_t^r$ appartient \`a $P$. 
 
\item\label{i.t=N} Si $t=N$, et si $\nu(k(N+1)+N)$ est strictement inf\'erieur \`a $m$, on pose  $\nu((k+1)(N+1))=\nu(k(N+1)+N)+1$~; c'est bien un entier de $\{0,\dots,m\}$. 
\item Si $t=N$ et $\nu(k(N+1)+N)=m$, la suite s'arr\^ete et on pose $s=k$.
\end{enumerate}

Remarquons que la suite $\nu$ est croissante (non n\'ecessairement strictement). De plus la sous-suite $\nu(k(N+1))$ est strictement croissante ce qui assure que cette suite est finie~: $s$ est inf\'erieur \`a $m$.   

Le lemme suivant conclut la preuve du lemme~\ref{l.suite2}.
\begin{lemm}La suite $\nu(0),\dots,\nu(s(N+1)+N)$ satisfait les conclusions du lemme~\ref{l.suite2}.
\end{lemm}
\begin{demo} Nous avons vu par construction que la suite est croissante. De plus l'item~(\ref{i.t=N}) de la construction de $\nu$ donne l'item (1) du lemme~\ref{l.suite2}~: pour tout $k\in\{1,\dots,s\}$, on a $\nu(k(N+1))=\nu(k(N+1)-1)+1$. 

Montrons d'abord que, pour tout $t\in\{0,\dots,N-1\}$, tout ensemble $P$ de la partition $\cP_t$ contient au plus un point de la forme $a_t^{\nu(k(N+1)+t)}$~: en effet, soit $k\in\{0,\dots,s\}$ le plus petit entier tel que $a_t^{\nu(k(N+1)+t)}\in P$. D'apr\`es l'item~(\ref{i.t<N}) de la construction de $\nu$, l'indice  $\nu(k(N+1)+t+1)$ est le plus grand entier $r$ tel que $a_t^r$ appartient \`a $P$.  La suite $\nu$ \'etant croissante, $\nu(k(N+1)+N)$ est sup\'erieur ou \'egal \`a $\nu(k(N+1)+t+1)$. L'item~(\ref{i.t=N}) de la construction de $\nu$ (et la croissance de $\nu$) montre $\nu((k+1)(N+1)+t)>\nu(k(N+1)+t+1)$. La suite $\nu$ \'etant croissante, pour tout $k'>k$ le point $a_t^{\nu(k'(N+1)+t)}$ n'appartient pas \`a $P$, ce qui montre notre affirmation. 

En cons\'equence, tout ensemble $P_j\in\cP$ contient au plus un point de la forme $f^t(a^{\nu(k(N+1)+t)})$, avec $k\in\{0,\dots,s\}$ et  $t\in\{0,\dots,N-1\}$.

Pour tout $k\in\{0,\dots,s\}$ et tout $t\in\{0,\dots,N-1\}$,  notons $B^k_t$ la boule centr\'ee au point $f^t\left(a^{\nu(k(N+1)+t)}_t\right)$ et  de rayon $\lambda d\left(f^t\left(a^{\nu(k(N+1)+t)}_t\right),f^t\left(a^{\nu(k(N+1)+t+1)}_{t+1}\right)\right)$. Afin de conclure la preuve du lemme~\ref{l.suite2}, il reste \`a montrer les items (2) et (3) c'est-\`a-dire que les boules $B_t^k$ sont deux \`a deux disjointes et que, pour tout $t\in\{0,\dots,N-1\}$, chaque boule $B_t^k$ est contenue dans $f^t(\cC)$. 

Soient $k\in\{0,\dots,s\}$ et $t\in\{0,\dots,N-1\}$.  Le point $f^t\left(a^{\nu(k(N+1)+t)}_t\right)$ appartient \`a un ensemble $P_j\in\cP$ associ\'e \`a une classe $E_j\in \cE$. Nous allons montrer~:
\begin{affi}\label{a.inclusion} Avec ces notations, la boule $B_t^k$ est incluse dans l'ensemble $B_j$ associ\'e \`a la classe $E_j$.
\end{affi}
Voyons d'abord que cette affirmation conclut la preuve du lemme~\ref{l.suite2}. En effet, nous avons vu que l'ensemble $B_j$ est contenu dans $f^t(\cC)$~; il en est donc de m\^eme pour $B_t^k$ (ceci montre l'item (3)). 

Pour montrer l'item (2), rappelons que les $B_j$ sont deux \`a deux disjoints, et que $B_j$ contient l'ensemble $P_j$. L'ensemble $B_j$ est donc disjoint de $P_{j'}$ pour $j'\neq j$. Comme chaque $P_j$ contient au plus un point de la forme $f^t(a^{\nu(k(N+1)+t)})$, il en est de m\^eme pour les ensembles $B_j$.  Utilisant l'affirmation~\ref{a.inclusion} et \`a nouveau que les $B_j$ sont deux \`a deux disjoints, un ensemble $B_j$ rencontre au plus un ensemble $B_t^k$. Les ensembles $B_t^k$ sont donc deux \`a deux disjoints. 

\begin{demo}[Preuve de l'affirmation~\ref{a.inclusion}] 
Estimons le rayon de la boule $B_t^k$. Pour cela on \'ecrit~:
\begin{equation}\label{e.rayon}
\begin{split}
d&\left(
f^t\left(a^{\nu(k(N+1)+t)}_t\right),f^t\left(a^{\nu(k(N+1)+t+1)}_{t+1}\right) \right)\leq\\
&\leq d\left(
f^t\left(a^{\nu(k(N+1)+t)}_t\right),f^t\left(a^{\nu(k(N+1)+t+1)}_{t}\right) \right)+
d\left(
f^t\left(a^{\nu(k(N+1)+t+1)}_t\right),f^t\left(a^{\nu(k(N+1)+t+1)}_{t+1}\right) \right).
\end{split}
\end{equation}

Le premier terme de la somme est une distance entre deux points de $P_j$~;  en effet, le point $f^t\left(a^{\nu(k(N+1)+t)}_t\right)$ appartient \`a $P_j$ par d\'efinition de $j$ et $f^t\left(a^{\nu(k(N+1)+t+1)}_{t}\right)$ appartient \`a $P_j$ par d\'efinition de $\nu(k(N+1)+t+1)$, voir l'item~(\ref{i.t<N}) de la construction de $\nu$. Ce premier terme est donc inf\'erieur au diam\`etre $\delta_{1,j}$ de $P_j$.

Le point $f^t\left(a^{\nu(k(N+1)+t+1)}_{t}\right)$ est un point $x_{i_0}$, $i_0\in I$, de la suite consid\'er\'ee au d\'ebut de cette section. Le second terme de la somme dans (\ref{e.rayon}) est le rayon $r_{i_0}$ correspondant \`a  $x_{i_0}$. Puisque $x_{i_0}$ appartient \`a $P_j$, ce second terme est donc inf\'erieur \`a $\delta_{2,j}$ (qui est le plus grand des rayons associ\'es aux points de $P_j$). 

On en d\'eduit~:
$$
\lambda d\left(
f^t\left(a^{\nu(k(N+1)+t)}_t\right),f^t\left(a^{\nu(k(N+1)+t+1)}_{t+1}\right) \right)\leq \lambda(\delta_{1,j}+\delta_{2,j})\leq \lambda\delta_j.
$$
Puisque $B_j$ est la r\'eunion des boules $B(x_i,\lambda\delta_j)$, pour $i\in E_j$, et que $f^t\left(a^{\nu(k(N+1)+t)}_t\right)$ peut s'\'ecrire de la forme $x_i$, $i\in E_j$, on obtient que $B_t^k$ est contenu dans $B_j$.

\end{demo}

Ceci termine la preuve du lemme~\ref{l.suite2}, et compl\`ete donc la preuve du th\'eor\`eme~\ref{t.connecting}.
\end{demo}

\vskip 1cm
\flushleft
{\bf Christian Bonatti} \ \  (bonatti@u-bourgogne.fr)\\
Institut de Math\'ematiques de Bourgogne,  UMR 5584 du CNRS\\
BP 47 870\\
21078 Dijon Cedex \\

France\\

\medskip
\flushleft{\bf Sylvain Crovisier} \ \  (sylvain.crovisier@u-bourgogne.fr)\\
Institut de Math\'ematiques de Bourgogne,  UMR 5584 du CNRS\\
BP 47 870\\
21078 Dijon Cedex \\
France\\
\end{document}